\documentstyle[epsf]{article}
\newcommand{\be}{\begin{equation}}
\newcommand{\ee}{\end{equation}}
\newcommand{\bea}{\begin{eqnarray}}
\newcommand{\eea}{\end{eqnarray}}
\newcommand{\beas}{\begin{eqnarray*}}
\newcommand{\eeas}{\end{eqnarray*}}
\newcommand{\bdm}{\begin{displaymath}}
\newcommand{\edm}{\end{displaymath}}
\newcommand{\ba}{\begin{array}}
\newcommand{\ea}{\end{array}}
\newcommand{\bi}{\begin{itemize}}
\newcommand{\ei}{\end{itemize}}
\newcommand{\ben}{\begin{enumerate}}
\newcommand{\een}{\end{enumerate}}
\newcommand{\bc}{\begin{center}}
\newcommand{\ec}{\end{center}}
\newcommand{\bfl}{\begin{flushleft}}
\newcommand{\efl}{\end{flushleft}}
\newcommand{\bfr}{\begin{flushright}}
\newcommand{\efr}{\end{flushright}}
\newcommand{\bd}{\begin{description}}
\newcommand{\ed}{\end{description}}
\newcommand{\bq}{\begin{quote}}
\newcommand{\eq}{\end{quote}}
\newcommand{\bfg}{\begin{figure}}
\newcommand{\efg}{\end{figure}}
\newcommand{\bt}{\begin{table}}
\newcommand{\et}{\end{table}}
\newcommand{\btb}{\begin{tabular}}
\newcommand{\etb}{\end{tabular}}
\newcommand{\btg}{\begin{tabbing}}
\newcommand{\etg}{\end{tabbing}}

\newcommand{\kslash}
           {\mbox{$ k \hspace{-1.1ex} \mbox{/} \hspace{-0.07ex} $}}
\newcommand{\lslash}
           {\mbox{$ l \hspace{-0.9ex} \mbox{/} \hspace{-0.15ex} $}}
\newcommand{\pslash}
           {\mbox{$ p \hspace{-1ex} \mbox{/} \hspace{-0.08ex} $}}
\newcommand{\qslash}
           {\mbox{$ q \hspace{-1.1ex} \mbox{/} \hspace{-0.05ex} $}}

\newcommand{\wrt}{\mbox{\epsfysize=4mm \epsfbox{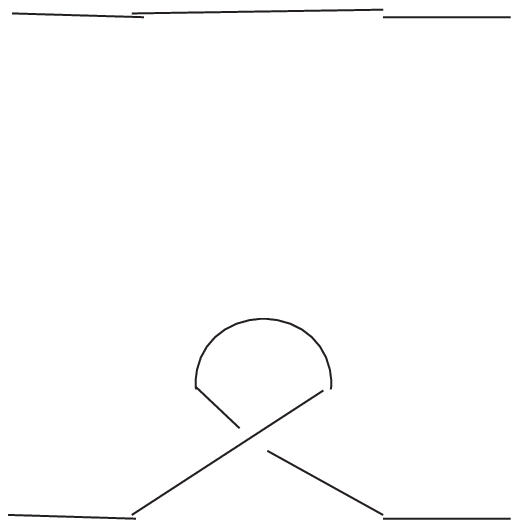}\hspace{-0.2em}}}
\newcommand{\cil}{\mbox{\epsfysize=4mm \epsfbox{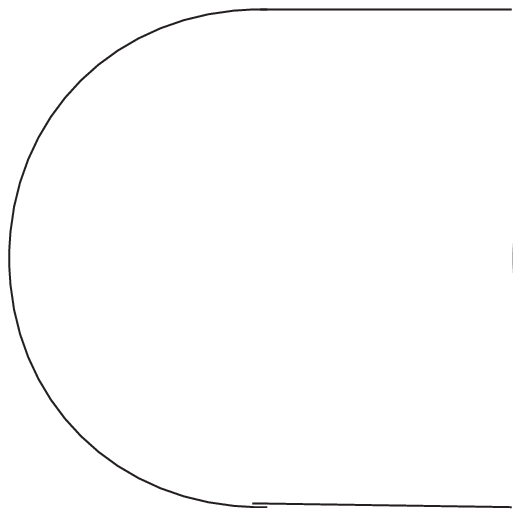}\hspace{-0.2em}}}
\newcommand{\cir}{\mbox{\epsfysize=4mm \epsfbox{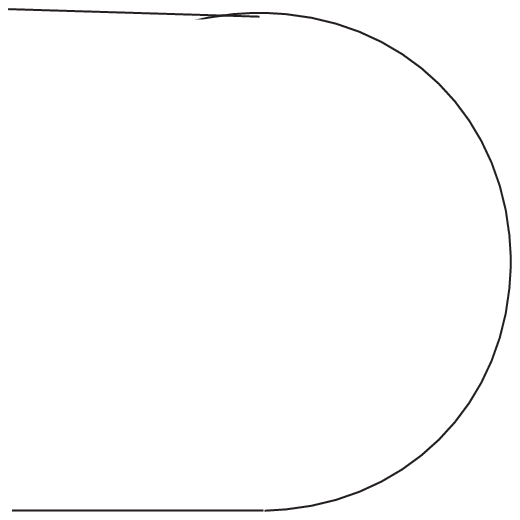}\hspace{-0.2em}}}
\newcommand{\crs}{\mbox{\epsfysize=4mm \epsfbox{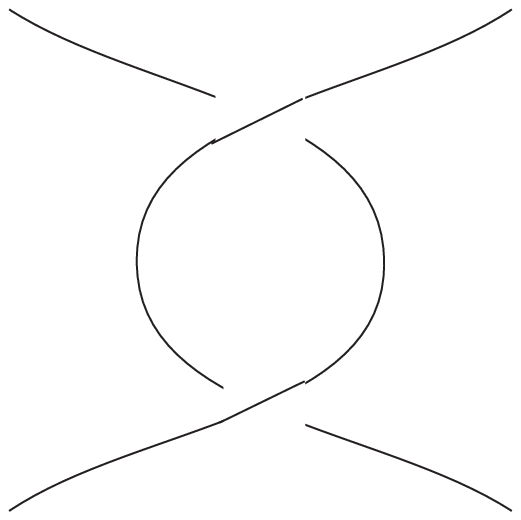}\hspace{-0.2em}}}
\newcommand{\cic}{\mbox{\epsfysize=4mm \epsfbox{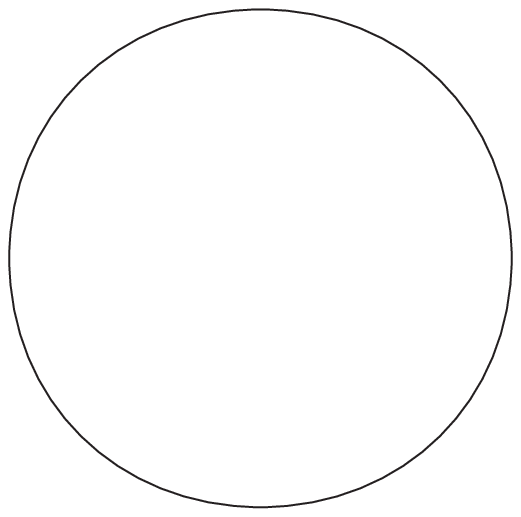}\hspace{-0.2em}}}
\newcommand{\sip}{\mbox{\epsfysize=6mm \epsfbox{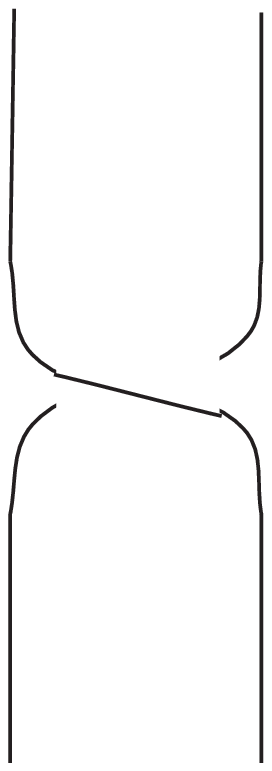}\hspace{-0.2em}}}

\newcommand{\lil}{\mbox{\epsfysize=6mm \epsfbox{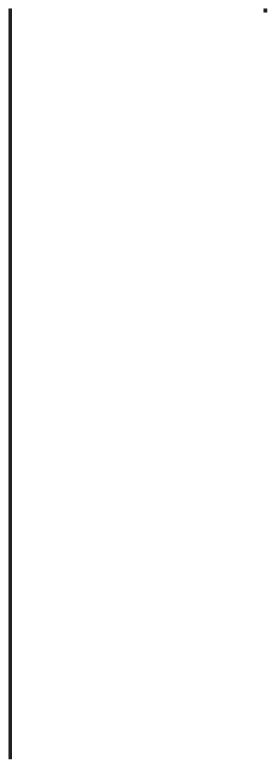}\hspace{-0.2em}}}
\newcommand{\lir}{\mbox{\epsfysize=6mm \epsfbox{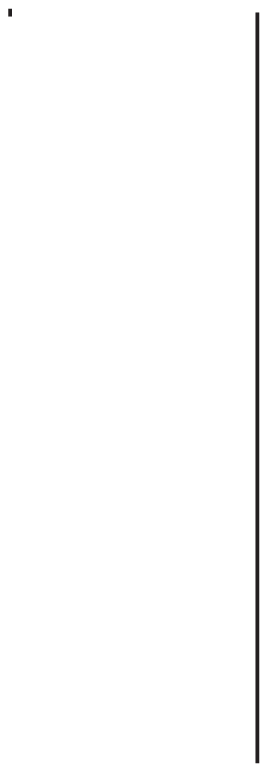}\hspace{-0.2em}}}

\begin{document}
\title{\raisebox{2cm}{\makebox[0pt][r]{\small MZ-TH//96-18}} 
{\Large Renormalization and Knot Theory}\footnote{Habilitationsschrift, 
to appear in {\em Journal of Knot Theory and its Ramifications.}}}
\author{Dirk Kreimer\thanks{email: kreimer@dipmza.physik.uni-mainz.de}\\
Dept.~of Physics\\
Mainz Univ. \hfill and \hfill Tasmania Univ.\\
Staudingerweg\hfill GPO Box 252C\\
55099 Mainz\hfill Hobart TAS 7001\\
Germany\hfill Australia
}
\maketitle
\begin{abstract}
We investigate to what extent renormalization can be understood as an
algebraic manipulation on concatenated one-loop integrals. We find that
the resulting algebra indicates a useful connection to knot theory
as well as number theory
and report on recent results in support of this connection.
\end{abstract}
\tableofcontents
\section{Introduction}
Renormalization theory is one of the fundamental topics of perturbative 
Quantum Field Theory (pQFT). Renormalizability was one of the guiding
principles in the construction of the Standard Model (SM), our current
understanding of particle physics phenomenology. From first insights
into the problem of ultraviolet (UV) divergences in a QFT 
\cite{schweber} to the understanding of renormalization of gauge theories
\cite{thooft}, including spontaneous symmetry breaking
and the beautiful interplay with BRST identities \cite{Lee}, our 
understanding has vastly improved. 


The first complete account to renormalization as such
has been given by Dyson, Salam and
Weinberg \cite{Dy,Salam,Weinberg}. One can derive all properties of a renormalizable
QFT from the knowledge of the Dyson Schwinger equations, an
understanding of the problem of overlapping divergences
\cite{Salam} and Weinberg's theorem \cite{Weinberg}.
Later, the subject was reinvestigated thoroughly by 
the authors of \cite{BPHZ}.
Their approach was helpful especially in understanding parametric
representations, but often the structure of Green functions is less
transparent  there. Nevertheless, the BPHZ method has become the standard text 
book  approach to renormalization.
Few would argue with the view that
the structure of a renormalizable pQFT is most successfully encoded 
in structures like the forest formula \cite{Zimm}, the BPHZ formalism
 and the $R$-operation \cite{BPHZ}.


Yet, as we will show, there seems to be still a more fundamental layer
behind these structures. This paper is concerned with the exploration
of an algebraic structure inherent in all renormalizable quantum field 
theories, which reveals a connection to knot theory and similar algebraic
structures usually assigned to the investigation of braids. We will focus
on this connection. We hope to show that these structures are not only
of interest in their own right but also allow for a useful simplification of
actual computations of renormalized quantities.


This paper is based on the material in \cite{org}.
We added an introduction into the subject, refined and extended
the arguments 
at various places and added a section reviewing the results obtained 
recently.

The paper is organized as follows. In section two we will 
introduce ideas and notions of perturbative
quantum field theory and renormalization theory.
This section is by no means a complete account, but should provide 
for the reader unfamiliar with
renormalization the necessary background.
We refer to \cite{schweber,Kaku,Collins} for a more detailed 
discussion of the subject. This section can be skipped by a reader
who is familiar with the topics of Feynman diagrams and renormalization
theory. 

In section three, 
we will restrict ourselves
to vertex corrections of ladder-type diagrams. This simple topology
mainly serves to fix our notations and to exhibit the basic idea.
We will introduce a notation which allows us to express results
for the simplest topology, the ladder topology, in terms
of one-loop functions. The results obtained are not new as
such, but allow to prove a remarkable fact: these simple
topologies add under renormalization to Z-factors free of
transcendental numbers. Sections three to six are mainly
concerned with establishing this one-loop algebra in
all relevant situations. Section three considers the ladder topology
with nested divergences, section four establishes the same
for two-point functions, section five considers the case of
iterated nested and disjoint topologies. Finally, section six  addresses the
study of overlapping ladder topologies, remarkably establishing
the same one-loop algebra in this case as well. 
In section three and six
we shortly disgress to  link diagrams, exploiting the fact that
the one-loop algebra obtained so far fits into a pattern
motivated by skein relations \cite{knots} and braid algebras.
The core of this paper is section seven, where we present an argument that
any ladder topology, realized by what renormalizable theory and
Feynman graph whatsoever, is free of trans\-cendental numbers for the
$Z$ factors: its UV-divergences have purely rational coefficients.



In fact, we will show that knot theory knows about the renormalization 
program in the sense that it generates all contributing forests with 
the help of the skein relation, revealing the proposed structure
behind the recursion which governs renormalization theory.


This opens the arena for  a link between the topology of Feynman
graphs,
expressed in terms of knots associated with the graph,
and transcendental coefficients appearing in its $Z$-factor
contributions. So, in section eight, we are concerned with collecting
empirical evidence for the appearance of transcendentals in accordance
with knots
identified in the graphs. We review the results \cite{plb,pisa,bdk,bgk,
db,newplb} which
have been achieved so far. 
These results point towards an interesting conjecture.
Positive knots obtained from Feynman graphs correspond to Euler sums.
Indeed, it seems that we can identify as many distinct
prime knots in a Feynman graph as we find independent transcendentals
in it. These transcendentals appear as generalized $\zeta$-functions -
Euler sums, and we collect the evidence for this unsuspected 
connection between knot theory, number theory and field theory in
section eight.


Conclusions and an outlook finish the paper.


We will restrict ourselves to MS schemes in the following and
we will use dimensional regularization throughout the paper.

\section{Perturbative Quantum Field Theory}

In this section we introduce perturbative Quantum Field Theory (pQFT),
the idea of renormalization and discuss the necessary technical
notions. None of the results derived in following
sections is a new result in renormalization theory or uses more than
the most basic facts of it. But looking at the structure
of renormalization from an new angle we will argue
that one can establish a useful connection to knot theory.
We are not yet able to provide a field-theoretic derivation of
this connection but rather collect empirical evidence
by investigating  and exploring well known results of perturbation theory.
 This is not the place to give a full account of renormalization theory, or
to investigate the origin of UV divergences from the viewpoint
of axiomatic field theory. Rather, we will sketch how pQFT appears
as the art of calculating Feynman diagrams, explain the idea of
renormalization and introduce the vocabulary.


\subsection{pQFT}
We now want to introduce pQFT and Feynman diagrams. Our presentation
follows \cite{Kaku}. Our aim is to motivate the derivation of
the Feynman rules. These rules give a pictorial approach to pQFT.
Free propagating particles are specified by lines (propagators)
and do interact at vertices, points in space-time where three or more
lines merge together. These points describe the coupling of
particles, and are specified by
a coupling constant determing the strength of the interaction.
In this way, each graph built out of propagators
and vertices describes a process in pQFT. 
Specifying a set of initial and final particles,
in pQFT one considers all possible graphs which allow
for these initial and final states. If there are no closed
loops in the graph, we speak of a tree-level graph. 
It delivers the lowest order contribution to the specified
process. These contributions
are of no further interest to us.
Higher order contributions in the coupling
constant demand
more interaction vertices.
The specified set of external particles is kept
fixed and thus the expansion in the coupling
constant becomes a loop expansion.  In the presence of closed loops
in the graph one has to take into account all possible
momenta for these internal particles.  This demands
to integrate internal loop momenta. 
For the most interesting theories, these integrations
are ill-defined. They diverge for large internal momenta.
Quite often one can handle these divergences in an iterative
procedure called renormalization. There, we claim,
a connection to knot theory appears: the singular contributions
to a pQFT are determined by knots derived from the topologies
of its Feynman graphs. 


So let us now motivate these Feynman rules.
Let our starting point be the consideration of a scalar
field theory; there is only
one self-interacting scalar particle present. 
More realistic examples involve particles sitting
in more complicated representations of the Lorentz group
(spin), but we omit these complication in
this elementary exposition.


We have a free Lagrangian given by
\begin{equation}
L=\frac{1}{2} (\partial_\mu \phi)^2 -\frac{1}{2}
m^2 \phi^2.
\end{equation}
Later we will add interaction terms of the form
$g\phi^n, n\geq 3$, to this Lagrangian. But first we  consider
quantization of the free theory.
The classical calculus of variation 
gives as the equation of motion
\begin{equation}
\partial_\mu\frac{\delta L}{\delta \partial_\mu \phi} 
-\frac{\delta L}{\delta \phi}=0 \Rightarrow\;(\Box+m^2) \phi=0.
\end{equation}
With the field $\phi(x)$ acting as the coordinate, parametrized
by the space-time point $x=({\bf x},t)$, we find as
its canonical conjugate momentum
\begin{equation}
\pi({\bf x},t)\equiv\frac{\delta L}{\delta \dot{\phi}}=\dot{\phi}({\bf x},t),
\end{equation}
and, accordingly, a Hamiltonian density
\begin{equation}
H=\frac{1}{2}[\pi^2+(\nabla\phi)^2+m^2\phi^2].
\end{equation}
We recall that the transition to quantum theory emerges when
we impose non-vanishing commutators between coordinates
and momenta. We consider the commutator at equal time.
We demand that commutators between expressions
separated by space-like distances vanish.\footnote{
This incorporates a notion of causality. Fields with
spacelike separation do not influence each other.} 
At equal time, this becomes
\begin{equation}
[\phi({\bf x},t),\pi({\bf x},t)]=i\delta^3({\bf x}-{\bf y}).
\end{equation}
Note that this separation is Lorentz invariant.
Lorentz transformations map spacelike distances to spacelike
distances.

\begin{figure}[ht]\epsfysize=2cm \epsfbox{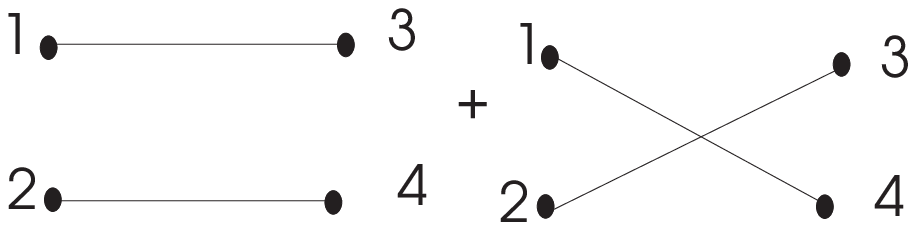} 
\caption[Free Particles.]{\label{hab21}\small It is easy to see that states like $<\!k_1,k_2|k_3,k_4\!>$
can be graphically represented as above. Their normalization
demands either $k_1=k_3,\;k_2=k_4$ or $k_1=k_4,\;k_2=k_3$.
In general, multiparticle states $<\!k_1,\ldots,k_n|q_1,\ldots,
q_m\!>$ do vanish for $n\not= m$, and for  $n=m$
their normalization delivers all pairings of the $k_i$
with the $q_j$. There is no interaction between
the straight lines corresponding to these pairings, which is indicated
by straight lines in the figure. The straight line represents an
incoming particle with momentum $k_i$ propagating undisturbed to
become an outgoing particle with momentum $q_j=k_i$.}
\end{figure}

We are now looking for fields $\phi(x)$ fulfilling
the Klein Gordon equation and the commutation relations.
One readily verifies that the {\em Ansatz}
\begin{equation}
\phi(x)=
\frac{1}{(2\pi)^{3/2}}\int d^4k \delta(k^2-m^2)
\Theta(k_0)[A(k)e^{-ikx}+A^\dagger(k) e^{ikx}],
\end{equation}
has the required properties.
We use
\begin{eqnarray}
\omega_k: & = & 
\sqrt{{\bf k}^2+m^2},\;e^{ikx}=\exp{(i\omega_k t-{\bf k\cdot
    x})},\nonumber\\
\int d^4k \delta(k^2-m^2)\Theta(k_0)
 & = & \int \frac{d^3k}{2\omega_k},\nonumber\\ 
a(k):=\frac{A(k)}{\sqrt{2\omega_k}}.
\end{eqnarray}
Upon Fourier transformation, we obtain expressions for $\phi,\pi$
where the operator valued Fourier coefficients $a,a^\dagger$
have to obey commutation relations
as well:
\begin{equation}
[a(k),a^\dagger(k^\prime)]=\delta^3({\bf k}-{\bf k}^\prime).
\end{equation}
Then we obtain the Hamiltonian as
\begin{equation}
H = \int d^3k\; \omega_k\; \left[a^\dagger(k) a(k) + \frac{1}{2}\right]=
\frac{1}{2}\int d^3k\; \omega_k\; \left[a^\dagger(k) a(k) 
+ a(k)a^\dagger(k)\right].
\end{equation}
The Hamiltonian is the generator of the time evolution,\footnote{$\phi(x)$ 
transforms as $U\phi(x)U^{-1}=\phi(x^\prime)$ which defines
it to be a scalar field.}
and accordingly we find the generator of spacelike
translations by replacing $\omega_k$ by the spacelike
components of the four vector $k_\mu=(\omega_k,{\bf k})$
\begin{equation}
P_i= \int d^3k\; k_i\; \left[a^\dagger(k) a(k) + a(k)a^\dagger(k)\right].
\end{equation}
By comparison with quantum mechanics
we recognize  the Hamiltonian as a sum (integral) over
harmonic oscillators, one for each mode $k$.
Accordingly, we assume the existence of a ground state $|0\!>$
annihilated by $a(k)$ for all $k$:
\begin{equation}
a(k)|0\!>=0.
\end{equation}
We also want to define normalized states
\begin{equation}
a^\dagger(k)|0\!>=|k\!>.
\end{equation}
But the energy of this ground state is ill defined
\begin{equation}
\int^\infty dk \;\omega_k=\infty.
\end{equation}
We discard the unobservable ground state energy by introducing
a normal ordering of creation operators
$a^\dagger(k)$ and annihilators
$a(k)$. We write all creation operators to the left,
and all annihilators to the right:
\begin{eqnarray}
:\!a(k)a^\dagger(k^\prime)\!:\;\;=0,\;
:\!a^\dagger(k) a(k^\prime)\!:\;\;=a^\dagger(k) a(k^\prime),\label{no}\\
\Rightarrow H=\int d^3k\; \omega_k\, :\!a^\dagger(k)a(k)\!:\;.\nonumber
\end{eqnarray}
Next we introduce states $|k_1,\ldots,k_n\!>$ as
\begin{equation}
|k_1,\ldots,k_n\!>=a^\dagger(k_1),\ldots,a^\dagger(k_n)|0\!>,\;
<\!k_1,\ldots,k_n|=<\!0|a(k_1),\ldots,a(k_n),
\end{equation}
with normalization 
\begin{equation}
<\!k^\prime|k\!>\;=\;<\!0|a(k^\prime)a^\dagger(k)|0\!>\;=\;
<\!0|[a(k^\prime)a^\dagger(k)-a^\dagger(k^\prime)a(k)]|0\!>\;=\;
\delta({\bf k}-{\bf k}^\prime).
\end{equation}
We interpret the states above as multiparticle states for particles
with momenta $k_i$.
There are no interactions present, cf.~Fig.(\ref{hab21}). 
By now we  have achieved the quantization of  the free Klein Gordon field.


To proceed,
we introduce a source term in the Klein Gordon equation:
\begin{equation}
(\Box+m^2)\phi=J(x).
\end{equation}
We look for a Green-function (the propagator) $\Delta_F(x-y)$ solution of
\begin{equation}
(\Box+m^2)\Delta_F(x-y)=-\delta^4(x-y).\label{dast}
\end{equation}
One easily sees that 
\begin{eqnarray}
\Delta_F(x-y) & = & \int \frac{d^4k}{(2\pi)^4}e^{-ik(x-y)}\Delta_F(k),
\nonumber\\
\Delta_F(k) & = & \frac{1}{k^2-m^2+i\eta},
\end{eqnarray}
is a solution, where the presence of the small imaginary part $i\eta$
specifies boundary conditions.\footnote{One can show that these
boundary conditions are in accord with causality: positive energy
solutions propagate forwards in time. Nevertheless, in
a relativistic covariant theory we also have negative energy solutions
propagating backwards in time, corresponding to antiparticles
propagating forwards in time.}
Using our Ansatz for the fields $\phi$ one finds that the propagator can be written with the help
of a time-ordered product $T$:
\begin{eqnarray}
\Delta_F(x-y) & = & <\!0|T[\phi(x)\phi(y)]|0\!>,\nonumber\\
T[\phi(x)\phi(y)] & := & \phi(x)\phi(y) \Leftrightarrow x_0>y_0,\nonumber\\
T[\phi(x)\phi(y)] & := & \phi(y)\phi(x) \Leftrightarrow y_0>x_0.
\end{eqnarray}
Having solved for the Green function of the Klein Gordon field
we are ready to introduce perturbation theory.
To get acquainted with it, we consider an example borrowed from
non-relativistic quantum mechanics.
We assume the existence of a Hamiltonian $H=H_0 + H_I$,
with a free part $H_0$ and an interaction term $H_I$.
The Schr\"odinger equation is
\begin{equation}
(i\partial_t-H)\psi=0\;\Rightarrow\;
(i\partial_t-H)G({\bf x},t;{\bf x}^\prime,t^\prime)=\delta^3({\bf x}-{\bf x}^\prime)
\delta(t-t^\prime),
\end{equation}
and we can use the Green function $G({\bf x},t;{\bf x}^\prime,t^\prime)$
to express the wave function
using Huygens principle:
\begin{equation}
\psi({\bf x},t)=\int d^3 x^\prime G({\bf x},t;{\bf x}^\prime,t^\prime)
\psi({\bf x}^\prime,t^\prime),
\;\;t>t^\prime.
\end{equation}

\begin{figure}[ht]\epsfysize=2cm \epsfbox{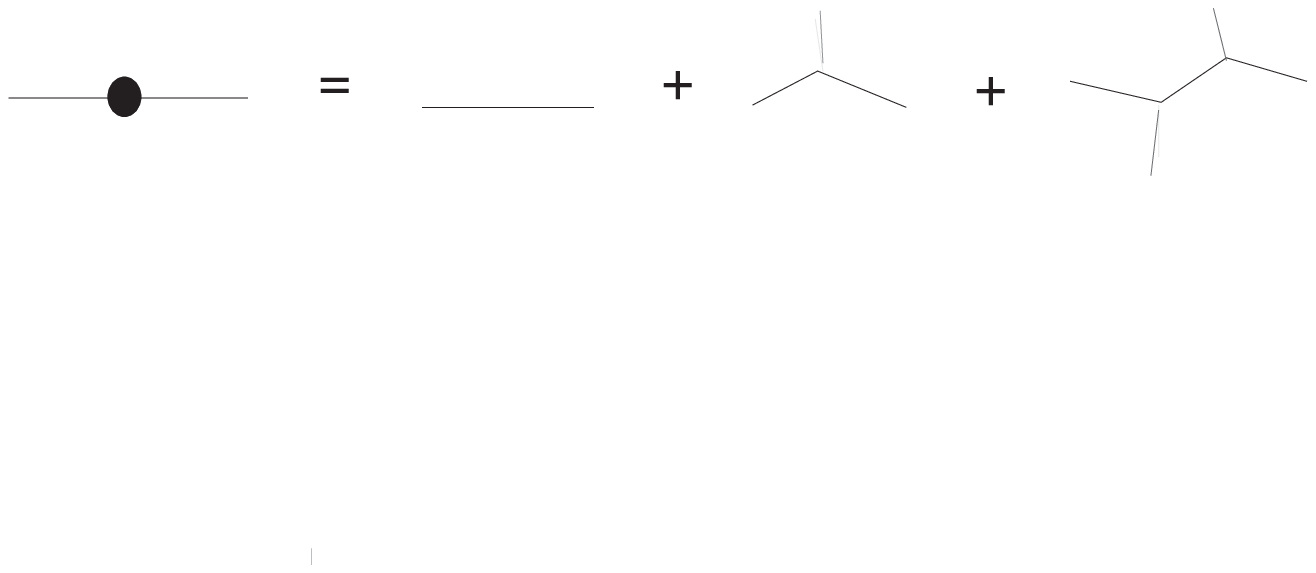} 
\caption[Expansion of the Green function.]{\label{hab22}\small 
The full Green function $G$ is a series starting with the bare
Green function $G_0$, denoted by a straight line in the above.
Interaction (for example with a classical background field $\cong$
dotted lines)
takes place at various points, which we describe by vertices.
We have to integrate over these points. Asymptotically,
we assume that we have incoming and outgoing plane waves.}
\end{figure}

Pondering the trivial identity
\begin{equation}
\frac{1}{A+B}=A^{-1}-A^{-1}BA^{-1} +A^{-1}BA^{-1}BA^{-1}+\ldots,
\end{equation}
and setting
\begin{eqnarray}
A=-H_0+i\partial_t & , & B=-H_I,\nonumber\\ 
G=\frac{1}{A+B} & , & G_0=\frac{1}{A},
\end{eqnarray}
we get a solution for the Green function in terms of iterated
interactions, graphically described in Fig.(\ref{hab22}):
\begin{eqnarray}
G({\bf x},t;{\bf y},s) & = & G_0({\bf x},t;{\bf y},s)\nonumber\\
 & &
+\int d^4x^\prime 
G_0({\bf x},t;{\bf x}^\prime,t^\prime)H({\bf x}^\prime,t^\prime)
G_0({\bf x}^\prime,t^\prime;{\bf y},s)+\ldots,\nonumber\\
G & = & G_0+G_0HG_0+\ldots\mbox{higher order in $H$}.
\end{eqnarray}


We can define the matrix element between an initial plane wave and a final plane wave:
\begin{equation}
S_{fi}:=\delta_{fi}+\int d^4x^\prime  \phi^\ast(x^\prime)H_I(x^\prime)
\phi(x^\prime)+\ldots,
\end{equation}
where we used that we can express the free Green function $G_0$
itself in terms
of plane waves. 


We now want to come back to relativistic quantum
field theory. We change from non-relativistic single-particle
quantum mechanics to relativistic multi-particle quantum field
theory. Instead of considering the interaction
with a classical source we will have to consider the full selfinteraction
of a quantized field.
Nevertheless, we will see that perturbation theory still allows
for a pictorial representation, and will now sketch
the remaining steps leading us to the celebrated Feynman rules.


We are interested in the $S$-matrix. It describes transitions
between incoming particles and final particles. This is the sitution
which we confront in particle physics: one has a set of
well defined incoming particles $|i\!>_{in}$ (prepared in some beam, say) 
assumed to be
non-interacting, one further has an interaction regime (a collision,
say, of two beams or a beam and a target) and detects a set
of outgoing particles $|f\!>_{out}$, again assumed to be non-interacting.


We define the $S$ matrix, which describes the probability of a transition
from an initial state $|i>_{in}$ to a final state $|f>_{out}$.
\begin{equation}
S_{f,i}:={}_{out}\!<\!f|i\!>_{in}.
\end{equation}
We assume that the asymptotic states ${}_{out }\!|f\!>$ and $|i\!>_{in}$ 
are related by a basis transformation:
\begin{equation}
|f\!>_{in}=R|f\!>_{out}, => S_{f,i}={}_{in}\!<\!f|R|i\!>_{in}
={}_{out}\!<\!f|R|i\!>_{out}.
\end{equation}
Asymptotic states are assumed to fulfil the free Klein Gordon
equation, and thus allow for a spectral decomposition similar
to a free field:
\begin{equation}
|q\!>_{in}=a^\dagger(q)|0\!>,
\end{equation}
where $a^\dagger(q)$ is the usual generator of a free one-particle
state with momentum $q$.
Our goal is to replace the fully interacting field $\phi$
everywhere by asymptotic fields,\footnote{Asymptotically, we assume
that the free fields and the full interacting fields
are connected in the weak sense: $\lim_{t\to \infty}<r|\phi|s>=Z^{\frac{1}{2}}
<r|\phi_{in}|s>$. At asymptotic times before any interaction takes place
the matrix elements of the full field approach the matrix elements
of the free field.}
as we did before in the quantum mechanics example.
 This is the so-called LSZ formalism \cite{LSZ}.
After a bit of algebra using the explicit representation
of $a^\dagger$ in terms of free fields
one finally obtains
\begin{eqnarray}
{}_{out}\!<\!p_1,\ldots,p_n|q_1,\ldots,q_m\!>_{in}=\nonumber\\
(iZ^{-\frac{1}{2}})^{n+m}\int d^4y_1\ldots d^4x_m \prod_{i=1}^n
\prod_{j=1}^m e^\ast_{p_i}(y_i) e_{q_j}(x_j)\nonumber\\
\times (\Box+m^2)_{y_1}\ldots (\Box+m^2)_{x_m}\nonumber\\
<\!0|T[\phi(y_1)\ldots \phi(x_m)]|0\!>,\label{sm}
\end{eqnarray}
where $e_q(x):=e^{iqx}/\sqrt{2\omega_q(2\pi)^3}$.
We have expressed the abstract S-matrix element as a vacuum expectation
value of fully interacting fields.
To find such expectation values we again refer to perturbation theory.
As an example, we consider an interaction Hamiltonian
\begin{equation}
H_I=\int d^3x\;{\cal H}_I=\int d^3x\;\frac{g}{4!}\phi^4.
\end{equation}
Assuming the existence of an unitary operator
$U(t)$ such that
\begin{equation}
\phi({\bf x},t)=U(t)\phi_{in}({\bf x},t)U^{-1}(t),
\end{equation}
one can show that $U(t)$ fulfils the requirements of a time evolution operator
\begin{eqnarray}
U(t):=U(t,-\infty),\; U(t_1,t_1)=1,\nonumber\\
U^{-1}(t_1,t_2)=U(t_2,t_1),\;
U(t_1,t_2)U(t_2,t_3)=U(t_1,t_3),
\end{eqnarray}
with boundary condition 
$\lim_{\to\infty}U(-t)|0\!>=|0\!>$.\footnote{The assumption that such an unitary operator connecting
the full field $\phi$ and the asymptotic field $\phi_{in}$ exists
cannot be maintained in a strict approach to QFT
and has to be relaxed somewhat. A good account
of the problem can be found in \cite{Haag}. Nevertheless, it allows a quick derivation of
perturbation theory. It is not understood why, then, perturbation theory works
so well.}
$U(t)$ being a time evolution operator, it fulfils a differential
equation:
\begin{equation}
\partial_t U(t)=-iH_I(t)U(t).
\end{equation}
Solving this equation as an exponential
one finally obtains a perturbative expansion for $U(t)$:
\begin{equation}
U(t)=1+\sum_{n=1}^\infty \frac{(-i)^n}{n!}\int_{-\infty}^t d^4x_1
\int_{-\infty}^t d^4x_2\ldots \int_{-\infty}^t d^4x_n
T[H_I(x_1)\ldots H_I(x_n)],
\end{equation}
which eventually serves to give the perturbation expansion for the
Green function 
\begin{equation}
G(x_1,\ldots,x_n):= <0|T[\phi(x_1)\ldots\phi(x_n)]|0>
\end{equation}
as
\begin{eqnarray}
& & G(x_1,x_2,\ldots,x_n) = \nonumber\\
& & \sum_{m=1}^\infty \frac{(-i)^m}{m!}
\int_{-\infty}^\infty d^4y_1
\ldots d^4y_m
<0|T[\phi_{in}(x_1)\ldots\phi_{in}(x_n)H_I(y_1)\ldots H_I(y_m)]|0\!>.
\end{eqnarray}
From these vacuum correlators one can obtain all
$S$-matrix elements.


We see that we obtain the Green function as a time ordered product
of free fields, with interactions $H_i$ inserted at various points.
We proceed by using the Wick theorem. This theorem
allows us to express the time-ordered product of fields as a
normal-ordered product (which has vanishing vacuum expectation
value by definition Eq.[\ref{no}]) and products of free propagators.
Using the spectral properties of free fields one
easily establishes
\begin{equation}
T[\phi_{in}(x_1)\phi_{in}(x_2)]=:\phi_{in}(x_1)\phi_{in}(x_2):+
<\!0|T[\phi_{in}(x_1)\phi_{in}(x_2)]|0\!>,
\end{equation}
and derives the vacuum expectation of free fields 
as 
\begin{eqnarray}
<\!0|T[\phi_{in}(x_1)\phi_{in}(x_2)\ldots \phi_{in}(x_n)]|0\!>
=\sum_{\mbox{perm}}
<0|T[\phi_{in}(x_1)\phi_{in}(x_2)]|0>\ldots\nonumber\\
<0|T[\phi_{in}(x_{n-1})\phi_{in}(x_n)]|0>.
\end{eqnarray}
Upon using
Eq.(\ref{dast})
one can discard the Klein-Gordon factors in Eq.(\ref{sm}) as they
become delta functions.

\begin{figure}[ht]\epsfxsize=5in \epsfbox{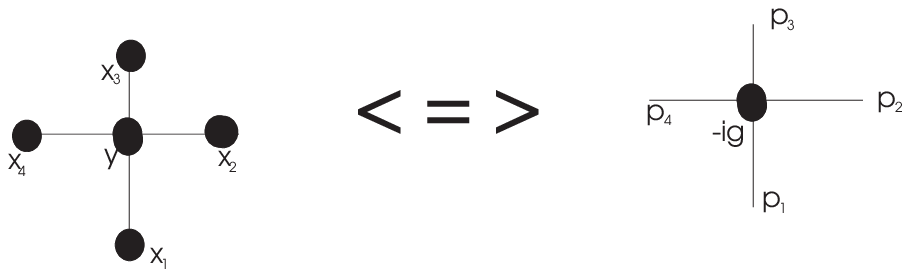} 
\caption[Interaction at a vertex point.]{\label{hab23}\small 
A graphical illustration of Eq.(\ref{gr}).
Upon Fourier-transformation we obtain overall momentum
conservation $p_1+p_2+p_3+p_4=0$, where $\Delta_F(x_i-y)=
\int d^4p_i e^{ip_i(x_i-y)}/(p_i^2-m^2+i\eta)$.}
\end{figure}

Finally we obtain to first order in $H_I$ for our example
of $H_I\sim \phi^4$
\begin{eqnarray}
G(x_1,x_2,x_3,x_4) & = & -\frac{ig}{4!}\int d^4y
<\!0|T[\phi_{in}(x_1)\ldots\phi_{in}(x_4)\phi_{in}(y)^4]|0\!>+\ldots
=\nonumber\\
 & & (-ig)\int d^4y \prod_{i=1}^4 (-\Delta_F(x_i-y))+\ldots,\label{gr}
\end{eqnarray}
which is described in Fig.(\ref{hab23}).
We are now prepared to introduce the Feynman rules of the theory.
We give them in momentum space. So we Fourier transform
all our Green functions and use the Fourier
transformed propagators $\Delta_F(q_i)$. We can now do the $y$ integral in
Eq.(\ref{gr}). It delivers a $\delta$-function
which guarantees momentum conservation: $\int d^4y\;
\exp{[i(p_1+p_2+p_3+p_4)y]}=
\delta^4(p_1+p_2+p_3+p_4)$. Note that we only Fourier transform
the arguments of the Green functions. Internal points
like $y$ are to be integrated. In momentum space this amounts to
an integration over internal loop momenta, which appear when we have
closed loops in the graph. 


So here come the Feynman rules for $\phi^4$ theory:
\begin{itemize}
\item For each internal scalar field line, associate a propagator given by:\\
$i\Delta_F(p)=\frac{i}{p^2-m^2+i\eta},$
\item at each vertex place a factor of:\\
$-ig,$
\item for each internal loop, integrate over:\\
$\int \frac{d^4q}{(2\pi)^4},$ where $q$ is a loop momentum.
\item Divide by the symmetry factor of the graph $G$, which is assumed to
contain $r$ vertices, say. 
It is given by the number of
possibilities to connect $r$ vertices to give $G$,
divided by $(4!)^r$.
\end{itemize}
The reader can easily check that the free propagator is a solution
of the free Klein-Gordon equation, and that the vertex
is mainly the coefficient of the interaction term in the Lagrangian.

\begin{figure}[ht]\epsfxsize=5in \epsfbox{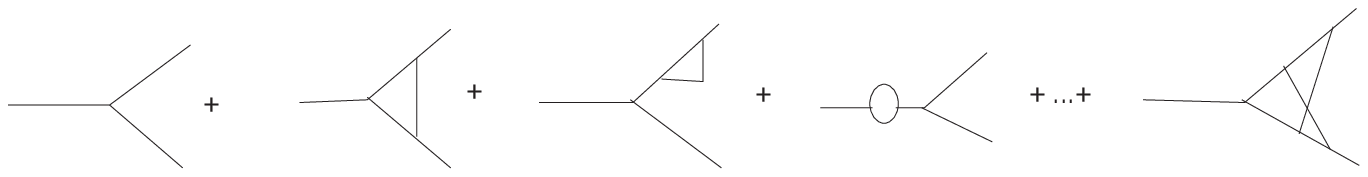} 
\caption[Loop expansion.]{\label{hab24}\small 
Increasing the number of interaction terms amounts to a loop
expansion of the theory. This time we consider an interaction
of the form $g \phi^3$, so that three free propagators
merge at a vertex. Each extra propagator provides a factor
$g^2$, so that a $n$-loop graph provides a correction
of order $g^{2n}$. These loop graphs have to be integrated over
internal loop momenta. This corresponds to a sum over all
possibilities for the momenta of the internal particles. 
Momentum conservation at each vertex, and the resulting overall
momentum conservation for external particles, are not able to fix these
internal loop momenta. They have to be integrated over the full
four dimensional $d^4q$ space.}
\end{figure}

We can now use these Feynman rules to construct the perturbative expansion
of the Green functions. Fig.(\ref{hab24}) gives examples, this time
for a theory with a cubic selfinteraction.

\begin{figure}[ht]\epsfysize=4cm \epsfbox{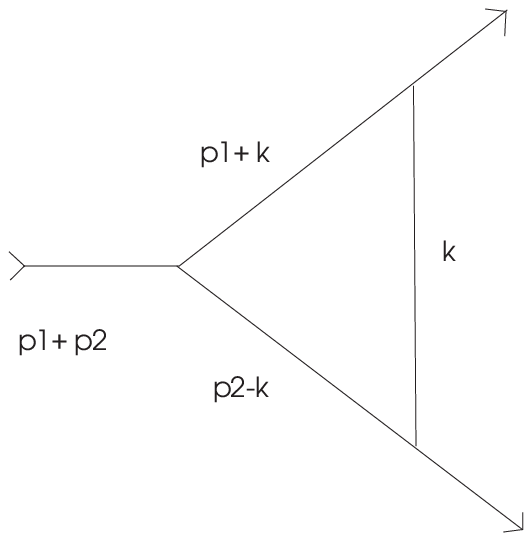} 
\caption[One-loop graph.]{\label{hab25}\small 
Such a one-loop graph demands an integration of
internal momenta. It turns out to be ill-defined. This
is a well-known obstacle in pQFT.}
\end{figure}

At this stage we pause. We can use the Feynman rules to build an
expansion
in graphs with increasing loop number, for a given Green function.
This amounts to a perturbative expansion in the coupling constant.
Such an approach has proven extremely successful in particle physics.
Nevertheless, there is a major obstacle: the presence of divergences
in these integrals with respect to the internal loop momenta.
So, having established the Feynman rules, we immediately run into
a problem. To see the obstacle consider a one-loop graph like 
Fig.(\ref{hab25}) 
for $\phi^3$ theory in six dimensions.
According to our Feynman rules we have to integrate the internal momentum.
The analytic expression is of the form
\begin{eqnarray}
\int d^6k \frac{1}{[k^2\!-m^2\!+i\eta][(k+p_1)^2\!-m^2\!+i\eta]
[(k-p_2)^2\!-m^2\!+i\eta]}\nonumber\\
\sim \int d^6k \frac{1}{k^6}\sim
\int^\infty \frac{dx}{x}. 
\end{eqnarray}
It is apparent that the integration is ill-defined.
In fact, these so-called UV-divergences are the focus of this paper.
They appear in all pQFT which are of interest in particle physics. 
The name stems from the fact that they appear at large internal
momentum, corresponding to large internal (ultraviolet) energies.
We want to understand how a theory plagued by such 
infinities can still be well-defined and able to deliver sensible 
predictions, how it can be {\em renormalized}.
\subsection{The idea of renormalization}
The theory we want to explore to understand renormalization
is $\phi^3$ theory in six dimensions ($\phi^3_6$). 
We will see soon that
precisely in six dimensions the theory provides all
typical properties of a renormalizable theory.
Other theories which we will refer to 
are given in the appendix.
We recall the Lagrangian of $\phi^3_6$:
\begin{equation}
L=\frac{1}{2}(\partial_\mu\phi)^2 
-\frac{1}{2}m^2\phi^2+\frac{g}{3!}\phi^3.
\end{equation}
Before we start to introduce renormalization which gives a meaning
to our momentarily ill-defined field theory we have to
introduce regularization. This serves the more modest
demand to give a meaning to the ill-defined loop integrals
themselves. In fact, the whole purpose is to have
a calculational tool which maps these loop integrals to some
Laurent series. This series is with respect to some
regularization parameter, and the pole terms in the series
should reflect the divergences of the integrals. 
The most common regulator is dimensional regularization.
For our purposes, it is sufficient to define it by the following
formulas \cite{Collins}
\begin{eqnarray}
\int d^Dk \frac{[k^2]^a}{[k^2-m^2]^b} & = &
\pi^{D/2}[-m^2]^{D/2+a-b}\frac{\Gamma(a+\frac{D}{2})
\Gamma(b-a-\frac{D}{2})}{
\Gamma(\frac{D}{2})\Gamma(b)},\label{dr1}\\
\int d^Dk \frac{1}{[k^2]^a\;[(k+q)^2]^b} & = & 
\pi^{D/2}[q^2]^{\frac{D}{2}-a-b}
\frac{\Gamma(a+b-\frac{D}{2})\Gamma(\frac{D}{2}-a)
\Gamma(\frac{D}{2}-b)}{
\Gamma(a)\Gamma(b)\Gamma(D-a-b)},\label{dr2}\\
\int d^Dk (k^2)^\alpha & = & 0,\;\forall \alpha.\label{dr3}
\end{eqnarray}
Whenever the parameter $D$ equals the space-time dimension $N$
of the field theory under consideration we recover the original
integrals. For small deviations from this integer, we
obtain a Laurent series in $2\epsilon:=N-D$, which serves as our
regularization parameter. 
For example, the above graph in $\phi^3$ theory 
delivers\footnote{conventions
as in \cite{Collins}. In our conventions, we would drop
the $\frac{1}{2^D\pi^{D/2}}$-term. See the appendix for our conventions.}
\begin{equation}
G(p_1,p_2,m^2)=g^3\frac{1}{2^D\pi^{D/2}2\epsilon}+g^3F(p_1,p_2,m^2)+
{\cal O}(\epsilon),
\end{equation}
where $F(p_1,p_2,m^2)$ is independent of $\epsilon$.
Note that the pole term is independent of masses and
momenta $p_1,p_2$.
In general, one obtains a proper Laurent series in $D-N$ of degree
$n$ in a $n$-loop calculation. The pole terms reflect the 
UV-divergences of the Feynman graphs. Renormalization theory
dictates how to handle these divergences. For the
moment we want to have a first look at renormalization
theory in action by handling the divergences present up to
order $g^2$ in the coupling which is the one-loop case. 
Let us rewrite the Lagrangian
\begin{equation}
L=Z_2 \frac{1}{2}(\partial_\mu\phi)^2
-Z_m \frac{1}{2}m^2\phi^2+Z_1 \frac{g}{3!}\phi^3.
\end{equation}
We want to derive conditions for these $Z$-factors
so that they could render the theory finite.
First we note that they correspond to a change in the Feynman rules:
\begin{itemize}
\item vertex: $ig\rightarrow iZ_1 g$,
\item propagator: $\frac{1}{k^2-m^2+i\eta}\rightarrow
\frac{Z_2}{k^2-Z_m m^2 + i\eta}.$
\end{itemize}
Now, for simplification, we assume to live in a world where only
this vertex correction is divergent. 
Accordingly we set $Z_m=Z_2=1$ and check if we can 
determine $Z_1$ as to render the theory finite.
The new Feynman rules
change the result for the one-loop graph 
\begin{equation}
G^{(1)}(p_1,p_2,m^2)=
Z_1^3\left[ g^3\frac{1}{2^D\pi^{D/2}2\epsilon}+g^3F(p_1,p_2,m^2)+
{\cal O}(\epsilon)\right].
\end{equation}
We know already that the pole term of the one-loop expression
is independent of masses and momenta.


Our one-loop vertex function $G^{(1)}(p_1,p_2,m^2)$ 
is of order $g^3$, while the tree
level vertex ($\sim \Gamma^{(0)}$) is $g$.\footnote{To infer the
coupling constant from the vertex function we have to
strip off the propagators for the external particles.
The resulting Green function is called the
amputated vertex function $\Gamma(p_1,p_2,m^2)$. It is obtained
most easily by multiplication with the inverted propagators
in momentum space.} We want to render the theory finite
at the one-loop level which means that we want to have a finite result
to order $g^3$ for the vertex function. 
With our modified Feynman rules, to order $g^3$
we have
\begin{equation}
\Gamma(p_1,p_2,m^2)=Z_1 g + Z_1^3 (g^3 \frac{1}{2^D\pi^{D/2}2
\epsilon})+\ldots,
\label{vert}
\end{equation}
where now $\Gamma(p_1,p_2,m^2)$ refers to the amputated Green function
where the external propagators have been removed.
We see that $Z_1$ must have the following form
\begin{equation}
Z_1=1-\frac{g^2}{2^6\pi^3 2\epsilon}.
\end{equation}
To make the one-loop renormalization
complete, we would have to include also the other $Z$-factors,
all of them of the form
\begin{equation}
Z_i=1-\sum_{j=1}^\infty \sum_{k=0}^j (g^2)^j 
\frac{c^{[i]}_{jk}}{\epsilon^k}.\label{zfac}
\end{equation}
Each monomial in the Lagrangian density obtains such a $Z$-factor.\footnote{
In our notation $c_{11}^{[1]}=\frac{1}{2(2^6\pi^3)}$.}
Theories which can be rendered finite by introducing one
$Z$-factor for each monomial in the original 
Lagrangian are called renormalizable. We more or less focus
on these theories for our purposes.\footnote{
One can relax this restriction somewhat. Calculable
theories may have an increasing number of terms
in the Lagrangian. At each loop order the
need to compensate divergences introduces new monomials
in it. As long as the number of these monomials
remains finite at each loop order, one still can give these theories
some sense. It is an open research problem to what extent the results
obtained in this paper generalize to such theories.}
Higher order corrections -loop corrections- are obtained
by considering the Green-functions for each of these monomials,
with $<\!0|T[\phi_{in}(x_1)\phi_{in}(x_2)]|0\!>$
generating divergences to be absorbed by $Z_2$ and $Z_m$,
while the divergences in $<\!0|T[\phi_{in}
(x_1)\phi_{in}(x_2)\phi_{in}(x_3)|0\!>$ demand $Z_1$ as well.

\begin{figure}[ht]\epsfysize=4cm \epsfbox{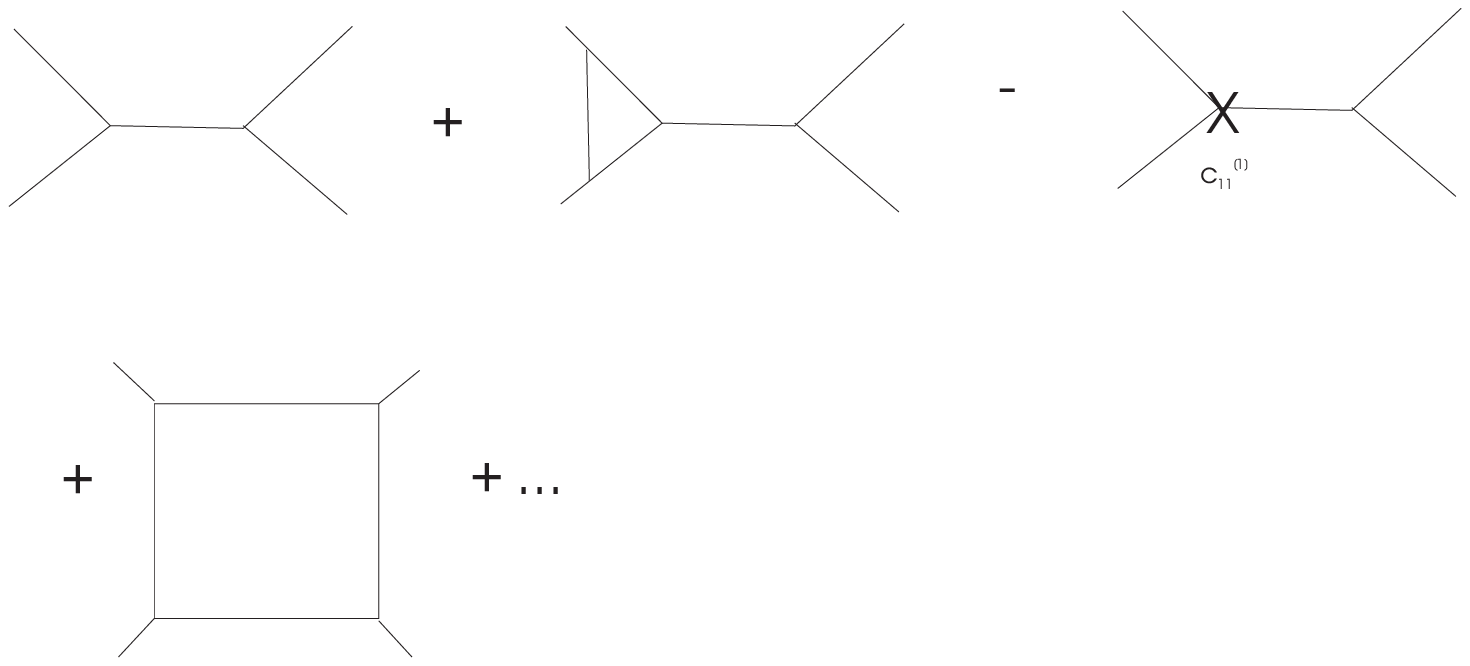} 
\caption[Localizing divergences.]{\label{hab26}\small 
All divergences of the four-point Green function
$G(x_1,x_2,x_3,x_4)$ in $\phi^3_6$ stem from
divergences of two- and three-point functions.
Our $Z$-factors remove these divergences, but the result
for the four-point function depends on the finite values of our 
{\em counterterms}. We only indicated divergences
resulting from the vertex function Eq.(\ref{vert})
and indicate the contribution from the $Z_1$ factor by a cross.
Each divergent vertex function has a compensating
counterterm at the same order in the coupling.}
\end{figure}

Let us come back to our example. Assume that
we also have determined $Z_m,Z_2$ so that the two- and three-point
Green functions are finite at the one-loop level. 
But these finite values are not determined yet.
They depend on our choice for the $Z$ factors, according to 
Eq.(\ref{zfac}).
Assume further that we now calculate some $n$-point Green function,
$n \geq 4$, at the 
one-loop level. Fig.(\ref{hab26})
gives examples for the four-point function.
Finally assume that we determine the $Z$ factors so that the values
for masses and coupling constants inferred from
Green functions $G(x_1,x_2)$ and $G(x_1,x_2,x_3)$ coincide
with the experimental values for masses and coupling
constants. Specifying such a set of {\em renormalization conditons}
determines the finite part in the $Z$-factors.
Renormalization conditions as just described are known as on-shell
conditions. We do not use them here.
We are solely interested in the divergences of the theory.
In this paper we use $MS$ conditions,
for which all $c_{j0}^{[i]}=0$. With these conditions
one only removes the pure UV divergences, so that MS $Z$-factors
provide a book-keeping of the UV singular structure
of the theory - with this choice we
subtract only the pure pole part of the Green functions.


We make this choice
only for the Green functions which suffer from divergences.
Their finite values cannot be predicted by a renormalizable
theory. They have to be taken from experiment, if one wants to test the
theory.
But then all other Green functions can be calculated in a 
unique manner.
So the result for the four-point function in Fig.(\ref{hab26}) is all
fixed. Success or failure of pQFT then depends on the answer
to the following question: Do the S-matrix elements
for processes involving Green functions like the above
four-point function (involved in $2->2$ particle scattering for example)
correctly describe experimentally observed cross sections
for these processes? Will our pQFT be
predictive for general $S$-matrix elements, once
a finite set of parameters is taken from
experiment? The answer is affirmative;
any observed particle scattering process is
successfully described by a renormalizable pQFT.
This is the motivation for pQFT, and the motivation
for the enormous effort put into loop calculations these days. 


Let us now understand renormalization in more detail.
We have to clarify the following points:
\begin{itemize}
\item 
We need a notion how to decide if a given Green function is plagued by
infinities - we have to introduce power-counting.
\item
We have to understand how subdivergences are entangled
- we have to introduce nested and overlapping subdivergences.
\item
We have to understand how to iterate the process
of cancellation of divergences with the help of $Z$-factors -
we have to understand the iterative structure of renormalization.
\end{itemize}

\begin{figure}[ht]\epsfysize=4cm \epsfbox{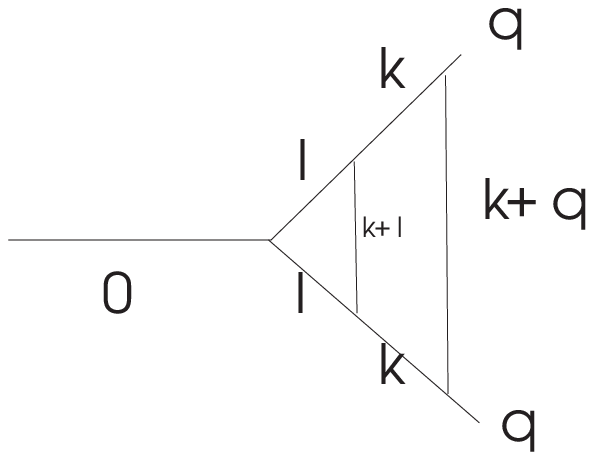}
\caption[Two-loop renormalization.]{\label{hab20}\small 
Momentum flow in the two-loop graph.
We recognize the one-loop vertex as a subgraph.}
\end{figure}

We will introduce these notions by simply considering a two-loop
example.
Consider Fig.(\ref{hab20}). It is a two-loop graph.
This is obvious by inspection. We can also infer this
from counting the number of propagators. The one-loop
vertex correction has three internal propagators.
In the next loop order we have to specify two
points somewhere in the one-loop graph, which
become internal vertices and are to be connected by a
propagator. Altogether this generates three more propagators.
In general for a theory with a cubic interaction
each loop order produces three more propagators.
So the number of propagators determines the loop order.
In our example we have realized a very special
two-loop graph. We see that the one loop graph
is nested in the whole graph as a subgraph.
So what are our expectations for the divergent structure of the graph?
We have two internal loop momenta. We route them as indicated in
the figure. The loop momentum $l$ runs solely through the one-loop
subgraph. We do expect to see a divergence when $l\to \infty$
simply because the one-loop subgraph is the divergent one-loop
vertex function of Fig.(\ref{hab25}), only this time with
momenta $0$,$k$,$k$ at its vertices.
So for any finite value of $k$ we expect a divergence from large $l$.
Now assume we keep $l$ fixed, with $k\to \infty$.
We observe that $k$ runs through a four-point function, and thus
through four propagators. It converges: $\int d^6k/k^8$ is finite
for large $k$. So we expect no divergence
from the sector ($l$ fixed, $k\to \infty$).


Let us next consider the sector where ($l\to \infty, k\to \infty$)
jointly, $k\sim l$. Altogether our two-loop graph
has six internal propagators. Each propagator behaves
as $1/l^2$ or $1/k^2$. In the region $l\sim k$
the whole integrand behaves as $1/l^{12}$, and  
the whole graph as $\int d^6k d^6l\frac{1}{(l^4(l+k)^2k^6)}\sim
\int dr/r$. So we expect a divergence for large
$l$ {\em and} $k$. We call it a logarithmic degree of divergence
for obvious reasons.


That this somewhat naive analysis is correct and the determination
of divergent sectors reduces to a counting of propagators
("edge counting") is due to Weinberg´s theorem \cite{Weinberg}. 
The resulting method we now describe
is known as power counting.


Any graph consists of a product of propagators. Let us assume
that our Feynman rules are such that all vertices are 
dimensionless.\footnote{Otherwise we have to count weighted
vertices as well in our powercounting. The weight of the 
vertices here is zero.}
Each propagator $P_i$ has a weight $c_i$. We determine this weight
as the degree of the inverse propagator,
for example for a scalar propagator $1/(p^2-m^2)$ the degree
is two, because the inverse propagator behaves as $p^2$ for
large $p$. We consider a $n$-loop Feynman graph.
Each loop contributes an $N$-dimensional integration measure,
with $N=6$ in our $\phi^3_6$ example.
Assume we have $n_i$ propagators of weight $c_i$ in the graph,
and $r$ different sorts of propagators altogether.
Then we call a graph overall divergent if
\begin{equation}
\omega:=nN-\sum_{j:=1}^r c_j n_j \geq 0.
\end{equation}
For Fig.(\ref{hab25}) and Fig.(\ref{hab20}) we have $\omega=0$.
Such graphs we call logarithmic divergent, while
$\omega=1$ or $2$ refers to linear or quadratic divergent graphs.


A few remarks are in order. Note that an addition of extra loops
will not change the degree of divergence in our $\phi^3_6$
theory. Any new loop generates three more propagators,
while on the other hand we have one more integration to do.
So the net change vanishes: $6-2\times 3=0$.
So in six dimension any three point vertex function in
$\phi^3$ theory is overall logarithmic divergent, while any
two-point function is overall quadratic divergent and
any Green function with four or more external particles
remains overall convergent. This depends
on the dimension. In $d<6$ dimensions for example
even the vertex function in $\phi^3$ is overall convergent,
and addition of further loops reduces the degree of
divergence by $d-6$. On the other hand, in $d>6$ dimension
any $r$-point Green function becomes divergent for a 
sufficiently high loop-number, as $d-6$ now worsens
the degree of divergence with each loop number.


We call $d=6$ the critical dimension of $\phi^3$ theory.
It is renormalizable in this dimension. 
Precisely in this critical dimension the coupling
constant is dimensionless. Only the
Green functions corresponding to fundamental monomials
in the Lagrangian evolve divergences in loop integrals.
In $d>6$ dimensions the theory becomes non-renormalizable.
Any Green function diverges at a sufficiently high loop
number. The coupling constant has negative dimension.
In $d<6$ dimension the theory is superrenormalizable.
Divergences appear only at a sufficiently low loop number
and the coupling constant has positive dimension.


These notions of criticality and renormalizability extent
to other pQFT in an obvious manner. All other theories we refer
to in this paper are critical in four dimensions.
Particle physics seem to favour pQFTs which are critical in four
dimensions.


Back to our two-loop example in Fig.(\ref{hab20}).
Power counting tells us that it has an logarithmic overall 
degree of divergence. We also expect a divergence from its
one-loop subgraph. 
We can indeed extend the notion of powercounting in an obvious
manner. Let us in our mind put a little box around any divergent 
subgraph $g_i$ which we detect in a Feynman graph $G$. We do a
powercounting for the Feynman graph in this box, and say that
the degree of divergence of the subgraph $g_i$ 
is its overall degree of divergence
$\omega_i$.
We then say that $G$ has an overall degree of divergence $\omega$ 
and subdivergences
contributed by subgraphs $g_i$ of degree of divergence $\omega_i$.

\begin{figure}[ht]\epsfysize=4cm \epsfbox{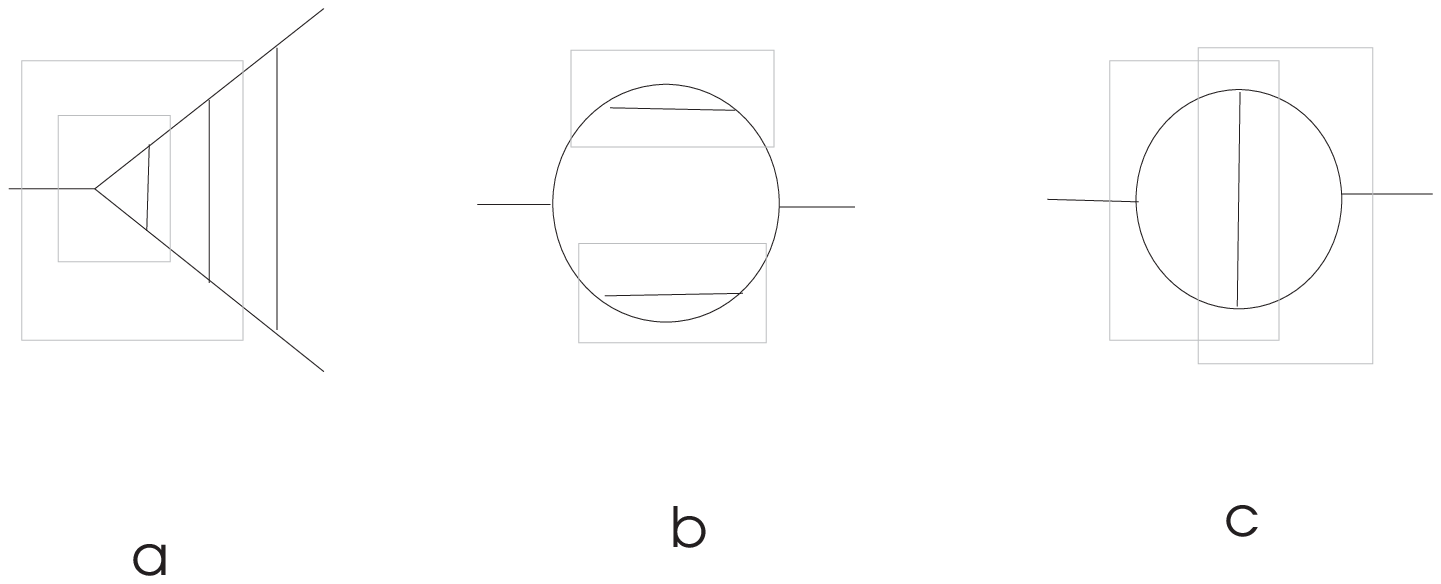} 
\caption[Classes of subdivergences.]{\label{hab27}\small 
There are various ways subdivergences can appear.
They can be nested as in (a), disjoint as in (b) or overlapping
as in (c). We put every subdivergence in a dashed box. These boxes are
called forests. Forests which are not nested in a larger one
are called maximal.}
\end{figure}

Fig.(\ref{hab27}) provides examples for various different
ways how subdivergences can appear. They will all
be considered in the following sections.
They can come in nested
form and as disjoint subdivergences, and they can appear to be overlapping.
This last form is the most complicated one.
Renormalization theory is an iterative procedure 
to eliminate the divergences of the theory. One first
eliminates the subdivergences, and then hopes to control
the remaining overall degree of divergence.
This iteration procedure can be summarized in a forest formula
\cite{BPHZ,Itz} but we now rather study our two-loop example
to see how the iteration works.

\begin{figure}[ht]\epsfysize=3cm \epsfbox{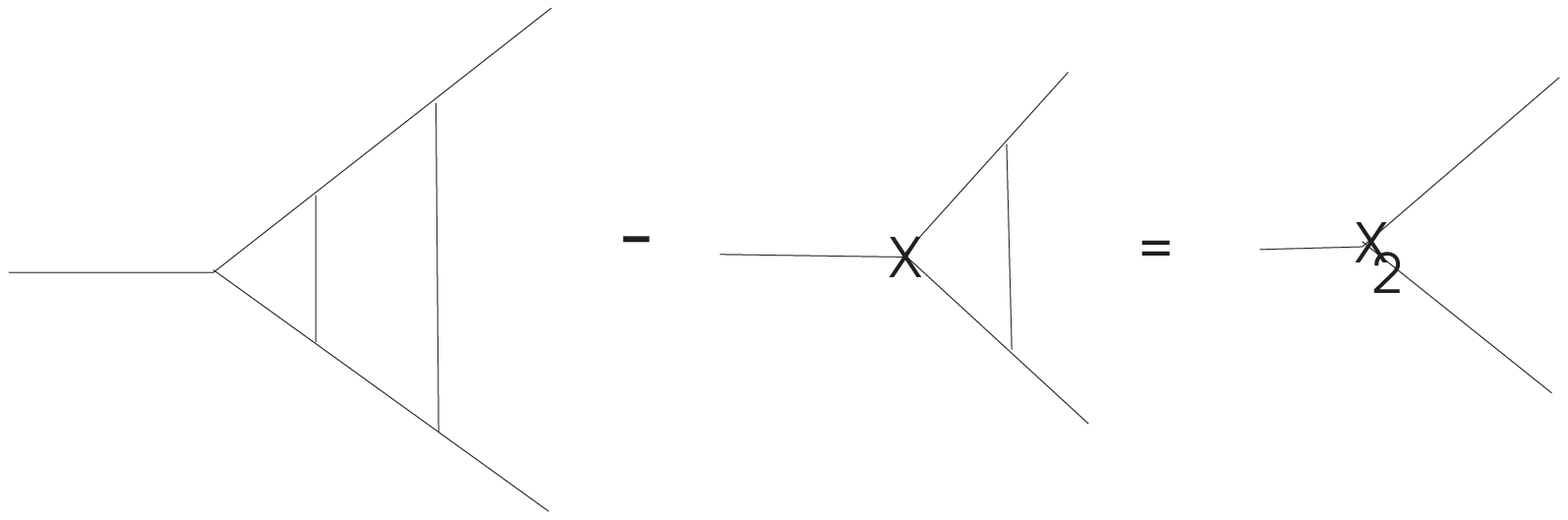} 
\caption[Adding counterterm contributions.]{\label{hab28}\small 
We add a counterterm graph. It absorbs the subdivergence
in the subgraph, present for $l\to \infty$, $k$ fixed.
We should also add counterterm contributions at the
other two vertices. They would absorb the subdivergences of the two-loop
graphs where the one-loop subgraph is located at these other vertices.
After subtraction of the one-loop counterterm graph we can determine
the two-loop $Z$-factor, indicated by $X_2$.}
\end{figure}

We have already added a one-loop $Z_1$-factor to our Lagrangian.
We are interested now to obtain a finite result at order $g^5$.
At this order we not only have the two-loop graph, but also a
one-loop graph with a counterterm insertion from the renormalization
of the one-loop divergence. 
This term appears due to our replacement
$g\phi^3\rightarrow Z_1 g\phi^3=g\phi^3-
g^3\frac{c_{11}^{[1]}}{\epsilon}\phi^3$. 
We have to add this counterterm graph
to the two-loop graph, as it is done in Fig.(\ref{hab28}).
The counterterm was determined to render the one-loop vertex
function finite. This function appears as a subgraph in the full
two-loop graph, and it is indeed so that for any fixed $k$
the difference of the full two-loop contribution and the counterterm
contribution is finite. The only remaining divergence
appears when both loop momenta tend to infinity jointly.
This is the overall divergence. This is true in general.
After elimination of all subdivergences by counterterms
the remaining divergences appear when all loop momenta
tend to infinity jointly.


So now our renormalization program continues as for the one-loop
case. The difference of the full two-loop graph and counterterm
contribution has the form
\begin{equation}
g^5[\frac{c_2}{\epsilon^2}+\frac{c_1}{\epsilon}+\ldots]
\Rightarrow c^{[1]}_{21}=-c_1, c^{[1]}_{22}=-c_2.
\end{equation}
and the finite parts $c^{[1]}_{20}$ must again be determined
by the renormalization conditions, and do
vanish in the $MS$ scheme.\footnote{For $\phi^3_6$
one finds explicitely $c_1=\frac{1}{16(2^{2D}\pi^D)}$ and $c_2=\frac{-1}{
8(2^{2D}\pi^D)}$.}


It is only in this combination of the full two-loop graph
with the counterterm contribution that the $c^{[1]}_{2i}$, $i\in \{1,2\}$,
come out as pure numbers as they should. The two-loop
graph itself will contain divergences of the
form $ \log{q^2}/\epsilon$, where the divergence comes from
the divergent subgraph, and multiplies a $\log{q^2}$ which
appears in the finite part of the second loop integral.
In fact, much more complicated non-polynomial functions
of $q^2$ and masses could be present. Fortunately,
after addition of the counterterm contribution, these
non-local contributions cancel, and the $c^{[1]}_{2i}$
appear as pure numbers. This is a general property as well:
after adding the appropriate counterterms for subdivergences
the overall divergences remain local.\footnote{Only
polynomials in momenta give local operators after Fourier transformation.}
We will see later explicitly how to handle disjoint
or overlapping cases. 


The whole renormalization procedure
 was remarkably successful in QED as well as in the standard
model of elementary particle physics and still provides the only
systematic mechanism for calculations in elementary particle physics.

\clearpage
\section{Planar Vertex Corrections}
In this section we explore algebraic structures
of nested divergences as generated by ladder graphs.
Ladder topologies are the easiest accessible and 
serve
as a convenient means to introduce our approach. 
Here we will for the first time draw a connection between
renormalization theory and link diagrams. Subsequent
sections will show that the connection goes via number theory:
the trivial nested topologies considered here provide rational coefficients
for the UV divergences, while transcendentals encountered in
more complicated examples faithfully describe the knots 
derived via link diagrams from our Feynman graphs.
\subsection{The vertex ladder}
We consider theories with renormalizable three-point couplings.
As pointed out in the previous remarks on renormalization theory
this means for us a
theory with vertex corrections of logarithmic degree of divergence,
corresponding to a dimensionless coupling constant.
Further, for simplification, we again restrict ourselves to
vertex corrections. We collect our conventions in an appendix.

\begin{figure}[ht]
\epsfysize=4cm \epsfbox{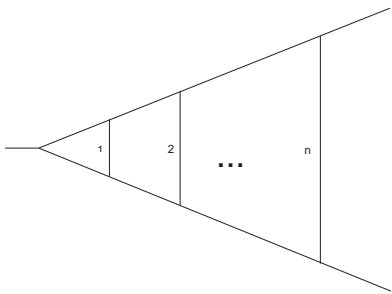}
\caption[The ladder topology.]{\label{hab31}
\small The ladder topology for vertex corrections.
We indicate loop momenta $l_1,\ldots,l_n$.
The ladder 
topology  might be realized by different Feynman graphs.
Our graphical representation of Feynman diagrams usually refers only
to the topology of the graph, so that lines may describe various 
different particles
according to the Feynman rules of the theory.}
\end{figure}
\begin{figure}[ht] \epsfysize=3cm \epsfbox{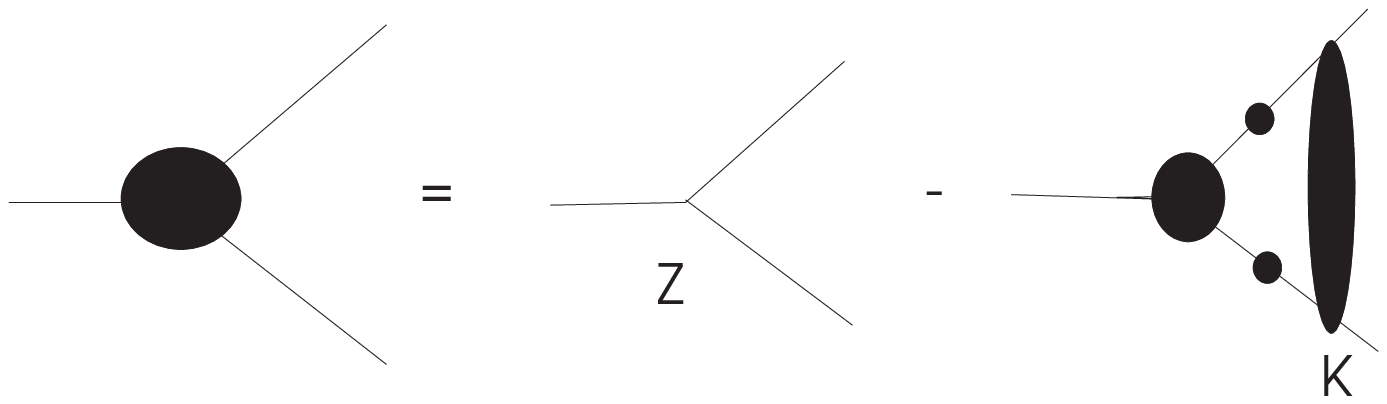}
\caption[The Schwinger Dyson equations.]{\label{hab32}\small 
We can generate the full three-point Green function
iteratively. An integration over the internal momentum in the loop
on the rhs is understood. The kernel $K$ allows a skeleton
expansion. Usually, all propagators and internal
vertices represent the full Green-functions, and the coupled
system of equation for these functions is hard to solve.
Full Green functions are indicated by black blobs.
In the ladder approximation considered here we restrict ourselves
to bare propagators (indicated by straight lines), 
bare internal vertices (no blobs) and to the
lowest order in the skeleton expansion.}
\end{figure}
\begin{figure}[ht] \epsfysize=3cm \epsfbox{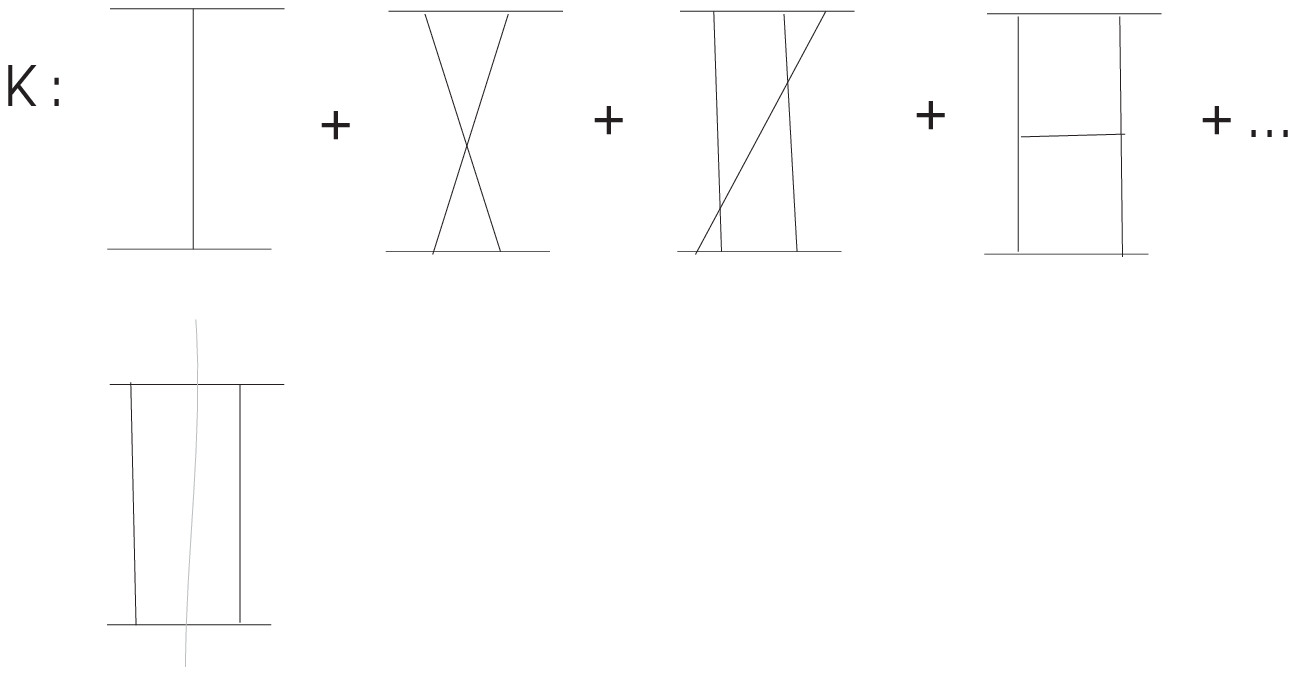}
\caption[The skeleton expansion.]{\label{hab33}\small
The skeleton expansion given for bare
internal propagators and vertices. In its dressed form (reinserting
blobs for very internal vertex or propagator) it generates
the kernel $K$ of the Schwinger Dyson equations.
Note that the lowest order skeleton is the one-particle exchange
graph of order $g^2$.
We also indicate a graph which is not a skeleton.
It is an iteration of the lowest order skeleton.
It separates in two pieces when we cut it a two propagators,
and thus it is two-line reducible.}
\end{figure}

Let us simplify for a start even further 
to topologies of the form given in Fig.(\ref{hab31}).
We have omitted all dressing of the
internal vertices and propagators here. 
In general, the iteration of loops is controlled
by the Schwinger-Dyson equations.
We explain these equations in Fig.(\ref{hab32}).
These equations nicely organize the perturbative expansion.
Each Feynman diagram consists of propagators and
vertices. These are themselves Green functions which allow
for a perturbative expansion. The Schwinger Dyson equations
are coupled integral equations which iterate these
Green functions in terms of themselves. Knowing their
full solution would correspond to knowing the full non-perturbative
Green functions. A peculiar role is played by the skeleton diagrams.
They furnish four-point functions which serve as a kernel
in the Schwinger Dyson equations.
These skeleton diagrams are Feynman diagrams which are the basic blocks
of iteration in the equations. Any proper skeleton contribution
is not two-line reducible. When we close two of its lines
we get a vertex function free of divergences, as long as
the vertices and propagators of the skeleton were undressed,
cf.~Fig.(\ref{hab33}). 


In our limit, the Schwinger Dyson equation generates the
ladder graphs by iteration of the lowest order in the
skeleton:
\begin{eqnarray}
\Gamma(p_1,p_2) & = & Z_g \Gamma^{0}\nonumber\\
 & & -\int d^Dk\;
\Gamma^0 D(k)\; \Gamma^0 \;D(k+p_1)\;\Gamma(k+p_1,k-p_2)\;D(k-p_2),
\label{sd}
\end{eqnarray}
here $\Gamma^{0}$ is the tree-level vertex, and $D$ the bare propagator,
which is the inverse of the free Klein-Gordon equation.\footnote{Or
the inverse of any other appropriate equation of motion.}
The equations thus iterate the full vertex $\Gamma$ in terms of itself
and the propagator function $D$. 

\begin{figure}[ht] \epsfysize=5cm \epsfbox{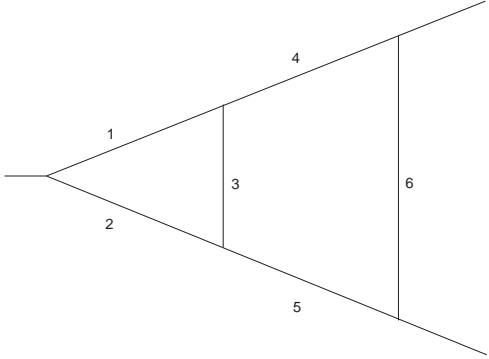}
\caption[Two-loop ladder.]{\label{hab34}\small The two-loop function
of the previous section. Our considerations
are now valid for any Green function in a pQFT realizing
this topology. We denote propagators by numbers for later
reference.}
\end{figure}

For a start let us simplify even more and consider the two-loop
case as in Fig.(\ref{hab34}).


Normally, one would calculate this graph, then calculate its counterterm 
graphs (thereby incorporating the one-loop $Z$-factor),
and finally add all contributions. This is what we did in the
previous section. In the sum one would observe all
the expected cancellations necessary to give a local two-loop $Z$-factor.
One might do better by subtracting out the 
subdivergences from the very beginning.
Thus, subtract from this vertex correction $\Gamma^{(2)}$ another 
vertex function,
forming a new function $\Gamma^{(2)}_{\Delta}$:
\bea
\Gamma^{(2)}_{\Delta} := \Gamma^{(2)} - \tilde{\Gamma}^{(2)},
\label{e2}\eea
where we define
\bea
\tilde{\Gamma}^{(2)} := \Gamma^{(2)}\mid_{m_1=m_2=m_3=p=0};
\label{e3}\eea   
in $\tilde{\Gamma}^{(2)}$ we have set all masses in the subdivergence
to zero and evaluate the vertex at zero momentum transfer.
Due to the fact that the subdivergence itself is independent of
the masses and exterior momenta,
$\Gamma^{(2)}_\Delta$ is a finite Green function: 
\bea
\Gamma^{(2)}_{\Delta}=
\int (\Gamma^{(1)}- \tilde{\Gamma}^{(1)})K,\label{fin}
\eea
is finite. This is also clear as the iteration
of the integral kernel
$K$, -which is a finite four-point function-, with the finite 
$(\Gamma^{(2)}- \tilde{\Gamma}^{(2)})$ has a subtracted form.

\begin{figure}[ht] \epsfysize=3cm \epsfbox{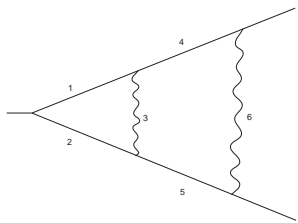}
\caption[A QED example.]{\label{hab35}\small The QED graph corresponding to the 
example Eq.(\ref{e5}).}
\end{figure}

\noindent Three examples:\footnote{For all these examples
one can consider the difference in Eq.(\ref{fin}). For power counting
purposes one considers the equation on a common denominator
and easily verifies its finiteness.}\\
for $\phi^3$ in six dimensions we have
\bea
\tilde{\Gamma}^{(2)}=\int d^Dld^Dk\; \frac{1}{l^4(l+k)^2(k^2-m^2)^2((k+q)^2-m^2)},
\label{e4}\eea
with a tree-level vertex $\Gamma^{0}=1$,
while for massive QED in four dimensions, in the Feynman gauge
\bea
\tilde{\Gamma}^{(2)}_{\sigma}=
\int 
d^Dl d^Dk\; \gamma_\mu \frac{1}{\kslash-m}\gamma_{\rho}
\frac{1}{\lslash}\gamma_{\sigma}\frac{1}{\lslash}
\gamma^{\rho}\frac{1}{\kslash-m}\gamma^{\mu}\frac{1}{(l+k)^2}
\frac{1}{(k+q)^2},
\label{e5}\eea
corresponding to the Feynman graph in Fig.(\ref{hab35}).
Here, the tree level vertex is $\gamma_\sigma$.
A convenient example is also provided by 
Yukawa theory, the coupling of a fermion to a scalar particle.
For a massive fermion, we obtain
\begin{equation}
\tilde{\Gamma}^{(2)}=
\int 
d^Dl d^Dk\;  \frac{1}{\kslash-m}
\frac{1}{\lslash}\frac{1}{\lslash}
\frac{1}{\kslash-m}\frac{1}{(l+k)^2}
\frac{1}{(k+q)^2},
\end{equation}
again with a tree-level vertex $\Gamma^{0}=1$. 
Note that in the above examples the propagators of the innermost
loop are massless, according to our
definition of $\tilde{\Gamma}$.
We see that in all examples the $l$-integration over the subdivergence 
becomes
particularly simple. It depends solely on the other loop momentum.


Apart from the UV divergences at large internal momentum
there are no other divergences present. In particular,
we do not encounter
IR singularities (which appear for small
internal momenta) in any of our considerations.


Now consider $\Gamma_{\Delta}^{(2)}$. It is UV convergent, both in its 
subdivergent and overall divergent behaviour, which can as well be concluded
from power-counting. Our subtraction  at the
subdivergence improves the power-counting
for the inner loop momentum by one degree, but it improves also the
overall degree of divergence by the same amount, 
as the subtracted term was only
modified in exterior parameters like masses and momenta, and the overall
degree of divergence is independent of these parameters \cite{mpla}.


Accordingly, we have shown that the UV divergences of $\Gamma$ are located
in a simpler Green function $\tilde{\Gamma}$, where we could set some
masses  to zero and evaluate the whole function
at vanishing momentum transfer.


From now on we make use of the following notation:
\bea
<\Gamma>= <\Gamma_{\Delta}+\tilde{\Gamma}>
=<\Gamma_{\Delta}>+<\tilde{\Gamma}>=<\tilde{\Gamma}>.
\label{e6}\eea
Here we introduce a projector $<\ldots>$ onto the UV-divergences 
(the proper part of a Laurent expansion of $\Gamma$ in
$\epsilon$, where $\epsilon$ is the DR regularization parameter)
so that, by definition,
\bea
<\mbox{UV-finite expression}>=0.
\label{e7}\eea
We also use the $n$-th order complete $Z$-factor and various
Green functions defined as follows:
\bea
{\bf Z}_1^n & := & 1-c_1^{(1)}-c_1^{(2)}\ldots-c_1^{(n)},\nonumber\\
{\bf \Gamma}^n & := & \Gamma^{(0)}+\Gamma^{(1)}\ldots+\Gamma^{(n)},\nonumber\\
\tilde{\Gamma}^{(n)} & := & \Gamma^{(n)}|_{m_1=m_2=m_3=p=0},\nonumber\\
\bar{\Gamma}^{(n)} & := & \tilde{\Gamma}|_{m_i=0\forall i}.
\label{e8}\eea
We redefined our notation for the $Z$ factor somewhat.
Compared with Eq.(\ref{zfac}) 
the $c_i$ incorporate all the pole terms to a given loop order.
Note that we have two modified Green functions: 
$\tilde{\Gamma}^{(n)}$, which has only a modification
in its innermost loop, and $\bar{\Gamma}^{(n)}$, where all
masses are set to zero. We consider the case of zero momentum
transfer throughout. 


For the ladder diagrams, the index $n$ coincides with the loop order.
Now, having absorbed the UV-divergence of $\Gamma^{(2)}$ in 
$\tilde{\Gamma}^{(2)}$
let us add the counterterm graph
\bea
\tilde{\Gamma}^{(2)}\rightarrow
Z_1^{(2)}:=<\tilde{\Gamma}^{(2)}-Z_1^{(1)}\tilde{\Gamma}^{(1)}>,
\label{e9}\eea
$Z_1^{(2)}$ contains $<\Gamma^{(2)}>$ plus
its counterterm graph. We will become more acquainted
with its use in the following.


With the above definitions we have
\bea
<\tilde{\Gamma}^{(2)}-\bar{\Gamma}^{(2)}-Z^{(1)}
(\tilde{\Gamma}^{(1)}-\bar{\Gamma}^{(1)})>=0,
\label{e10}\eea
since the overall divergent behaviour of
\beas
(\tilde{\Gamma}^{(2)}-Z_1^{(1)}\tilde{\Gamma}^{(1)}),
\eeas
and
\beas
(\bar{\Gamma}^{(2)}-Z_1^{(1)}\bar{\Gamma}^{(1)}),
\eeas
are the same.
This is exactly what we expect. Both expressions have the same 
asymptotic behaviour with regard to their overall degree of divergence.
They further have the same subdivergence $Z_1^{(1)}$. 
After eliminating
this subdivergence by an appropriate subtraction in the first step
the difference of the two terms is necessarily UV convergent. 
This is nothing else than the statement that
$Z_1^{(2)}$ must be polynomial in exterior masses and momenta. 
Powercounting confirms that such trivial subtractions render
our two-loop examples finite.


By use of Eq.(\ref{e10}) we have absorbed all UV 
singularities in the expression
\bea
Z_1^{(2)}=\bar{\Gamma}^{(2)}-Z_1^{(1)}\bar{\Gamma}^{(1)}.
\label{e12}\eea
It is significant that {\em 
only massless three-point functions at zero momentum transfer} 
appear in the above expression. The significance stems
from the fact that now
each loop integration depends only on one parameter: the
momentum of the next loop integration. Only the last loop
integration depends on the exterior momentum $q$.


The one-loop $Z$-factor from the vertex correction,
 $Z_1^{(1)}$, has the form
\bea
\bar{\Gamma}^{(1)}=
\Delta(\epsilon)\;(q^2)^{-\epsilon}
\Gamma^{(0)}\Rightarrow
Z_1^{(1)} = <\bar{\Gamma}^{(1)}>,
\label{e13}\eea
where we used the fact that $\bar{\Gamma}^{(1)}$ scales like
$(q^2)^{-\epsilon}$ in DR. 
This equation means the following: $\bar{\Gamma}^{(1)}$
can only depend on one variable: $q^2$, the square of the
external momentum. Its dependence is determined by the
scaling properties of the Green function. $\bar{\Gamma}^{(1)}$
also factorizes the tree level vertex $\Gamma^{(0)}$.
The $q^2$ dependence is a trivial scaling $\sim
(q^2)^{(D-N)/2}$, so that $\bar{\Gamma}^{(1)}(q^2)=(q^2)^{(D-N)/2}
\bar{\Gamma}^{(1)}\mid_{q^2=1}$, where the exponent $(D-N)/2$
follows from dimensional considerations. 
$\frac{\bar{\Gamma}^{(1)}\mid_{q^2=1}}{\Gamma^{(0)}}$ 
is then a sole function
of $\epsilon=(D-N)/2$, which is the content of the 
previous equation.


As we are calculating in a MS scheme we were allowed to project onto the
divergent part in the Laurent expansion of  
$\bar{\Gamma}^{(1)}\mid_{q^2=1}$.
The fact that all the dependence on the exterior momenta
$q$ is in the scaling behaviour $(q^2)^{-\epsilon}$
is important. It allows to define a function
$\Delta$ which only depends on the DR parameter $D$ and incorporates
all the information to obtain MS $Z$-factors. 


Let $F(k,q)$ be the integrand for the calculation of
$\bar{\Gamma}^{(1)}$.
We see that (using Eq.(\ref{e13})) 
\bea
\bar{\Gamma}^{(2)}=
\Delta \int d^Dk (k^2)^{-\epsilon} F(k,q)=
\Delta {}\;{}_1\!\Delta\;\;(q^2)^{-2\epsilon}\;\Gamma^{(0)},
\label{e16}\eea
where we define for later use
${}\;{}_j\!\Delta$ by
\bea
\int d^Dk (k^2)^{-\epsilon j} F(k,q)
=: {}\;{}_j\!\Delta \;\;(q^2)^{-
\epsilon(j+1)}.
\label{e17}\eea


We consider the functions ${}_j\!\Delta$ as modified one-loop functions.
For any renormalized theory they can be obtained from the corresponding standard
one-loop integral by a change in the measure
\beas
\int d^Dk \rightarrow \int d^Dk (k^2)^{-\epsilon j}.
\eeas


Throughout we tacitly assumed that the massless one-loop vertex function
at zero momentum transfer only reproduces the tree-level vertex 
$\Gamma^{(0)}$ as 
a form factor (Eq.(\ref{e13})). This is not necessarily so.
Even a theory as simple as massless QED delivers  a
counterexample, by providing two formfactors, generated
by the two spin structures $\gamma_\sigma, q_\sigma \qslash$.
But it turns out that
the modifications due to this complication do not
spoil our general reasoning. We will present the necessary
changes
in an appendix.


We denote the parameter $j$ above as the `writhe number' in the following,
for reasons which will become clear below.


Finally, for the two-loop case
\bea
Z_1^{(2)} = <\Delta  {}\;{}_1\!\Delta-<\Delta>\Delta>,
\label{e18}\eea
by virtue of Eq.(\ref{e13},\ref{e16}).



Our next step will be to generalize the above approach 
to the $n$-loop case. We still pretend that the world consists only 
of ladder-type vertex corrections, in which case our approach
would deliver the full $Z$-factor for the vertex correction. 

\begin{figure}[ht] \epsfysize=5cm \epsfbox{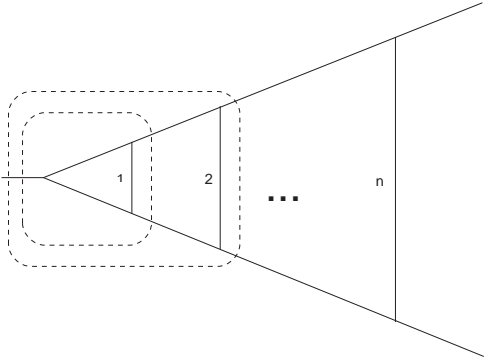}
\caption[Divergent sectors.]{\label{hab36}
\small  The dashed boxes denote the divergent sectors of the ladder topology.
These give rise to the restricted set of counterterms we are interested in.
The divergent sectors are, from the left to the right, the domains
$\lambda\to\infty$ in $\{\lambda l_1\},\{\lambda l_1,\lambda l_2\},
\ldots,\{\lambda l_1,\ldots,\lambda l_{n-1}\}$.}
\end{figure}

Again we define 
\bea
\Gamma^{(n)}_\Delta:=\Gamma^{(n)}-\tilde{\Gamma}^{(n)},
\label{e19}\eea
We still evaluate at zero momentum transfer. By setting the masses 
$m_1,m_2,m_3$ of the
inner loop to zero in the second term we remove {\em all}
subdivergences of the graph.
It follows that $<\Gamma^{(n)}_{\Delta}>=0$.
This is obvious as all existing forests include
the innermost loop as a nest. All divergent sectors in Fig.(\ref{hab36}),
which are all logarithmic divergent, thus have an
improvement in power counting when we have such an
improvement for this inner loop. It follows as well by
considering $\Gamma_\Delta=\int\ldots\int (\Gamma^{(1)}-\tilde{\Gamma}^{(1)})
K\ldots K$, and writing out the iteration of the bubble in the Dyson-Schwinger
equations explicitly. 
To calculate $Z_1^{(n)}$ we have to consider
\bea
Z_1^{(n)} = <\bar{\Gamma}^{(n)}
-\sum_{i=1}^{n-1}Z_1^{(i)}\bar{\Gamma}^{(n-i)}>.
\label{e20}\eea
The sum on the rhs adds the appropriate counterterms which
compensate all subdivergent sectors of the full Green function.
We are allowed to consider only massless quantities here 
as the overall degree of divergence is mass independent, and the term in 
angle brackets is free of subdivergences by construction.

\begin{figure}[ht] \epsfysize=7cm \epsfbox{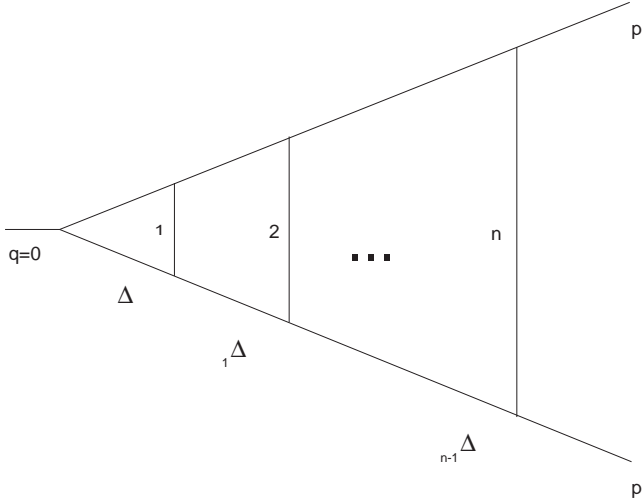}
\caption[Concatenated one-loop functions.]{\label{hab37}\small 
By evaluating at zero momentum transfer the calculation
of a massless ladder becomes a concatenation of ${}_j\!\Delta$ 
functions.}
\end{figure}

We have  
\bea
\bar{\Gamma}^{(n)}=
(q^2)^{-\epsilon n}\prod_{i=0}^{n-1}{}\;{}_i\!\Delta,
\label{e21}\eea
according to our definition Eq.(\ref{e17}).
Here sequential expressions like $\Delta{}_1\!\Delta\ldots$ denote
a concatenation of one-loop functions as shown in Fig.(\ref{hab37}). 
Again, 
in case that there is  more than
one form factor in this concatenation  some technical
subtleties are involved which are given in the appendix.
By investigating counterterm expressions and Green-functions we see
that only two type of operations on ${}_j\!\Delta$ are
needed to give the $Z$-factor: a concatenation of these
functions for the Green functions, and a projection onto
their divergent parts to find counterterms.  
This justifies to define two operators
\bea
B^k(\Delta) & := & \prod_{i=0}^{k}{}\;{}_i\!\Delta,\nonumber\\
A^r(\alpha) & := & \Delta<<<\ldots<\alpha >\Delta>\ldots \Delta>,\nonumber\\
 & & \mbox{r $<>$ brackets}\nonumber\\
\Rightarrow A^r(B^k(\Delta)) & = &  
\Delta<<<<\ldots<\prod_{i=0}^{k}{}\;{}_i\!\Delta
>\Delta>\ldots\Delta>,\nonumber\\
B^k(A^r(\Delta)) & = & 
<<<\ldots<\Delta >\Delta>\ldots \Delta>
\prod_{i=0}^{k}{}\;{}_i\!\Delta,\nonumber\\
B^0(\Delta)=A^0(\Delta) & = & \Delta.
\label{e22}\eea
$B$ acts by concatenating massless one-loop functions with increasing
writhe number, $A$ by projecting on the divergent part 
iteratively.


We can now give the general result for
$Z_1^{(n)}$:
\bea
Z_1^{(n)}= <[-A+B]^{(n-1)}(\Delta)>.
\label{e23}\eea

\begin{figure}[ht] \epsfxsize=5in \epsfbox{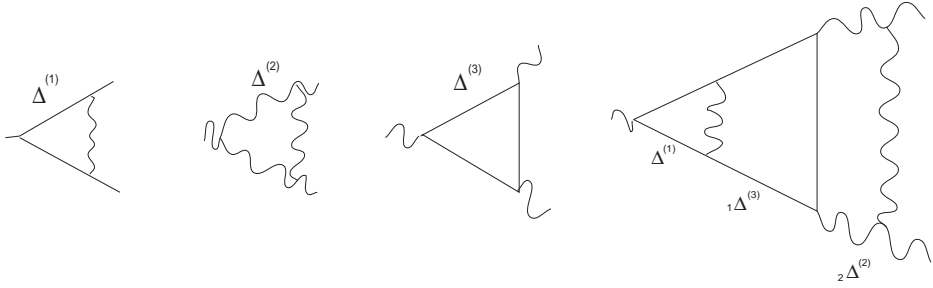}
\caption[Various different $\Delta^{(i)}$-functions.]{\label{hab38}\small
Various different 
$\Delta^{(i)}$'s. We understand that in our result Eq.(\ref{e23}) we have
to consider the appropriate combination of $\Delta$'s as arguments for
the $A$ and $B$ operators. Examples here are taken from
non-abelian gauge theories which provide various kinds of
couplings. Here we see cubic selfcouplings of the gauge boson
and the usual coupling of a fermion to a gauge boson.}
\end{figure}

For $n=3$ this gives explicitly
\bea
Z_1^{(3)}=
<\Delta{}\;{}_1\!\Delta{}\;{}_2\!\Delta-
\left[<\Delta{}\;{}_1\!\Delta>
-<<\Delta>\Delta>\right]\Delta\nonumber\\
-<\Delta>\Delta{}\;{}_1\!\Delta>.
\label{e24}\eea
It should be clear that 
in the above expressions we have on the rhs the proper three-loop
Green-function, next the two-loop $Z$-factor times the one-loop
Green function, and the one-loop $Z$ factor times the
two-loop function.
We remind ourselves that the functions ${}_j\!\Delta$ are simple
modifications of the corresponding three-point one-loop function,
evaluated in the massless limit at zero momentum
transfer. The above concatenation properties are universal
for all these functions, as long as they are generated
by a renormalizable theory. We note that there might be in fact a whole
`bouquet' of these functions at our disposal.
This is so as in more realistic
theories our ladder topology might consist of various different
$\Delta^{i}$'s, as in the  example of Fig.(\ref{hab38}).

\begin{figure}[ht] \epsfysize=5cm \epsfbox{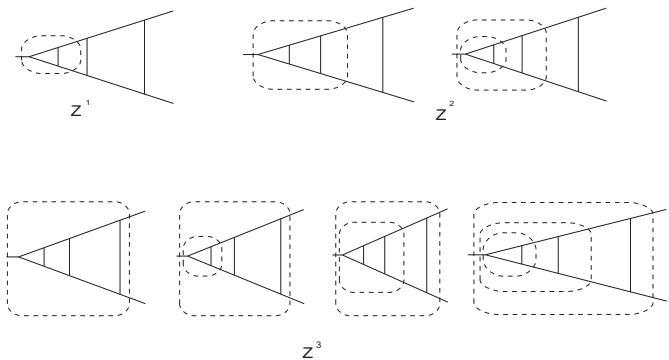}
\caption[Divergent sectors.]{\label{hab39}\small 
We have calculated $<\Gamma^{(n)}>$ so far, including
the counterterms for the dashed divergent
sectors, taking into account all possible forests.}
\end{figure}

Recalling our definition for the full vertex function  
\bea
{\bf \Gamma}^{(n)} & := & \Gamma^{(0)}+\sum_{i=1}^{n}\Gamma^{(i)},\nonumber\\
\Rightarrow <{\bf \Gamma}^{(n)}> & = & <{\bf \bar{\Gamma}}^{(n)}>,  
\label{e25}\eea
we have its MS renormalized form,
\bea
{\bf \Gamma}_R^{(n)}={\bf Z}_1^{(n)}{\bf {\Gamma}}^{(n)}.
\label{e26}\eea
We can easily check that ${\bf Z}_1^{(n)}$ renders ${\bf \Gamma}_R^{(n)}$
finite.
Indeed, we only have to check that ${\bf \bar{\Gamma}}_R^{(n)}$ is
finite and easily calculate this expression to be
(using Eq.(\ref{e21},\ref{e23}))
\bea
 \bar{{\bf \Gamma}}_R^{(n)} & = & 
 (1+\sum_{i=0}^n [B]^i(\Delta))(
1-\sum_{i=0}^{n-1}<[-A+B]^i(\Delta)>),\nonumber\\
&  = & 1+\sum_{i=0}^{n-1}([-A+B]^i(\Delta)-<[-A+B]^i(\Delta)>),
\label{e27}\eea
which is evidently finite.
In some obvious shorthand notation we can write this
result to all orders as
\bea
 {\bf \bar{\Gamma}}_R=\frac{1}{1-[-A+B]}(\Delta)-<\frac{1}{1-[-A+B]}
(\Delta)>.
\label{e28}\eea
In our conventions we never have written the coupling constants
explicitly. But each application of $A$ or $B$ increases the
loop number by one, and thus the order of our expansion.
Thus an expansion
in the coupling constant is understood in Eq.(\ref{e28}).
($A^rB^m$ is of order ${\cal O}(g^{2(r+m)})$).
Before we extend this rather academic example to more realistic situations,
 including the other topologies, self-energies and the like
let us summarize what we have achieved so far: mainly 
we have calculated the 
ladder-type
vertex corrections of Fig.(\ref{hab37}), 
including the counterterms
as depicted in Fig.(\ref{hab39}).
Our intention is to generalize this to all possible 
topologies, including dressing of internal lines and vertices, and,
crucially, including overlapping divergences. But before we do so
in the next sections, let us have another look at our ladder.
At this stage we have introduced a sufficient amount of notation
so that we can start to consider our example again, this time turning to 
knot- and link-theory as another approach to obtain our results Eq.(\ref{e23}).
\subsection{Links and ladders}
In the previous section we reduced the renormalization procedure of
ladder graphs to algebraic operations on concatenated
one-loop functions. We now want to compare this result with 
the behaviour of link diagrams in the presence of a skein
relation acting on the diagram. 
It will turn out that we identify a link diagram
with the Green-function plus all its counterterms.
Individual terms in renormalization theory should correspond to individual
terms obtained after skeining the link diagram. In the course
of this skeining, one generates from a $n$-component link diagram
in particular a one-component link diagram, a knot, amongst other terms.
Later we see that this knot gives us some information about
the numbers we should expect when we calculate the divergences of the 
corresponding graph. We use some elementary notions
of link and knot theory, as provided for example in \cite{knots}.


According to the rules we give below, the individual components 
in the link diagrams of interest are all trivial circles (unknots
from the viewpoint of knot theory).
Only through the entanglement do we 
generate non-trivial topological structure
in the link diagrams. Using a skein relation to disentangle the diagrams
we will then generate knots in intermediate steps of this algorithm;
in fact applying the skein relation $n-1$ times to a $n$ component
link will generate  the one-component knots which we are after.
These knots classify our Feynman graph 
and determine their UV divergences.

\begin{figure}[ht]
\epsfysize=7cm \epsfbox{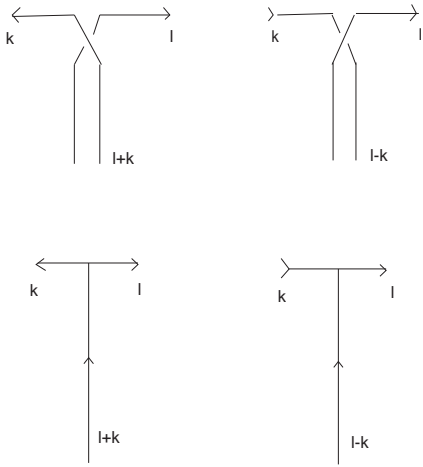}
\caption[From a vertex to a crossing.]{\label{hab310}\small
The replacement of a three-point vertex by an 
overcrossing. When we reverse the orientation of
lines at the vertex, we have also to reverse the orientation of lines
in the link diagram, and, accordingly,
exchange the over- to an undercrossing. All cases follow from
strict obedience to a "traffic rule": the momentum from the right
is overcrossing.}
\end{figure}
\begin{figure}[ht] 
\epsfxsize=5in \epsfbox{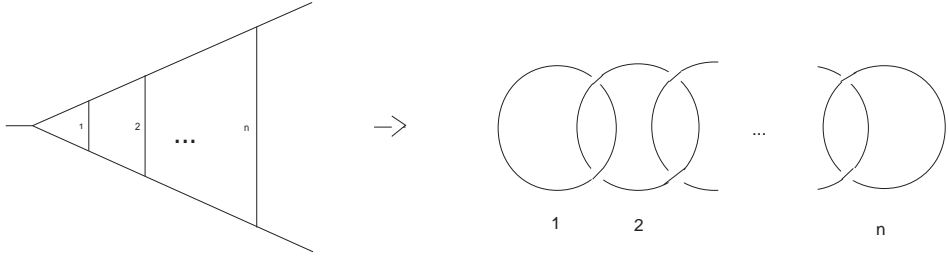}
\caption[From a Feynman graph to a link diagram.]{\label{hab311}\small 
The translation from a Feynman diagram to a link diagram.
Each vertex is replaced by an over/undercrossing according
to the momentum flow at the vertex. We follow the convention
to have the momentum flow in each loop counterclockwise.
Here we used a momentum routing so that each propagator
$P_i, i=1,\ldots,n$,
appearing as a rung in the Feynman graph
above, carries loop momentum $l_i-l_{i+1}$.}
\end{figure}

For a start, let us introduce the following two 
rules to map any Feynman diagram into a link 
diagram, as shown in Fig.(\ref{hab311}).
\begin{itemize}
\item
Every loop in the Feynman diagram corresponds to a link.
Correspondingly, a $n$-loop diagram will map to a link diagram consisting of
$n$ components.
\item
The links are oriented according to the flow of loop momenta, 
and follow the rule
that at every vertex the momentum coming from the right is 
overcrossing as in Fig.(\ref{hab310}).
\end{itemize}
We ignore exterior momenta in this process. 
For the ladder topology
we also understand that each crossing in our link diagram should
correspond to a vertex in the Feynman graph.
Our nested loop structure results in  a sequence of (Hopf-) links
as in Fig.(\ref{hab311}). 


For the momentum routing as envisaged in Fig.(\ref{hab311})
we can now start moving these Hopf-links into different
positions, using Reidemeister moves. It is not too difficult
to see that this corresponds to other momentum flows, still
in accordance with our rules given above.
We then use the skein relation, and disentangle the
link diagram until we end with a collection of disjoint
unknots, where the unknots may have non-trivial curl, see below.

The states which are generated in this manner we identify
with the terms demanded by renormalization theory. 
We will survey in a later section results which concern
Feynman diagrams of a much more complicated topology.

The empirical results to be discussed there justify the
following remarks on the way how to assign link diagrams to 
Feynman diagrams, and how to incorporate the skein relation.

According to the momentum flow, we assign a link
diagram to a Feynman diagram. The demand that
every vertex in the Feynman diagram should correspond to a crossing 
in the link diagram sets a lower boundary for the number of
crossings in the link diagram. Further crossings will be
generated solely by the topological complexity of the Feynman diagram.
Such topological complexity is absent from the simple ladder topologies
considered here, but will be present in the examples discussed
in the chapter on knots and transcendentals.

These extra crossings in the link diagram shall be of the least number
possible in accordance with the two rules above.
This is a non-trivial restriction on the allowed Reidemeister moves,
which we hope to set in a more axiomatic context in future
work. So far, we embodied it in our investigations
by regarding the Feynman diagram as built up
from a ladder diagram: for every Feynman diagram
it is true that when we remove a sufficient number of propagators,
we obtain a simple (ladder) topology.
Reinserting the other propagators step by step amounts to 
adding link component after link component. The basic ladder
topology is assumed to correspond to the Hopf link  as
described here. The extra components are then added so as to 
minimize the number of extra crossings, that is we add the next
propagator which specifies the location of two more crossings in the
diagram, and then try to connect these in the way which avoids as
many crossings as possible. The reader will find it instructive
to study Figs.(50,51,53,56) in the light of these considerations.

Once the set of link diagrams assigned to a Feynman diagram
is found, we will apply the skein relation to those crossings which 
correspond to proper vertices in the diagram.
For a $n$-loop Feynman diagram we apply the skein relation
$n-1$ times.

We do so as we are particularly interested in the one-component 
link diagrams, {\em knots}, which one obtains from
the link diagram during this procedure. They cannot 
appear before the skein relation is applied $n-1$ times.
For the simple examples here, one always finds at this stage
unknots, which, as we will see, is reflected in number theoretic
properties of the counterterms (absence of knots = absence
of transcendentals).

In general, we stop skeining at this stage and find that the knots
resulting  after skeining $n-1$ times faithfully predict
the transcendentals encountered in the calculation.
Once we have applied our skein relation and disentangled
the link diagram as described here, we do allow for 
Reidemeister moves in the generated knots.

Recent results \cite{BK15,4TR} point towards a connection 
of this process with chord diagrams, and indicate that the
whole procedure can be understood in the context of a four-term
relation.

\begin{figure}[ht] 
\epsfysize=3cm \epsfbox{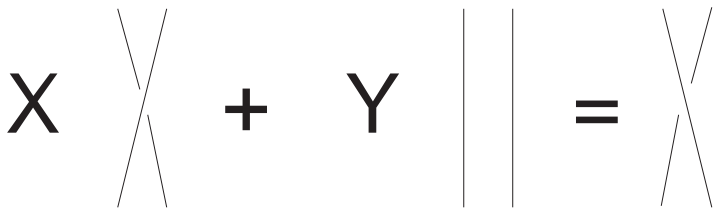}
\caption[The skein relation.]{\label{hab312}\small
The skein relation, an exchange identity which allows the
disentangling of the link diagram. $X$ and $Y$ have to be regarded
as operators to be identified with $A$ and $B$ in an appropriate 
manner. As usual, we assume that the operations appear locally
on a specified crossing in the link diagram.}
\end{figure}

We now continue to presume that some sort of braid structure underlies 
the algebraic systematics of renormalization theory. So we assume
that we can establish a skein relation of the form as indicated
in Fig.(\ref{hab312}). 
This assumption is at the moment based only on the vague evidence that
the ladder topologies renormalize according to Eq.(\ref{e23}), which, as we
will see below, fits into the pattern suggested by a skein relation.
In the course of the following sections we will look for further evidence to
justify this assumption.

\begin{figure}[ht] 
\epsfysize=7cm \epsfbox{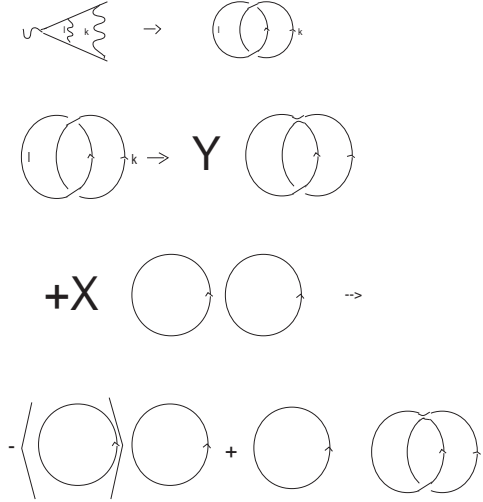}
\caption[A two-loop link.]{\label{hab313}\small 
An explicit  two-loop example.
The first line maps a Feynman diagram to a link diagram.
This link diagrams then, under the skein operation, gives two
terms. We compare these two terms with the terms which we
need to renormalize the diagram. We identify the
$Y$ term of the skein relation with the $B$ operator,
and the $X$ term with the removal of subdivergences as
provided by the operator $A$.}
\end{figure}

Let us consider Fig.(\ref{hab313}), a two-link diagram
and its disentanglement. It gives us a clue how to interprete
the coefficients $X$ and $Y$ in the skein relation. 
The three-link diagram, 
corresponding to $\Gamma^{(3)}$, 
\beas
\cil\crs\crs\cir\;\rightarrow\\
-<\cic>\cil\crs\cir\;\;+\;\cic\;\;\cil\wrt\crs\cir\;\rightarrow\\
+<<\cic>\cic>\cic\;\;-\;<\cic\;\;\cil\wrt\cir\;>\;\;\cic\\
-<\cic>\;\;\cic\;\;\cil\wrt\cir\;\;+\;\;\cic\;\;\cil\wrt\cir\;\;
\cil\wrt\wrt\cir
\eeas
confirms this in this graphical calculation for the three-loop case.
We already identified the operators $X,Y$ and
we obtain the correct result.
We conclude that the identification 
{\Large \beas
X^{r-1}(\cic\;\;\cic\;\ldots\;\cic\;\;) & \Rightarrow & [-A]^{r-1}(\Delta),\nonumber\\
Y^{r-1}(\;\cil\wrt\wrt\ldots\wrt\cir_{r-1}\;)
 & \Rightarrow & B^{r-1}(\Delta),
\eeas}
delivers our previous results.
We identify the unknot with an appropriate one-loop function ${}_0\!\Delta$
and links of the form 
\beas
\cil\wrt\wrt\ldots\wrt\cir_j\;
\mbox{\small A link with writhe number $j$.}
\eeas
with the corresponding function ${}_j\!\Delta$. 
Note that this implies that we have no invariance under
Reidemeister type $I$ moves. Accordingly,  we work with a regular isotopy.
Furthermore, the above result is independent of the routing of
momenta in the Feynman graph, as all possible routings
will generate
only pairwise concatenated rings, so that under the action
of the skein relation we would always obtain the same result.
In fact, we add to our rules the further demand that each crossing
in the link diagram should either correspond to a vertex in
the Feynman diagram, or must necessarily be generated by a non-trivial
topological structure of the diagram. 
In this respect ladder diagrams are trivial.


We conclude that the $n$-link diagram of the form of Fig.(\ref{hab311})
delivers, via the skein relation and appropriate
identification of the operators $A$ and $B$,
the $Z$-factor, $Z_1^{(n)}$, as drawn below.
\beas
\cil\crs\crs\crs\ldots\cir_n =\\
-<\cic\;>\;\;\cil\crs\crs\ldots\cir_{n-1}\\
+\cic\;\;\cil\wrt\crs\ldots\cir_{n-1} =\\
+<<\cic\;>\cic\;>\;\;\cil\crs\ldots\cir_{n-2}\\
-<\cic\;>\;\cil\wrt\crs\ldots\cir_{n-2}\\
-<\cic\;\;\cil\wrt\cir\;>\cil\crs\ldots\cir_{n-2}\\
+\cic\;\;\cil\wrt\cir\;\;\cil\wrt\wrt\crs\ldots\cir_{n-2}\\
\ldots\;
\mbox{\small The skein tree in general.}
\eeas

\begin{figure}[ht] \epsfysize=4cm
\epsfbox{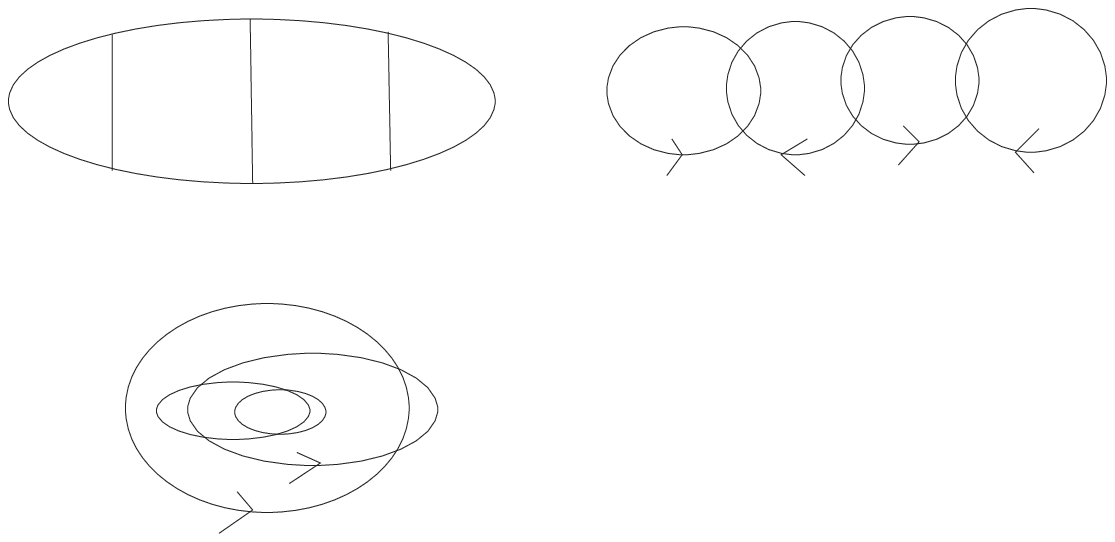}
\caption[A braid diagram with block structure.]{\label{hab315}
\small One way to obtain the braid of Eq.(\ref{brn4})
is to route momenta as in this diagram.}
\end{figure}

Let us briefly discuss how the above approach looks
from the point of view of a braid group approach.
Every closed link entangled in some other closed link generates a braid 
diagram with a very peculiar topology: 
\bea
\mbox{\small The braid diagram.} & & \sip\lir\lir\nonumber\\
\mbox{\small A closure of all strands
is always understood.} & & \sip\lir\lir\nonumber\\
\mbox{\small The braid would be generated in the case
$n=4$.} & & \lil\sip\lir\nonumber\\
 & & \lil\sip\lir\nonumber\\
& & \lil\lil\sip\nonumber\\
& & \lil\lil\sip\nonumber\\
& & =\sigma_1^2 \sigma_2^2 \sigma_3^2\label{brn4}
\eea
We infer such a braid diagram from our standard link diagram 
by the help of Fig.(\ref{hab315}).
Every crossing of a braid $i$ with a braid
 $(i+1)$ in Eq.(\ref{brn4})  corresponds
to the action of a braid group generator $\sigma_i$. 
We follow the convention to orient our loops consistently
so that all braid generators have positive powers. We will
see that empirical evidence supports this convention, so that
only positive knots determine the UV-divergent structure of a 
perturbative quantum field theory.


The $n$-link diagram Fig.(\ref{hab311}) we identify with $Z_1^{(n)}$
\bea
 \sigma_1^2\ldots \sigma_{n-1}^2 & \leftrightarrow & Z_1^{(n)},\nonumber\\
\sigma_i^2 & = & Y \sigma_i+X {\bf 1},
\label{e31}\eea
which implies the usual identification of the skein relation with the
Hecke
algebra relation $\sigma=X\sigma^{-1}+Y{\bf 1}$,
so that we would recover our previous expressions, e.g.
\bea
Z_1^{(3)} & \leftrightarrow & \sigma_1^2\sigma_2^2\Rightarrow\nonumber\\
& = & X^2-XY\sigma_1-YX\sigma_2+Y^2\sigma_1\sigma_2\Rightarrow\nonumber\\
& = & [-A+B]^2(\Delta).
\label{e32}\eea
Here we used a Hecke algebra representation of the braid group.
At this point the algebraic structure is in fact not
fully developed. Our $n$-link 
corresponds to words containing nothing else than products of the
form $\prod_i \sigma_i^2$, so that we can draw the braid diagram in
the simple block form of Eq.(\ref{brn4}).
Other routings of loop momenta map to expressions containing the
same braid generators in various orders, but reduce to the
same set of terms after skeining.
Due to this simple structure we do not, at this level,
see anything of the proper braid algebra
structure, encoded in relations as
\bea
\sigma_i\sigma_{i+1}\sigma_i=\sigma_{i+1}\sigma_i\sigma_{i+1}.
\label{e33}
\eea
\beas
\ldots\sip\lir\ldots & & \ldots\lil\sip\ldots\\
\ldots\lil\sip\ldots & & \ldots\sip\lir\ldots\\
\ldots\sip\lir\ldots & & \ldots\lil\sip\ldots\\
\mbox{\small The Reidemeister type $III$ move corresponding to 
Eq.(\ref{e33}).}
\eeas
We will not dwell on this point here but will have more to say
later on.
Now we establish similar algebraic structures in
more general circumstances.
\subsection{A generalization to ladder cables}
This generalization shows mainly that also more complicated vertex
graphs follow a behaviour which fits into the pattern established so far.
We do now allow for higher terms in the skeleton expansion of the vertex.
They correspond to more complicated link diagrams not
encountered yet. 
Non-trivial terms in the skeleton expansion
will have a non-trivial link structure, and the
question is how these links are entangled. Here
we do not investigate the link structure of each term,
but only pursue how renormalization theory  concatenates them.

\begin{figure}[ht] 
\epsfysize=3cm \epsfbox{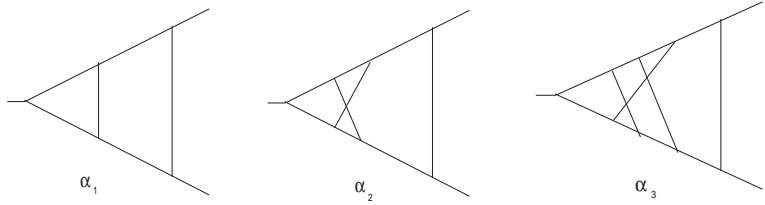}
\caption[The general ladder topology.]{\label{hab320}\small 
The general ladder topology. Every cable of lines
defining a new subdivergence is labelled according to its
topology. The examples given here have different topologies
$\alpha_i$ for their subdivergences. We replaced the
one-loop subdivergence by the first terms in the bare skeleton
expansion for the vertex, cf.~Fig.(\ref{hab33}).}
\end{figure}
\begin{figure}[ht] \epsfysize=4cm
 \epsfbox{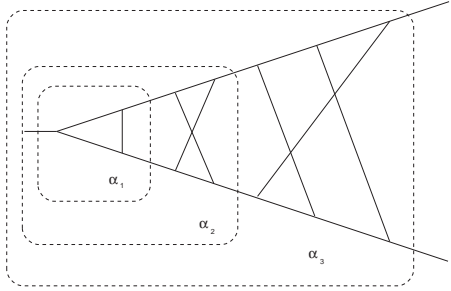}
\caption[Different sectors, different cables.]{\label{hab321}\small 
The cabling of loops into divergent sectors.}
\end{figure}

We denote different topologies by $\alpha_i$.
A Feynman graph of topology $\alpha_i$ has
$n(\alpha_i)$ loops.
For the time being we regard any  topology $\alpha_i$
as corresponding to a cable (this cable being a 
$n(\alpha_i)$-component link itself) 
and ignore the `fine structure' of these cables.
The sole question is how higher terms in the skeleton
expansion are concatenated from the viewpoint of renormalization
theory. Later we might have a conjecture how the corresponding
cables are knotted.
So let us now slightly generalize the ladder topology of Fig.(\ref{hab311}).
We want to include all sorts of vertex corrections which are themselves
free of subdivergences. This corresponds to a full skeleton expansion
of the vertex as in Fig.(\ref{hab320}). 
We label the different topologies by indices $\alpha_i$ but
still omit the dressing of internal vertices and propagators.
So the subscript $i$ labels different terms in the skeleton
expansion.
Note that every topology gives rise to a simple pole in $\epsilon$ only
\bea
\Delta^{\alpha_i}=\frac{c^{\alpha_i}}{\epsilon}+d^{\alpha_i},
\label{e34}\eea
due to the fact that it has no subdivergences by the very
definition of a skeleton expansion.
A Feynman graph as in Fig.(\ref{hab321}) 
will then result in a $Z$-factor contribution
\bea
Z & = & \prod_{i=1}^{3}[-A+B](\Delta^{\alpha_i}),
\label{e35}\eea
with the obvious definitions
\bea
\prod_{i=s}^{r} A(\Delta^{\alpha_i}) & := &
\Delta^{\alpha_r}
<\ldots <<\Delta^{\alpha_{s}}>\Delta^{\alpha_{s+1}}>
\ldots \Delta^{\alpha_{r-1}}>,\nonumber\\
\prod_{i=s}^r B(\Delta^{\alpha_i}) & := & 
\Delta^{\alpha_s} \;{}_{n(\alpha_s)}\!\Delta^{\alpha_{s+1}}
\ldots \;{}_{[n(\alpha_s)+\ldots+n(\alpha_{r-1})]}\!\Delta^{\alpha_{r}}.
\label{e36}\eea
Note that the concatenation in the $B$ operator takes the loop number of each
topology into account. That the same algebraic structure appears
results from the fact that we still have only nested divergences,
with only one maximal forest appearing, as indicated in Fig.(\ref{hab321}).


From the viewpoint of braids we would still map this situation to a 
$n$-cable link, so that we associate a cable (a collection of links)
to each $\Delta^{\alpha_i}$.
In braid generator language we would obtain generators 
$\sigma$, each acting on links in a different 
representation $\alpha_i$ so to 
speak, where in the `$B$' part of the skein relation we have 
a concatenation  of writhe numbers as shown above
and taking the various $n(\alpha_i)$'s into account.
As an example of this Fig.(\ref{hab321}) delivers
\bea
Z_1^{(3)}(\alpha_1,\alpha_2,\alpha_3)=
<\Delta^{\alpha_1}{}\;{}_{n(\alpha_1)}\!\Delta^{\alpha_2}{}\;{}_{
n(\alpha_1)+n(\alpha_2)}\!\Delta^{\alpha_3}\nonumber\\
-(<\Delta^{\alpha_1}{}\;{}_{n(\alpha_1)}\!\Delta^{\alpha_2}>
-
<<\Delta^{\alpha_1}>\Delta^{\alpha_2}>)\Delta^{\alpha_3}\nonumber\\
-<\Delta^{\alpha_1}>\Delta^{\alpha_2}{}\;{}_{n(\alpha_2)}
\!\Delta^{\alpha_3}>.
\label{e37}\eea
The $\Delta^{\alpha_i}$ are implicitly 
defined
as the dimensionless function of the
regularization parameter $\epsilon$ which multiplies the scaling
$(r^2)^{-\epsilon n(\alpha_i)}$ of the corresponding $n(\alpha_i)$-loop
Green function.


Note that our presentation here was done under the assumption that 
the $\Delta^{\alpha_i}$ are known. This is true for all possible
topologies up to the four-loop level, thanks to the major progress
in massless two-point functions, obtained by various authors 
\cite{chet}. A summary of the situation in this area can be found
in \cite{smirnov}.
Later we will see that our knot-theoretic approach suggests a
way to obtain these functions to all loop orders for all
topologies by associating them with various knots.


To get the final $Z$-factor we just have to add the results for the various
topologies.
Define the total loop order $n_t$ by 
\bea
n_t=\sum_i n(\alpha_i).
\label{e39}\eea
We have for the $Z$-factor generated from all topologies contributing
in a given loop order $m$
\bea
{\bf Z}_1^n = 1-\sum_{m=1}^n \sum_{I_m} Z_1^{(n_t=m)}(I_m),
\label{e40}\eea 
where $\{I_m\}$ denotes a complete 
set of topologies such that $n_t=m$, and the sum over
$I_m$ runs over all these topologies. Note also that with different 
topologies the
number of different types of graphs per topology proliferates and has to be
taken into account by an appropriate choice of basic functions
$\Delta^{\alpha_i}$.
We stress that at this stage we have not specified the actual link
diagrams associated to the different cables, nor the 
entanglement of the cables with each other. This is
already a very difficult problem, as it demands the concatenation
of topologically complicated terms in the skeleton, sitting
as subdivergences in other such terms. We will learn soon about
the complexity of such a problem, but still have to continue
our investigations of simple topologies.
This finishes our considerations of the vertex ladder topology and we now turn
to two-point functions.

\clearpage
\section{Corrections at the Propagator: planar, nested}

\begin{figure}[ht] 
\epsfysize=3cm \epsfbox{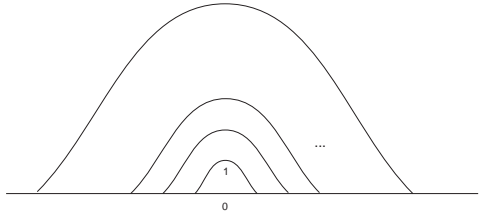}
\caption[The rainbow topology.]{\label{hab41}\small 
We call this the rainbow topology. For reference,
we gave numbers to the various propagators.}
\end{figure}

Again we consider only nested divergences.
The whole situation is largely a repetition of the situation
in the previous section. Nevertheless, we have to comment on some new features.
Consider the Feynman  graph in Fig.(\ref{hab41}).
Note that when we insert a zero-momentum coupling 
we would recover the situation of the previous section. 
In the spirit of it we would like to define
an expression of the type $[-A+B]^n(\Omega^1)$. We expect 
$\Omega^1$ to be the one-loop massless two-point function,
and would also need to define ${}_j\!\Omega^1$, its version with
increased writhe number. In fact, this is indeed correct. 
But this is only so when we restrict ourselves to two-point 
functions which are at most linearly divergent.
Otherwise single subtractions which we have considered so far
might not be sufficient
to subtract subdivergences. We will briefly comment on this difficulty
below and give some more details in an appendix.
But let us consider Fig.(\ref{hab41}) for linear divergent cases
and denote the graph by
$\Omega^{(n)}$.
Define 
\bea
\tilde{\Omega}^{(n)}:=\Omega^{(n)}\downarrow_{m_0=m_1=0}.
\label{e41}\eea
The notation $\ldots\downarrow_{m_0=m_1=0}$ means the following.
Consider the Feynman graph in the form where all propagators have
the usual quadratic denominators, and spin structures determine
the numerator expression. Then we set $m_0=m_1=0$ in the
first two propagators in the denominator.
We nullify the masses of the inner loop.
One easily shows by powercounting that 
\bea
<\Omega^{(n)}-\tilde{\Omega}^{(n)}>=0,
\label{e42}\eea
and it is sufficient to investigate $\tilde{\Omega}^{(n)}$. We add its 
counterterm
graphs. Next we want to set all masses to zero as the remaining
overall divergence is mass independent. 
For a massive theory  we typically have to take into account
one further subtlety: 
that there are linear and logarithmic divergent terms to be
considered (distributed over two formfactors
usually) for the mass and wave-function
renormalization. This reflects the fact that the 
two-point function $<0|T[\phi_{in}(x_1)\phi_{in}(x_2)]|0>$
corresponds to two monomials in the Lagrangian usually,
the mass term $-m^2\phi^2/2$ and the kinetic term $(\partial_\mu\phi)^2/2$.
Note that for a massless theory one directly obtains 
formulas equivalent to the results in the previous section.
For massive particles $\Omega$ is  a sum
of two terms, one providing the mass renormalization,
the other one the so-called wave function renormalization.
Together they achieve that the propagator-function  and its derivative
with respect to the exterior momentum are finite.\footnote{Higher order
derivatives are finite anyhow. Upon Taylor expanding one verifies
that for a Green function of overall degree of divergence $\omega$
the first $\omega+1$ terms in its Taylor expansion in an
external momentum are divergent. Accordingly, for a logarithmic divergent
vertex function only the first -zeroth- order term
is divergent, while a linear divergent propagator function also has
a divergent derivative. This can be confirmed by applying the derivative
to the integrand and subsequent powercounting.}
For example the fermion propagator
delivers
\begin{equation}
\Omega= a(q^2,m^2) \qslash + m\; b(q^2,m^2) {\bf 1}.
\end{equation}
 
 
Therefore this problem resembles the
one considered in the appendix, where the case of various
form factors is studied. There the reader will find a worked out
procedure for the general case.

\begin{figure}[ht] 
\epsfxsize=8cm \epsfbox{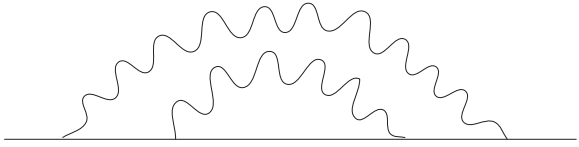}
\caption[Mass and wave function
renormalization.]{\label{hab42}\small 
An example for a mass and wave function 
renormalization.}\end{figure}

Here we only describe this general procedure as follows:
In the numerator, keep all terms which are overall linearly or
logarithmically divergent. That is, consider the numerator as a polynomial
in masses and the exterior momentum. Abandon all terms which are of degree
two or higher, as they provide only finite contributions. 
The terms which are now linear in masses  are the terms contributing 
to the mass renormalization, the other terms are necessarily linear in
the exterior momentum at the end, so they will give the wave-function
renormalization.
As an example let us consider massive QED at the two-loop order,
Fig.(\ref{hab42}).
The integration of the inner loop gives\footnote{
Depending on the chosen renormalization conditions
 these $Z$-factors are not yet
the standard mass and wave function renormalization, which can be inferred 
from a linear combinations of our $Z$-factors.} 


\bea
\Omega^{(1)} & = & (\Omega\kslash+ m \Omega_m{\bf 1})
(k^2)^{-\varepsilon},\nonumber\\
Z_2^{(1)} & = & <{\Omega}\kslash+ m {\Omega}_m{\bf 1}>
=\kslash<\Omega>+m<\Omega_m>{\bf 1}\nonumber\\
 & = :& Z_{2,w}^{(1)}  \kslash+Z_{2,m}^{(1)} m {\bf 1}.
\label{e43}\eea


Adding the two-loop graph and its counterterm $Z_{2,w}^{(1)}\kslash 
+Z_{2,m}^{(1)}m{\bf 1}$ gives then
\bea
Z_2^{(2)}: &  & 
[<\Omega\;{}_1\!\Omega>-<<\tilde{\Omega}>\Omega>]\kslash,\nonumber\\
  & + & [<2\Omega\;{}_1\Omega_m+\Omega_m\;{}_1\Omega_m>\nonumber\\
  & &  - <<2\tilde{\Omega}+\tilde{\Omega}>\Omega>]m{\bf 1}.
\label{e44}\eea
The generalization to higher loop orders is obvious and incorporates
that the concatenation of terms generated by the $B$ operator will
mix the different one-loop functions, as discussed in the appendix.


So far we assumed that the overall degree of divergence was not
worse than linear. The case of quadratic divergences can be usually
handled by appropriate differentiation with respect to the exterior 
momentum or by further subtractions. 
In practical applications one also often uses gauge symmetries, which tend
to improve the quadratic degree of divergence for a vector boson
propagator to a logarithmic one. Nevertheless we sketch how to treat
a two-loop calculation involving a quadratic subdivergence in an appendix.
Typically, quadratic divergences involve an integral of the form
\begin{equation}
\int d^Dl \frac{(l^2)^{-\varepsilon}}{l^2-m^2}=
\int d^Dl \frac{m^2(l^2)^{-\varepsilon}}{(l^2-m^2)l^2},
\end{equation}
where we used a typical identity (Eq.[\ref{dr3}])
in dimensional regularization
to transfer a quadratic divergence to a logarithmic one.
Note that at no stage do we encounter new infrared singularities
as we carefully avoid oversubtractions in this approach \cite{mpla}.


Summarizing we note that we observe in this section similar algebraic
structures as in the previous one. All the Feynman graphs
considered here can be obtained from the vertex graphs of the previous 
section by deletion of an external propagator. 
Consequently the results
as far as link diagrams are concerned remain unchanged.
For the rainbow topologies considered here the counterterm
contributions match the skein algebra for concatenated Hopf links
in the same way as before.

\clearpage
\section{Simple Entanglements}
So far we have considered two cases, ladder and rainbow topologies.
In both cases we restricted ourselves to a very simple
situation as far as renormalization is concerned: we only
had one maximal forest, with subdivergences strictly
nested into each other, and all divergences localized at the
same vertex (or propagator, for the case of the previous section).
Both cases give link diagrams of a simple topology, chains of
pairwise concatenated links. The corresponding braid expressions were of
second degree in each braid group generator. Correspondingly,
after applying skein relations, so far we never met a knot.
The most complicated figure generated was the unknot with writhe number
$n-1$, for a $n$-loop graph.


In this section we like to study  more general
Feynman graphs. We still restrict ourselves
to Feynman graphs which are iterated from only the first
term in the skeleton expansion. We will have
ladder and rainbow topologies. But this time we allow
them to be localized at arbitrary different points
in the graph. 


These graphs are still calculable from our elementary generalized
one-loop functions. But we will see that when we localize 
subdivergences at different points, we will  be 
able to generate more complicated link diagrams. 


In a later section we obtain the following
result: As long as we have simple iterated ladder
(or rainbow) topologies localized at not more than two different points,
the pole terms in our Laurent series have purely rational
coefficients. We will learn to associate the appearance of
transcendentals with the appearance of knots in the 
link diagrams. Here we will address the question if 
such transcendentals appear when we localize simple subdivergences
at various different points, and if the link diagrams can
contain non-trivial structures.

\begin{figure}[ht]\epsfysize=4cm \epsfbox{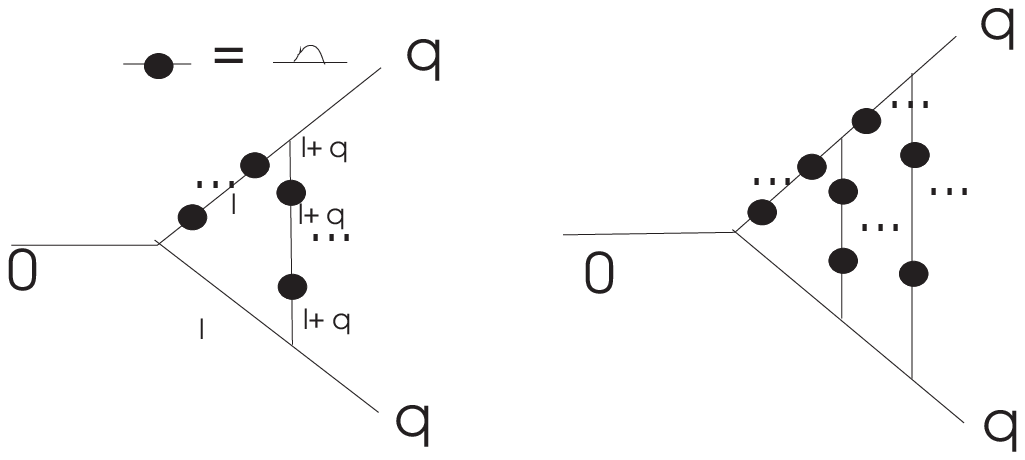}
\caption[Dressing a one-loop graph.]{\label{hab51}
\small The simplest possibility for divergences located at different
points appears with dressings of a one-loop skeleton graph.
The dressings, indicated by the black blobs, are
assumed to be with one-loop rainbow diagrams.
Generalizations as on the rhs are still calculable
in a finite set of generalized one-loop functions. Now that we
have disjoint  subdivergences new phenomena may occur.}
\end{figure}

Let us combine the results of the previous two sections. We still
consider vertex corrections to provide the basic skeleton 
graph and begin to dress internal
propagators 
with rainbow diagrams, Fig.(\ref{hab51}).
Let us start with the simplest possible example,
a one-loop graph dressed with one-loop bubbles at various points,
as in Fig.(\ref{hab51}) on the lhs.
\begin{figure}[ht] \epsfysize=2in \epsfbox{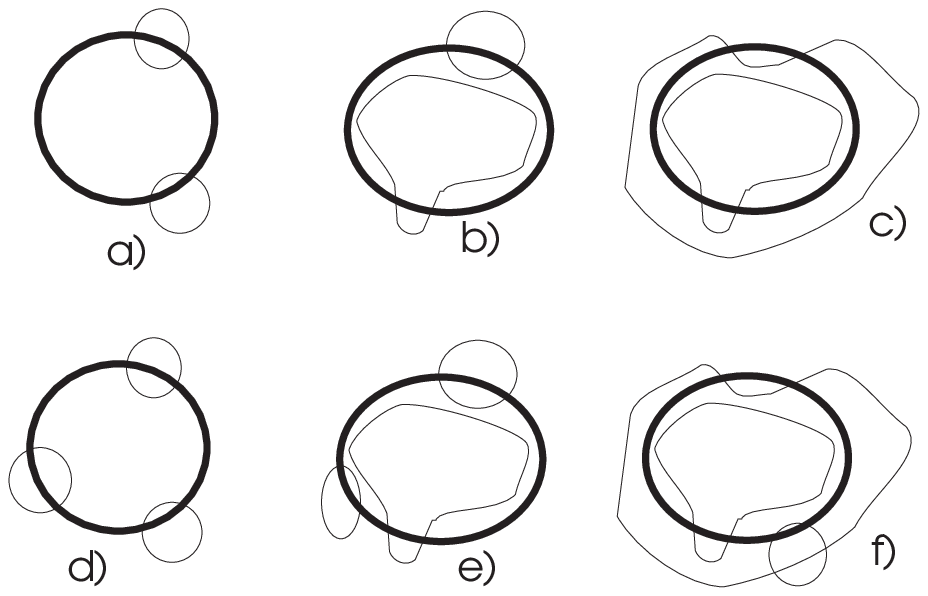}
\caption[Link diagrams for disjoint subdivergences.]{\label{hab52}
\small Link diagrams for disjoint subdivergences for a three
and four loop example. We follow the guiding
principle that all components which correspond to subdivergences
shall be concatenated with the appropriate skeleton, indicated
by a thick line. We indicate
different possible routings of momenta. For the three loop example
all routings are equivalent. In the four-loop case  
we obtain non-trivial entanglements for the first time.}
\end{figure}
This time we have disjoint subdivergences.
We can easily calculate the graph and its counterterm expressions.
At zero momentum transfer, the graph has two different
propagators. One of them carries loop momentum $l$ while the other one
carries momentum $l+q$.
Let us denote the graph with $i_1$ insertions at the lines
carrying momentum $l$,
and $i_2$ at the line carrying momentum $l+q$, by $G(i_1,i_2)$.
In Table(\ref{t1}) we give the results of a calculation of these
diagrams. We indicate the values for the pure
graph $G(i_1,i_2)$ as well as the results after addition of the counterterms,
denoted as $\bar{G}(i_1,i_2)$. 
{\footnotesize
\begin{table}
\[
\begin{array}{c|c|c}
(i_1,i_2) & {\bf G(i_1,i_2)} & {\bf \overline{G}(i_1,i_2)}\\ \hline
(2,0) & 
(-\frac{1}{2}\zeta(2)+\frac{3}{2}\gamma_E^2-8\gamma_E+\frac{44}{3})
\frac{1}{\epsilon}+(-\gamma_E+\frac{8}{3})\frac{1}{\epsilon^2}+
\frac{1}{3\epsilon^3} 
& \frac{1}{3\epsilon^3}
             -\frac{1}{3\epsilon^2}
             -\frac{1}{3\epsilon}\\ \hline
(1,1) &
(-\frac{1}{2}\zeta(2)+\frac{3}{2}\gamma_E^2-8\gamma_E+\frac{44}{3})
\frac{1}{\epsilon}+(-\gamma_E+\frac{8}{3})\frac{1}{\epsilon^2}+
\frac{1}{3\epsilon^3} 
 & \frac{1}{3\epsilon^3}
             -\frac{1}{3\epsilon^2}
             -\frac{1}{3\epsilon}\\ \hline
(3,0) & 


{(-\frac{83}{6}\zeta(3)+2\zeta(2)\gamma_E-\frac{11}{2}\zeta(2)
-\frac{8}{3}\gamma_E^3+22\gamma_E^2-79\gamma_E+\frac{475}{4})\frac{1}{
\epsilon} \atop
+(-\frac{1}{2}\zeta(2)+2\gamma_E^2-11\gamma_E+\frac{79}{4})\frac{1}{
\epsilon^2}
+(-\gamma_E+\frac{11}{4})\frac{1}{\epsilon^3}+\frac{1}{4\epsilon^4}}
 & -\frac{1}{4\epsilon^4}
             +\frac{1}{4\epsilon^3}
             +\frac{1}{4\epsilon^2}
             +(\frac{1}{4}-\frac{1}{2}\zeta(3))\frac{1}{\epsilon}\\ \hline
(2,1) & 


{(-\frac{59}{6}\zeta(3)+2\zeta(2)\gamma_E-\frac{11}{2}\zeta(2)
-\frac{8}{3}\gamma_E^3+22\gamma_E^2-79\gamma_E+\frac{475}{4})\frac{1}{
\epsilon} \atop
+(-\frac{1}{2}\zeta(2)+2\gamma_E^2-11\gamma_E+\frac{79}{4})\frac{1}{
\epsilon^2}
+(-\gamma_E+\frac{11}{4})\frac{1}{\epsilon^3}+\frac{1}{4\epsilon^4}}
& 
-\frac{1}{4\epsilon^4}
             +\frac{1}{4\epsilon^3}
             +\frac{1}{4\epsilon^2}
             +(\frac{1}{4}-\frac{1}{2}\zeta(3))\frac{1}{\epsilon}
\end{array}
\]
\caption[Results for one-loop bubble insertions.]{\label{t1}\small
Results for one-loop bubble insertions. We note that $G(3,0)\not=G(2,1)$,
while $\bar{G}(3,0)=\bar{G}(2,1)$. We also note that 
$\bar{G}(3,0)=\bar{G}(2,1)$ contains the transcendental $\zeta(3)$.}
\end{table}}
To obtain these results one only has
to consider combinations of Laurent series generated by
\begin{equation}
I(i_1,i_2):=\int d^Dk \frac{1}{[k^2]^{1-i_1\epsilon}
[(k+q)^2]^{1-i_2\epsilon}},
\end{equation}
in the case of Yukawa theory which was used here 
as an example.\footnote{Similar results can be obtained for all
renormalizable theories, and are discussed in general in \cite{newplb}.}
A few remarks are in order. First we observe that 
the graphs themselves depend on the way how the bubbles
are distributed over the graph,
for example $G(3,0)\not= G(2,1)$. This dependence vanishes when
one considers the graph together with its counterterm contributions,
$\bar{G}(3,0)=\bar{G}(2,1)$.
This is a general property. In \cite{newplb} the reader will find a proof
for this result to all orders.
These findings are in agreement with the association of link diagrams
to the graphs. The link diagram should correspond to the graph including
the counterterm expressions, according to our experience in the 
previous sections. But the link diagram cannot distinguish between
bubble insertion at various different
points, as demonstrated in Fig.(\ref{hab53}).
\begin{figure}[ht] \epsfysize=3cm \epsfbox{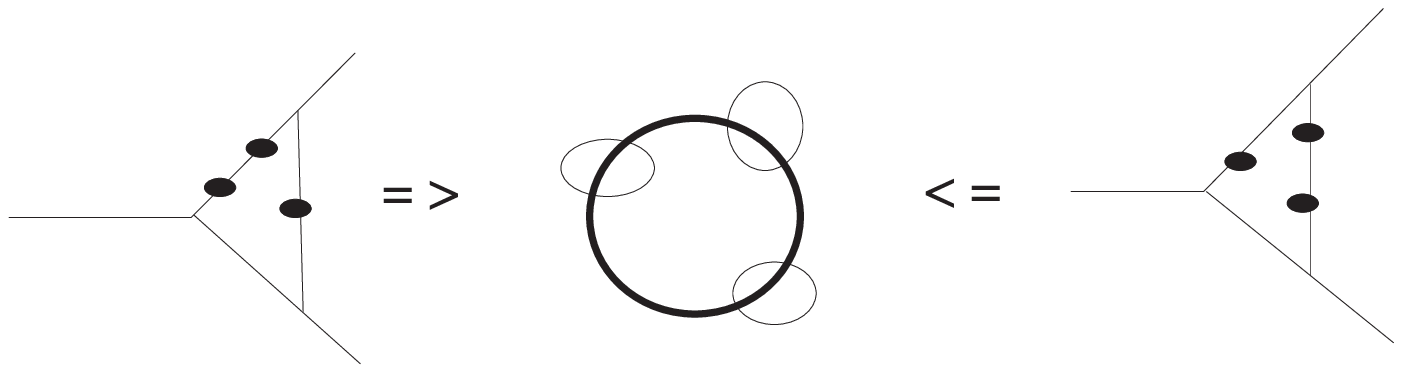}
\caption[Identities from link diagrams.]{\label{hab53}
\small Link diagrams demand certain identities between diagrams.
We only give some very basic examples, while such identities can be
verified in a much more general context.}
\end{figure}
As a consequence, we should have
identities for the function $\bar{G}(i_1,i_2)$, which is the
function $G(i_1,i_2)$ with subtracted subdivergences:
\begin{equation}
\bar{G}(i_1,i_2)=\bar{G}(i_1^\prime,i_2^\prime),\;\; 
\mbox{for}\;i_1+i_2=i_1^\prime+i_2^\prime.
\end{equation}
We stress that these identities are not fulfilled for $G(i_1,i_2)$
while $\bar{G}(i_1,i_2)$ abides by them, as it should.
In \cite{newplb} it is shown that such identities are in fact fulfilled for
a much larger class of diagrams. 
It occurs that  any Feynman diagram is in accordance
with the demands of link theory.
This very much supports our ideas of mapping Feynman diagrams
to link diagrams. In particular, field theory itself does not 
restrict diagrams to obey such identities, so that this 
reflects a profound connection between renormalization
theory and link theory indeed.

\begin{figure}[ht]\epsfysize=2cm \epsfbox{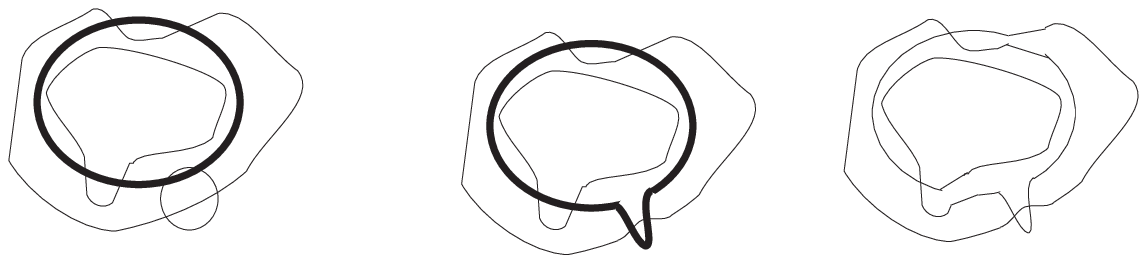}
\caption[The trefoil knot from disjoint subdivergences.]{\label{hab54}
\small Explaining the appearance of
the trefoil knot.}
\end{figure}
\begin{figure}[ht]\epsfysize=4cm \epsfbox{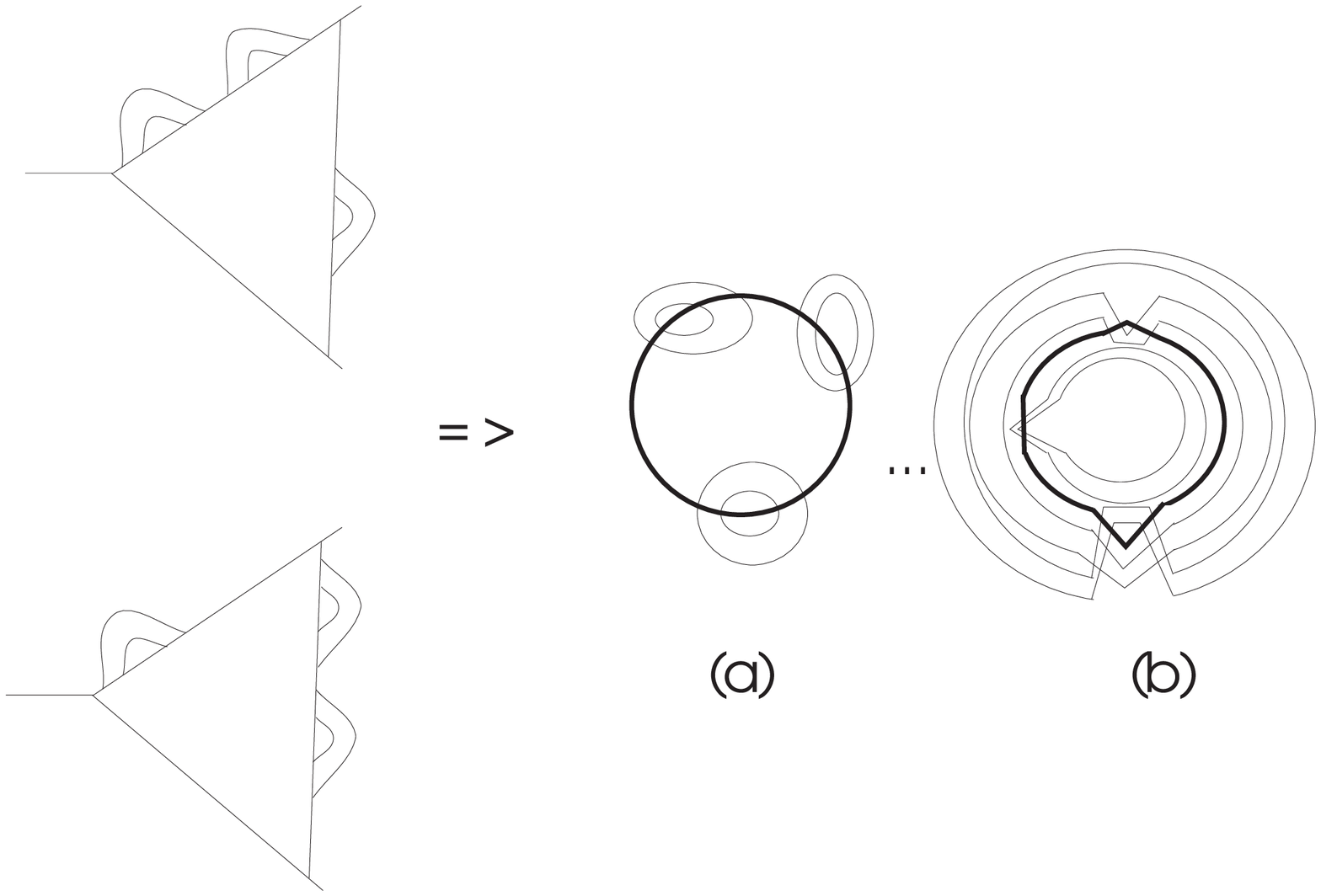}
\caption[Examples for disjoint
entanglements.]{\label{hab56}\small
The Feynman graph  provide two-loop rainbow dressings
of the basic one-loop skeleton. They both generate the same
set of link diagrams. From the analysis of these link
diagrams we expect rational numbers, provided by the knot free
link diagram (a), as well as transcendentals, provided by the
much more knotted link diagram (b).}
\end{figure}

Another fascinating observation is contained in 
Table(\ref{t1}). At the four loop level we have link diagrams
which have a simple structure like the link diagrams in previous sections,
but also we have a contribution where a more entangled diagram
occurs, see Fig.(\ref{hab52}f).
On the other hand, in Table(\ref{t1}) we see at the four loop level
the transcendental $\zeta(3)$ appearing. These two
observations are not unrelated. In anticipation of results we will
report on in a later section we assign the transcendental $\zeta(3)$
to the trefoil knot, while purely rational counterterms
indicate unknottedness of the corresponding link diagrams.
We find this trefoil knot in the link diagram of Fig.(\ref{hab52}f).
This is demonstrated in Fig.(\ref{hab54}).
But then in Fig.(\ref{hab52}) we see that at the four loop
level also simpler link diagrams appear, Figs.(\ref{hab52}d,e), 
and thus we expect
the total contribution to obtain a term involving $\zeta(3)$
as well as a term which is a pure rational number.
Table(\ref{t1}) demonstrates this indeed.
All these results are confirmed to much higher loop numbers in
\cite{newplb}. 
Let us also list  some explicit results for two-loop bubble insertions
in Table(\ref{t2}).
\begin{table}
\[
\begin{array}{c|c}
(i_1,i_2) & {\bf \overline{G}_2(i_1,i_2)}\\ \hline
(3,0) & 
{(-\frac{47}{336}\zeta(6)+\frac{11}{56}\zeta(5)+\frac{3}{140}\zeta(4)
+\frac{13}{70}\zeta(3)^2+\frac{9}{70}\zeta(3)+\frac{2}{7})\frac{1}{
\epsilon}
+(-\frac{11}{280}\zeta(5)-\frac{3}{28}\zeta(4)+\frac{1}{70}\zeta(3)
+\frac{1}{35})\frac{1}{\epsilon^2} \atop
+(\frac{3}{140}\zeta(4)-\frac{1}{14}\zeta(3)+\frac{1}{14})\frac{1}{
\epsilon^3}
+(\frac{1}{70}\zeta(3)-\frac{39}{280})\frac{1}{\epsilon^4}
+(\frac{47}{280})\frac{1}{\epsilon^5}
-(\frac{5}{56})\frac{1}{\epsilon^6}+\frac{1}{56\epsilon^7}}\\ 
\hline
(2,1) & 
{(-\frac{47}{336}\zeta(6)+\frac{11}{56}\zeta(5)+\frac{3}{140}\zeta(4)
+\frac{13}{70}\zeta(3)^2+\frac{9}{70}\zeta(3)+\frac{2}{7})\frac{1}{
\epsilon}
+(-\frac{11}{280}\zeta(5)-\frac{3}{28}\zeta(4)+\frac{1}{70}\zeta(3)
+\frac{1}{35})\frac{1}{\epsilon^2}\atop
+(\frac{3}{140}\zeta(4)-\frac{1}{14}\zeta(3)+\frac{1}{14})\frac{1}{
\epsilon^3}
+(\frac{1}{70}\zeta(3)-\frac{39}{280})\frac{1}{\epsilon^4}
+(\frac{47}{280})\frac{1}{\epsilon^5}
-(\frac{5}{56})\frac{1}{\epsilon^6}+\frac{1}{56\epsilon^7}}
\end{array}
\]
\caption[Two-loop bubble insertions.]{\label{t2}\small 
Two-loop bubble insertions follow the same pattern as described
before. Only the onset of knot-numbers is much more dramatic
as can be inferred from study of link diagrams, cf.~Fig.(\ref{hab56}).}
\end{table}
In Fig.(\ref{hab56}) we give some link diagrams for the examples considered 
in Table(\ref{t2}).
We find a consistent description, if we demand that all subdivergences
connect to their skeleton. For the two-loop bubble insertions
we stress that the inner bubble sees the outer one as its skeleton.
It is then this component which is entangled in the basic one-loop skeleton
drawn as a thick line. Such a concatenation of skeletons governs how
to concatenate our link diagrams in the presence of subdivergences.

\begin{figure}[ht]\epsfysize=4cm \epsfbox{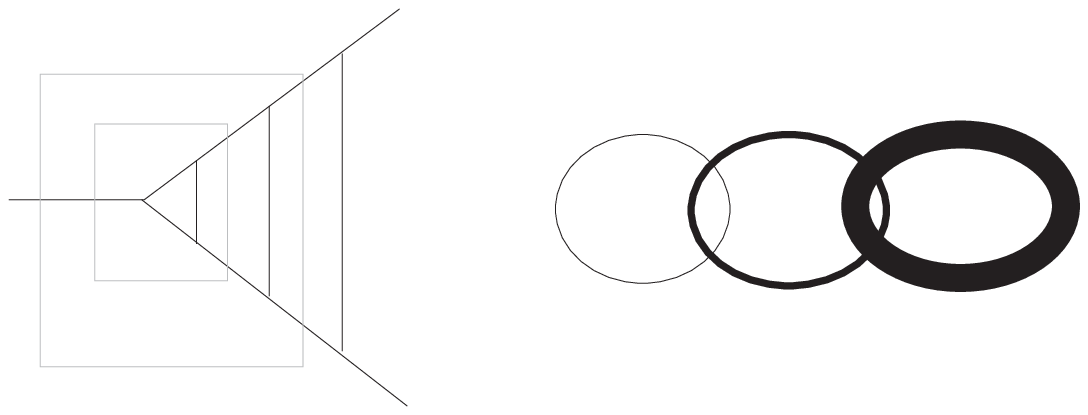}
\caption[Iterated skeletons.]{\label{hab57}\small
The cases of the previous section were an iterated set
of skeletons, determing the
entanglement of links. Thicker lines serve as the skeleton
for thinner ones.}
\end{figure}

All this is in agreement with our previous sections.
In the nested cases considered there
the $i$-th link had the $(i+1)$-th link as its skeleton
as indicated in Fig.(\ref{hab57}).
This singles out the standard assignment of concatenated
Hopf links as the most appropriate
link diagram. For nested divergences this results in a unique 
link diagram so that the component
$i$ connects to the component $i+1$ while here, for disjoint subdivergences, 
the demand to connect all subdivergences to the same component can
result in very different link diagrams, compare 
Figs.(\ref{hab52},\ref{hab56}).

\begin{figure}[ht]
\epsfysize=3cm \epsfbox{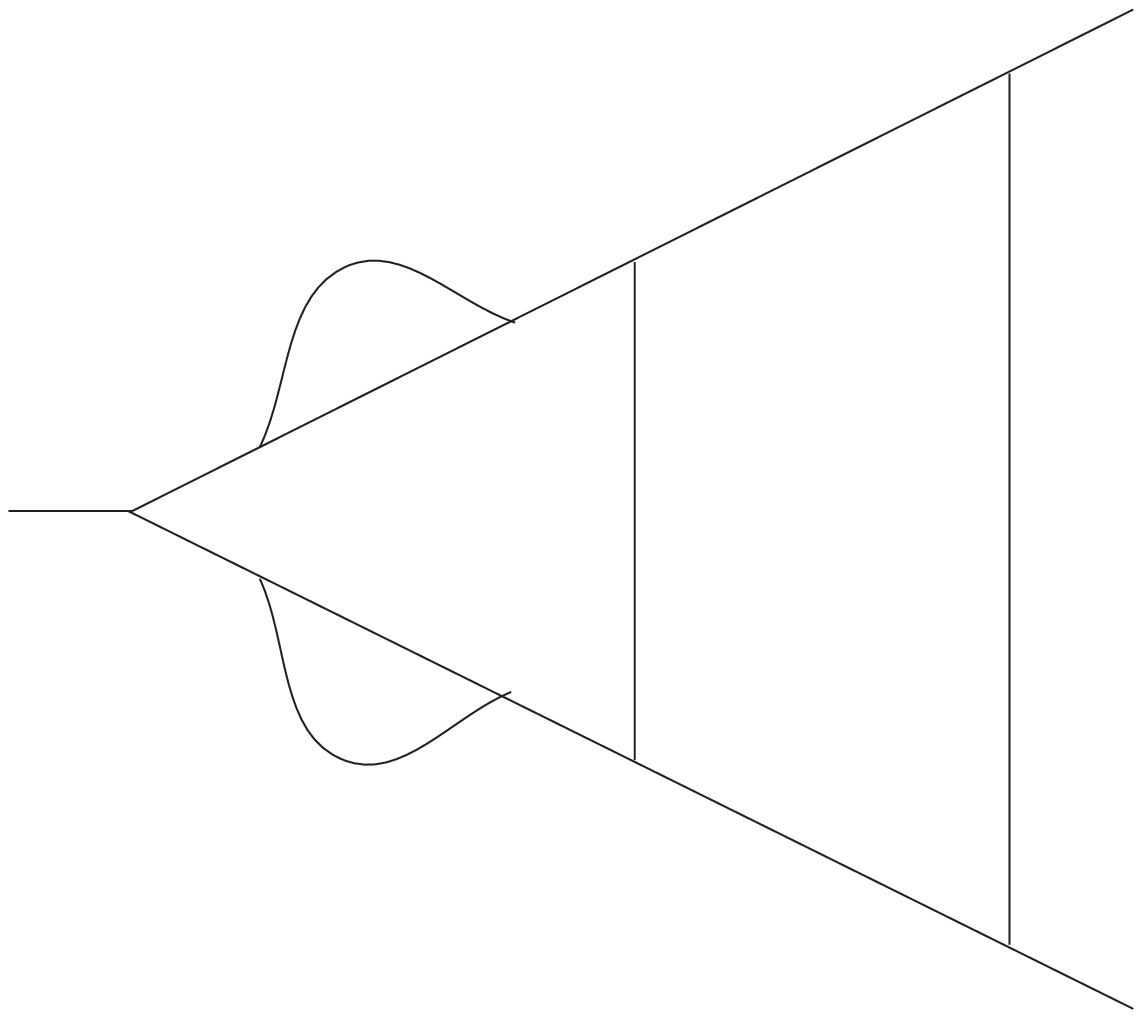}
\caption[Disjoint subdivergences.]{\label{hab59}\small
A dressing generating disjunct divergences.}\end{figure}

We now consider some examples to show how the differences between
nested and disjoint divergences look in terms of our
generalized one-loop functions.
We start with two disjoint subdivergences, Fig.(\ref{hab59}).
Note that the dressing does not influence our argument used to reduce the
considerations to massless functions at vanishing momentum transfer. This is so 
because we consider not only the graph but
the graph with all its counterterms. They compensate for the divergences
generated by the dressing, so that the arguments used in sections
three and four remain unchanged. Thus, our basic functions
$\Delta,\Omega$ encountered so far still serve as the appropriate
set to express all results.

\begin{figure}[ht] 
\epsfysize=3cm \epsfbox{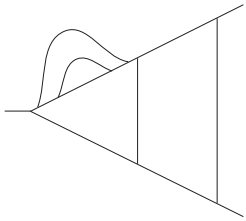}
\caption[A nested case.]{\label{hab510}\small  Compare this nested case with the previous
example.}\end{figure}

Back to our example above,
renormalization theory tells us that the correct answer is
\bea
\mbox{Fig.(\ref{hab59})} & \rightarrow &
\Omega \Omega\;{}_2\!\Delta\;{}_3\!\Delta -
 2 <\Omega>\Omega\;{}_1\!\Delta
\;{}_2\!\Delta\nonumber\\ 
& & + <\Omega><\Omega>\Delta\;{}_1\!\Delta
 -<\Omega\Omega\;{}_2\!\Delta\nonumber\\
  & & -2<\Omega>\Omega\;{}_1\!\Delta+
<\Omega><\Omega>\Delta>\Delta,
\label{e46}\eea
in our notation.
We see that the presence
of disjoint subdivergences is reflected by the presence
of expressions like $<\Omega><\Omega>$.
They do not come in concatenated forms  like
$<<\Omega>\Omega>$ or $\Omega\;{}_1\!\Omega$ generated by the $A$
and $B$ operators. These concatenations  were due to the total 
nested structure
of the subdivergences considered previously. 
Let us compare the above example
with Fig.(\ref{hab510}).
Combining our results on ladders and rainbows it delivers
\bea
\mbox{Fig.(\ref{hab510})} & \rightarrow &
\Omega \;{}_1\!\Omega\;{}_2\!\Delta\;{}_3\!\Delta -
  <\Omega>\Omega\;{}_1\!\Delta\;{}_2\!\Delta 
  \nonumber\\
   & & -<\Omega\;{}_1\!\Omega>\Delta\;{}_1\!\Delta
    + <\Omega<\Omega>>\Delta\;{}_1\!\Delta\nonumber\\
  & & -<\Omega\;{}_1\!\Omega\;{}_2\!\Delta 
 - <\Omega>\Omega\;{}_1\!\Delta\nonumber\\
 & & -<\Omega\;{}_1\!\Omega>\Delta+
<\Omega><\Omega>\Delta>\Delta,
\label{e47}\eea
which is obvious from sections three and four.
Both cases are free of transcendentals. All associated link diagrams
are free of knots.

\begin{figure}[ht] 
\epsfysize=8cm \epsfbox{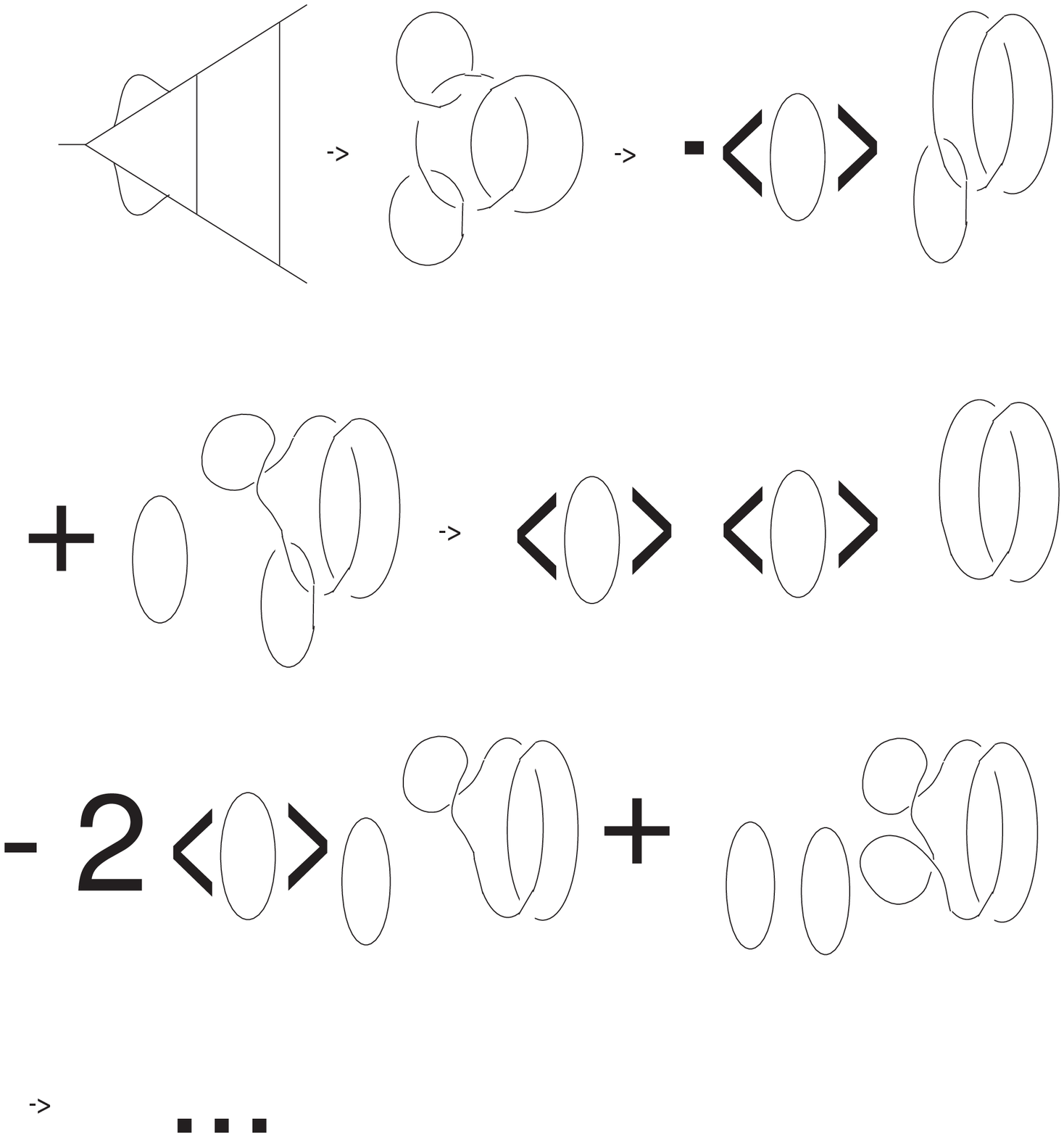}
\caption[An example calculation.]{\label{hab511}\small 
An example calculation using our graphical notation 
to demonstrate
the behaviour of $A$ and $B$ on disjunct subdivergences.
We indicated only one possible link diagram, as the others
produce equivalent results.}\end{figure}

In general we follow the rule 
that the $A$ and $B$ actions in the skein relation
concatenate the satellite links towards the skeleton link,
and we assume that opearators $A,B$ located
at different points act independently. 
We get
correct results, as in Fig.(\ref{hab511}).


We see the underlying principle: $A$ and $B$ treat disjunct
subgraphs as disjunct, that is they factorize
instead of concatenate the corresponding one-loop functions:
$\Omega\Omega$ instead of $\Omega{}_1\!\Omega$ and
$<\Omega><\Omega>$ instead of $<<\Omega>\Omega>$.
In \cite{newplb} examples up to ten loop were considered,
involving multiloop insertions at various places. 
The results always match the expectations of link theory,
fulfilling a large class of identities where, it seems,
no other explanation than an underlying connection to
link theory is available. 


This finishes our considerations of nested and disjunct divergences
and we now turn to overlapping divergences. Our hope is to find a similar
structure for topological simple graphs there. This is crucial for our
final
attempt to identify the topological nature
of a Feynman graph with certain properties of its divergent part.

\clearpage
\section{Corrections at the Propagator: planar, overlapping}
In this section we will address ourselves to overlapping
divergences. 
By now, we have some experience with nested and disjoint subdivergences.
In the next section we will see that ladder topologies like the
one considered in sections three and four deliver Laurent series
with rational coefficients in the proper pole part. 
This results from the special properties of
the $A,B$ operators of Eq.(\ref{e22}) in conjunction with the
structure of the generalize one-loop functions.
For this result to hold in general, we have to show that also
overlapping subdivergences allow for similar algebraic structures.
\subsection{The overlapping ladder}
\begin{figure}[ht] 
\epsfysize=2cm \epsfbox{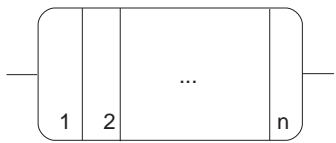}
\caption[The overlapping ladder.]{\label{hab61}\small  
The overlapping ladder. These link topologies are still
simple.}\end{figure}

\begin{figure}[ht] 
\epsfysize=2cm \epsfbox{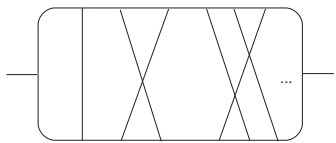}
\caption[The generalized case.]{\label{hab62}\small  The generalized case including various topologies for the
rungs.}\end{figure}

We will study graphs of the form as given in Fig.(\ref{hab61}).  
Further we briefly comment on their "cable generalizations" as in 
Fig.(\ref{hab62}).
So it is the main objective of this section to show that also the overlapping
divergences factorize in a manner similar to the cases studied so far.
This then allows to classify them by their corresponding link
algebra too. The main result of this section is that the 
overlapping ladder topology
gives indeed the same concatenations as before, 
with the only modification that
we have to sum over all possibilities to identify subdivergences in the
graph. Once more, this allows us in the next section to 
establish rationality also for simple topologies in the case of
overlapping divergences. This in turn defines transcendentality of
$Z$-factors as a sensible test for knots, that is non-trivial
topologies.

\begin{figure}[ht]
\epsfysize=3cm \epsfbox{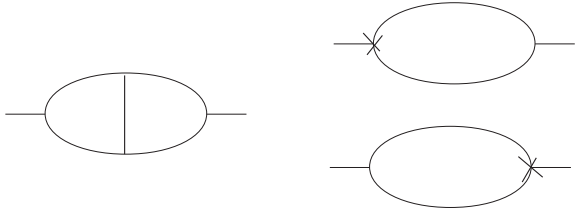}
\caption[The two-loop generic case.]{\label{hab63}\small  
The two-loop example. We have two subdivergences.}\end{figure}

Let us start with a simple two-loop example.
The graph has two overlapping subdivergences, so we have to calculate
the  expressions  indicated in Fig.(\ref{hab63}).
We have five propagators.
Let $N_i$  the product of all propagators involving
the loop momentum $l_i$. 
It is implicitly understood that there is some
numerator taking particle content and spin structures into
account. We operate in the following
only on the (scalar) denominator part
of the propagators under consideration,
without further mentioning the numerator structure. 
The equations given are literally
true for $\phi^3$ in six dimensions as an example,
but have obvious generalization to other renormalizable theories.
Nevertheless we remind the reader that the numerator structure affects
the power counting,
so that in the following one cannot use the exlicitly given 
denominator expressions for power counting purposes.


In the above example we
have $N_l=P_1 P_2 P_3$, $N_k=P_3 P_4 P_5$.
We define
\bea 
\tilde{N_i}:=N_i|_{q=0,m_j=0\;\forall j},\nonumber\\
{\bf N_i}:= \frac{\tilde{N_i}-N_i}{\tilde{N_i}}.
\label{e49}\eea
Then we have
\bea
\frac{1}{P_1 P_2 P_3 P_4 P_5} & = &
\frac{1}{P_1 P_2 P_3 P_4 P_5}{\bf N_l} {\bf N_k}
-\frac{1}{\tilde{N_l}P_4 P_5}-\frac{1}{P_1 P_2\tilde{N_k}}
+\frac{P_3}{\tilde{N_l}\tilde{N_k}},\nonumber
\eea
implying
\bea \Omega^{(2)} & \equiv & \Omega^{(2)}_f-\Omega^{(2)}|_l
-\Omega^{(2)}|_k+\Omega^{(2)}|_{l,k},
\label{e50}\eea
by construction.
Note that the  last term on the rhs vanishes in DR, according to
Eq.(\ref{dr3}). 


Let us now investigate Eq.(\ref{e50}). The first term on the rhs is
UV-convergent, as long as our overall degree of
divergence was not worse than linear; this is assumed to be the case.
In an appendix we comment on the general case.
In practice one has to be
cautious about the quadratic divergences of
vector-boson propagators. But they are  safe by gauge invariance, which 
reduces
the overall degree of divergence to  a logarithmic degree.


Back to our considerations of Eq.(\ref{e50}) we notice that 
the first term on the rhs will not contribute to our MS $Z$-factor,
due to its UV convergence:
\bea
<\frac{1}{P_1 P_2 P_3 P_4 P_5}{\bf N_l} {\bf N_k}>=0.
\label{e52}\eea
The remaining terms to be considered are
the second and third term on the rhs of Eq.(\ref{e50}).
The counterterms contribute
\bea
- Z_1^{(1)}\int d^Dl \frac{1}{P_4 P_5}
- Z_1^{(1)}\int d^Dl \frac{1}{P_1 P_2}.
\label{e53}\eea
Adding these contributions to the remaining two terms
on the rhs of Eq.(\ref{e50})
gives us a result free of subdivergences. As usual,
the remaining overall degree of divergence is independent of
masses in the process. We can set all
remaining masses to zero in this sum.
In the zero mass case we have $P_1=P_5,P_2=P_4$.
Our UV-divergences are thus contained  in
\bea
Z_2^{(2)}=2< \Delta
\int d^Dl \frac{(l^2)^{-\epsilon
}}{P_1|_{m_1=0}P_2|_{m_2=0}}\nonumber\\
-< \Delta>
\int d^Dl \frac{1}{P_1|_{m_1=0}P_2|_{m_2=0}} >.
\label{e54}\eea 
Defining
\bea
\;{}_j\!\Omega & := & \int d^Dl \frac{(l^2)^{-\epsilon j}}{
P_1|_{m_1=0}P_2|_{m_2=0}},\nonumber\\
\Omega & := & \;{}_0\!\Omega,
\label{e55}\eea
we have
\bea
Z_2^{(2)}=2< \Delta{}\;{}_1\!\Omega-< \Delta>\Omega >,
\label{e56}\eea
which has a striking similarity to our result Eq.(\ref{e18}).
The main difference is the factor of two which reflects the overlapping
structure in this simple example.
The apparent structure is typical for overlapping divergences.
We have two overlapping subdivergences. 
But we neither have a concatenation
like $<\Delta <\Delta>>$, as we had for nested subdivergences,
nor do we have a term $<\Delta>^2$, as it would characterize
disjoint subdivergences. 
Instead, we have a subtraction of subdivergences 
to the left and right. Due to the symmetry of the graph, this only 
generates
the factor of two in Eq.(\ref{e56}). Note that the two terms can be interpreted
as arising from two possibilities: each of both loops
can either be the skeleton or the subdivergence of the other one.
We sum over both possibilities.


In our above
result we used a somewhat
condensed notation. $Z_2^{(2)}$ stands for 
all the two-loop divergences of Fig.(\ref{hab63}).
It is a  notation for the  $Z$-factors corresponding
to a two-point function. As mentioned already usually there are two:
the mass, $Z_{2,m}$, and wave-function 
renormalization, $Z_{2,w}$. 
One would extract them by looking for their corresponding form factors,
e.g.~$\qslash$ and $m{\bf 1}$ for the fermion propagator, and by
Taylor expanding on the mass-shell. Again the reader interested in technical
details will find further comments in the appendix.
Graphically, our result is sketched below.
{\Large
\beas 
\left\{ \epsfysize=6mm \epsfbox{knot29.ps}\right\} =\\
{}_{{}_1\!\Omega}\!{\epsfysize=6mm \epsfbox{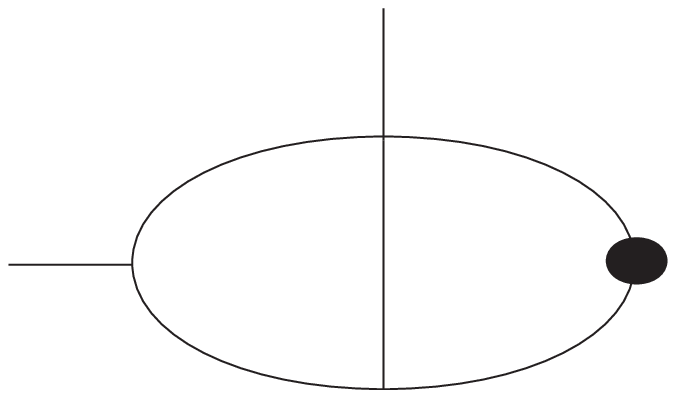}}_{\!\Delta}
+{}_{\!\Delta}\!{\epsfysize=6mm \epsfbox{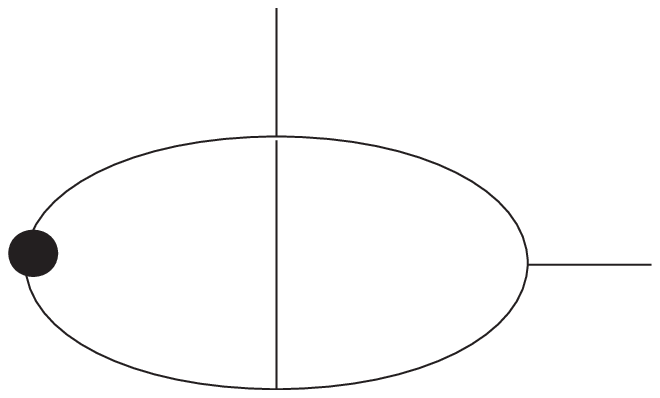}}_{{}_1\!\Omega}\\
-<\Delta>\epsfysize=6mm \epsfbox{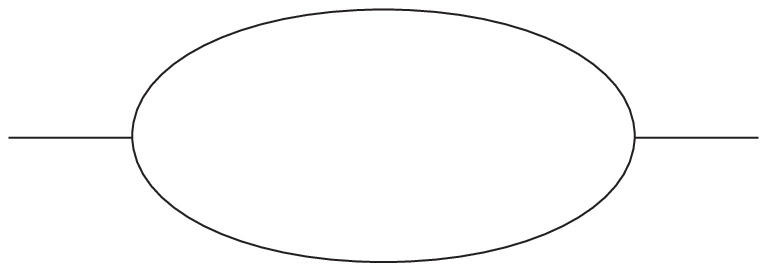}_\Omega
-\epsfysize=6mm \epsfbox{knot30s.ps}_\Omega\;<\Delta>.
\eeas
}
{\small Some notation emphasizing the fact that the result was just a
sum over all possibilities to interpret one loop as the skeleton 
$(\Omega)$ and the 
other one as the subdivergence $(\Delta)$.\\[5mm]}

\begin{figure}[ht] 
\epsfxsize=5in \epsfbox{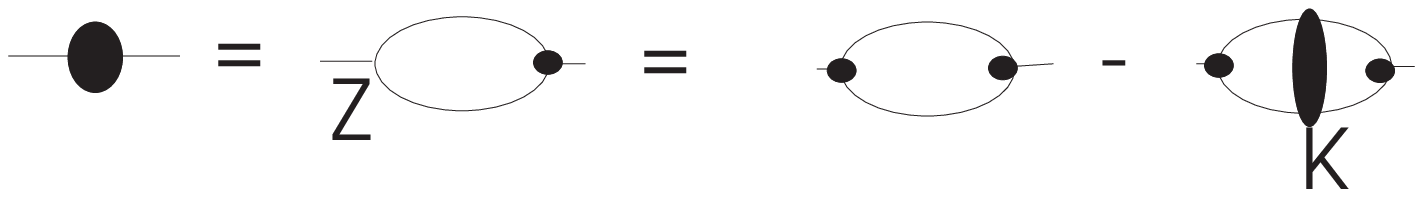}
\caption[The Dyson Schwinger equation for a propagator.]{
\label{hab64}\small The Dyson Schwinger equations for a propagator. 
In its second form we reexpressed the $Z$-factor of the vertex
with help of the Dyson Schwinger equation for this vertex.
Our subtractions in Eq.(\ref{e57})
on the lhs ($l_1$) and rhs ($l_n$) remove the divergences. So the
blobs have
a subtracted form, cf.~section 3.}\end{figure}

We want to generalize this approach to the $n$ loop case.
We are looking for a general prescription to convert overlapping
topologies back to products of one-loop graphs,
as we did before for the nested topologies. We then hope
to find a similar correspondence to braid structures as we found there.
Consider the $n$-loop ladder graph Fig.(\ref{hab61}).
We claim that 
\bea
\Omega^n_f & = & \Omega^n-\Omega^n|_{l_1}\-\Omega^n|_{l_n}+
\Omega^n|_{l_1,l_n}\nonumber\\
& = & \Omega^n N_{l_1}N_{l_n},
\label{e57}
\eea
is finite. The functions $\Omega^n|_{\ldots}$ used above are implicitly
defined by the second line in the above equation.
Note that $\Omega^n|_{l_1,l_n}$, the last term 
on the rhs in the first line, does not vanish if $n\not= 2$.
The assertion on $\Omega^n_f$ is justified by the fact that all possible 
divergent
sectors cancel in the above expression. This can also easily
seen by doing a powercounting for Eq.(\ref{e57}) or by investigating
the Dyson Schwinger equations of Fig.(\ref{hab64}).  
This is in agreement with the general result that an overlapping
divergence renormalizes by subtracting out all divergent subgraphs
in a manner different from disjunct or nested divergences \cite{Itz}.
Our two loop example above was generic in this respect.
While for disjoint cases the product of the operators removing the
subdivergences
will  appear, in the overlapping case this product structure is not
maintained. It is replaced by a sum over subtractions at all divergent
subgraphs.

\begin{figure}[ht] 
\epsfysize=6cm \epsfbox{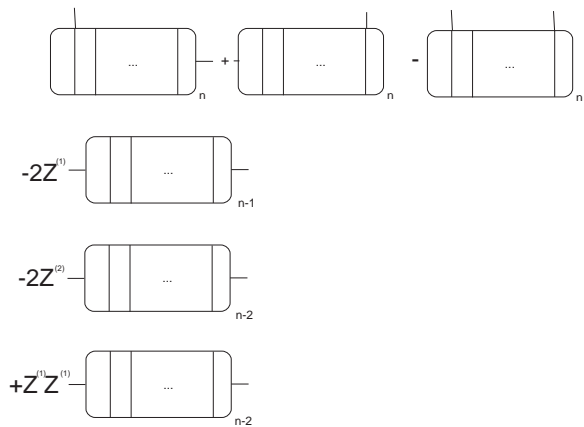}
\caption[Ladder with counterterms.]{\label{hab65}\small 
Our ladder with all its counterterms. In the
first line wee see the result of Eq.(\ref{e58}).
The next line gives a $n-1$ loop graph, multiplied by the
one-loop $Z$-factor. It appears twice, as we identify
such a one-loop subdivergence on the rhs as well as on the lhs of
the graph. The next lines incorporate the higher
contributions in the $Z$-factors.}\end{figure}
\begin{figure}[ht] 
\epsfysize=2cm \epsfbox{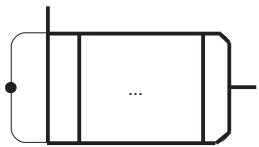}
\caption[Unfactorizable expressions.]{\label{hab66}\small 
These expressions still involve functions which do not
factorize into our basic $\Delta$'s and $\Omega$'s.
Typically, the expressions involve  functions to be evaluated
at non-vanishing momentum transfer, drawn in thick lines
above. The thin line furnishes a further loop
momentum flowing through the thick lines, and the exterior
momentum is still present.}\end{figure}
\begin{figure}[ht] 
\epsfysize=2cm \epsfbox{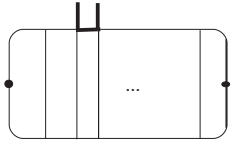}
\caption[1-states.]{\label{hab67}\small These `1-state' expressions are easy
to calculate. The thick line indicates the sole
propagator carrying the external momentum $q$.}\end{figure}

What now remains to be calculated are the last three terms on the rhs
of Eq.(\ref{e57}).
\bea
<\Omega^n_f>=0\Rightarrow <\Omega^n>=
<\Omega^n|_{1}+\Omega^n|_{n}-\Omega^n|_{1,n} >.
\label{e58}\eea
Let us add our counterterm expressions:
\bea
-2\sum_{i=1\ldots n-1}
Z_1^{(i)}\Omega^{(n-i)}\nonumber\\
+2\sum_{{i=1\ldots n-1 \atop {j=1\ldots n-1 \atop 1<i+j <n,i\not= j}}}
Z_1^{(i)}Z_1^{(j)}\Omega^{(n-i-j)}\nonumber\\
+\sum_{i=1\ldots n-1}(Z_1^{(i)})^2\Omega^{(n-2i)}.
\label{e59}\eea
The whole contribution then is given in the graphical notation
of Fig.(\ref{hab65}).  
This expression is
free of subdivergences and we can set all remaining masses to zero.
It is still not quite what we want as expressions like
those of Fig.(\ref{hab66}),  
involve loops which are not expressible in terms of our basic
functions $\Delta,\Omega$. Only expressions 
such as Fig.(\ref{hab67}) where the
exterior momentum flows just through one loop are amenable to our
procedure.  


The following constructive proof which shows
that we can express overlapping 
divergences in terms of concatenated one-loop functions is somewhat
technical. We nevertheless refrain from banning it to an appendix.
The fact that all graphs realizing a ladder topology 
follow the same pattern, even for overlapping divergences, 
is very important to us. We thus give the construction in some detail.


To see how to proceed let us recapitulate what we have done in
Eq.(\ref{e57}). We wrote the original function $\Omega^{(n)}$
as a difference between a function which was UV-convergent
$\Omega^{(n)}_f$ and some simpler functions $\Omega^{(n)}|_i$.
The guiding principle was Weinberg's theorem \cite{Weinberg}. 
In the form which is useful
for us it states that once an analytic expression
is finite by powercounting
for its overall degree of divergence {\em and} 
all its subdivergences, we are assured
that it is convergent. This allows us to conclude that
$<\Omega^{(n)}_f>=0$. But Weinberg's theorem also tells us that
an expression is convergent if it has a vanishing overall degree of 
divergence and either it has no subdivergences (the previous case) or 
it has all
its subdivergences subtracted by appropriate counterterms. 
This statement is just the underlying 
principle which determines the counterterms  
required to make a renormalizable theory finite. Vice versa, we know
that for every combination of loop momenta providing a subdivergence,
there exist an appropriate counterterm in the sum in Eq.(\ref{e59}).
We will show that this gives us an algorithm to continue the
process which was started in Eq.(\ref{e57}) to simplify $\Omega^{(n)}$.


But we will give the result first.
To this end let us introduce some notation.
\begin{itemize}
\item
The massless $n$-loop ladder graph: ${}^\mid \Box\ldots\Box^\mid_n$,
where ${}^\mid$ indicates the flow of the external momenta.
We call this a $(0,n)_n$ state, as the momenta flows through the
whole $n$-loop graph.
\item
The $(i,j)_n$ state:
$ \Box\ldots\Box_i^\mid\Box\ldots\Box_j^\mid\Box\ldots\Box_n$,
we agree to count $i$ from the left and $j$ from the right.
There are $n$-loops altogether. The external momentum flows
through $n-i-j$ of them.
\item
Counterterm graphs:
$\Box\ldots\Box_r] 
\Box\ldots\Box_i^\mid\Box\ldots\Box_j^\mid\Box\ldots\Box_{n-r-s}
[ \Box\ldots\Box_s$ denotes a $(n-r-s)$-loop graph in
a $(i,j)_{n-r-s}$ state multiplied by $Z^{(r)}Z^{(s)}$.
This is the same as the previous case, only that $r$ loops on the
left and $s$ loops on the right are replaced by the
corresponding $Z$ factor. 
Note that this is a condensed notation
for a whole $Z$-factor, e.g.: $\Box\Box] = 
(<\Box\Box-<\Box>\Box>)]$.
\end{itemize}


Our final claim is that the following expression, involving
1-states only, gives
the correct result for the $n$-loop $Z_2$-factor.
\bea
\sum_{i=0}^{n-1} \Box\ldots\Box_i^\mid\Box^\mid\Box\ldots\Box_n\nonumber\\
-\Box]\;\sum_{i=0}^{n-2} \Box\ldots\Box_i^\mid\Box^\mid\Box\ldots\Box_{n-1}\nonumber\\
-\sum_{i=0}^{n-2} \Box\ldots\Box_i^\mid\Box^\mid\Box\ldots\Box_{n-1}\;[\Box\nonumber\\
-\Box\Box]\;\sum_{i=0}^{n-3} \Box\ldots\Box_i^\mid\Box^\mid\Box\ldots\Box_{n-2}\nonumber\\
-\sum_{i=0}^{n-3} \Box\ldots\Box_i^\mid\Box^\mid\Box\ldots\Box_{n-2}
\;[\Box\Box\nonumber\\
\ldots\nonumber\\
+\Box]\;\sum_{i=0}^{n-3} \Box\ldots\Box_i^\mid\Box^\mid\Box\ldots
\Box_{n-2}\;[\Box\ldots\nonumber\\
=\sum_{j=0}^{n-1}\sum_{{i_1+i_2=j \atop i_1\geq 0,i_2\geq 0}}
\sum_{i=0}^{n-j} s(i_1,i_2) \underbrace{\Box\ldots\Box}_{i_1}]
\Box\ldots\Box_i^\mid\Box^\mid\Box\ldots\Box_{n-j}[
\underbrace{\Box\ldots\Box}_{j_2},\label{fini}
\eea
where $s(i_1,i_2)=1$ iff $i_1+i_2=0$ or $i_1i_2\not=0$,
and $s(i_1,i_2)=-1$ otherwise. 

Now, to prove it as promised, we use the fact that an overall 
convergent expression,
when dressed with internal vertex and self-energy corrections, can be
rendered finite by including the appropriate counterterms.
So what we have to do is to do the step in Eq.(\ref{e57})
at the same time for a $n$-loop state and some appropriate chosen
counterterm states so that in each step all subdivergences are
compensated.



We then want to  repeat the step in Eq.(\ref{e57}) $n-1$ times so that
we end up with expressions where the exterior momenta flows only
through one propagator, -1-states, for which we have:  
\beas
\Box\ldots\Box_i^\mid\Box^\mid\Box\ldots\Box_{n}
=B^{i-1}(\Delta)B^{n-i-2}(\Delta)\;{}_{n-1}\!\Omega.
\eeas
These 1-states are easy to calculate. All loops which are free
of the exterior momentum correspond to functions ${}_j\!\Delta$ and the
exterior momentum flow is in ${}_{n-1}\!\Omega$. So we obtain
the usual concatenations with the $B$ operator.
In the spirit of Eq.(\ref{e57}) we would like to write
\bea
\Omega^{(n)}|_1 & = & \Omega^{(n)}|_1 N_{l_2}N_{l_n}+\Omega^{(n)}|_{1,2}
+\Omega^{(n)}|_{1,n}-\Omega^{(n)}|_{1,2,n}.
\label{e60}\eea
Here $|_1$ denotes the loop which is free of the external
momentum already. 
If we now could push the analogy to Eq.(\ref{e57}) further and 
conclude 
\bea
<\Omega^{(n)}|_1 N_{l_2}N_{l_n}> =^? 0,
\label{e61}\eea
we would see the beginning of an algorithm leading to the final state
in Eq.(\ref{fini}). The problem is that now, in Eq.(\ref{e60}), we 
have removed the overall degree of divergence but not all
subdivergences. We have removed all subdivergences not involving
$l_1$, but the loop corresponding to $l_1$ produces a problem we have
not taken care of yet. We can do so by remembering Weinberg's theorem
and picking up the appropriate counterterm expression:
\bea
<\Omega^{(n)}|_1 N_{l_2}N_{l_n}-
Z_1^{(1)}\Omega^{(n-1)}|_1 N_{l_2}N_{l_{n}}>=0,
\label{e62}\eea
where $\Omega^{(n-1)}$ was expressed in $n-1$ loop momenta
$l_2,\ldots,l_{n}$.
We have in our $\Box$ notation:
\beas
<\Box_1^\mid\Box\ldots\Box_n^\mid\;\; - \;\;\Box]\;{}^\mid\Box\ldots\Box_{n-1}^\mid>\\
=<\Box\Box_2^\mid\Box\ldots\Box_n^\mid\;\; -\;\; 
\Box]\;\Box_1^\mid\Box\ldots\Box_{n-1}^\mid>\\
+\;\;<\Box_1^\mid\Box\ldots\Box_{n-1}^\mid\Box_{n}\;\; -\;\; 
\Box]\;{}^\mid\Box\ldots\Box_{n-2}^\mid\Box_{n-1}>\\
-\;\;<\Box\Box_2^\mid\Box\ldots\Box_{n-1}^\mid\Box_{n}\;\; +\;\;
\Box]\; \Box_1^\mid\Box\ldots\Box_{n-2}^\mid\Box_{n-1}>,
\eeas
as a legitimate step.


Applying the above mechanism
$n-1$ times for a $n$-loop graph gives us the sum
over 1-states predicted in the final result above.
At each step, by the very definition of
renormalizability, there is an appropriate set of counterterms available
so that the mechanism is justified.
It is easy to see that the signs in the relation Eq.(\ref{e60})
conspire in the right way to guarantee that each 1-state appears
exactly one time
in the final sum.


We give an example for the case $n=3$:
\beas
{}^\mid\Box\Box\Box^\mid\\
-\;\;\Box]\;{}^\mid\Box\Box^\mid\;\;-\;\;{}^\mid\Box\Box^\mid[\Box\\
-\;\;\Box\Box]{}^\mid\Box^\mid\;\;-\;\;{}^\mid\Box^\mid[\Box\Box
+\;\;\Box]^\mid\Box^\mid[\Box\\
=
{}^\mid\Box\Box^\mid\Box\;\;+\;\;\Box^\mid\Box\Box^\mid\;\;-\;\;\Box^\mid\Box^\mid\Box\\
-\;\;\Box]{}^\mid\Box\Box^\mid\;\;-\;\;{}^\mid\Box\Box^\mid[\Box\\
-\;\;\Box\Box]^\mid\Box^\mid\;\;-\;\;{}^\mid\Box^\mid[\Box\Box
\;\;+\;\;\Box]{}^\mid\Box^\mid[\Box\\
=
{}^\mid\Box^\mid\Box\Box\;\;+\;\;\Box^\mid\Box^\mid\Box\;\;
-\;\;\Box^\mid{}^\mid\Box\Box\\
-\;\;\Box^\mid\Box^\mid[\Box\;\;-\;\;
{}^\mid\Box^\mid\Box[\Box\;\;+\;\;\Box^\mid{}^\mid\Box[\Box\\
+\;\;\Box\Box^\mid\Box^\mid\;\;+\;\;\Box^\mid\Box^\mid\Box\;\;-\;\;
\Box\Box^\mid{}^\mid\Box\\
-\;\;\Box]{}^\mid\Box^\mid\Box\;\;-\;\;
\Box]\Box^\mid\Box^\mid\;\;+\;\;\Box^\mid{}^\mid\Box[\Box\\
-\;\;\Box\Box]{}^\mid\Box^\mid\;\;-\;\;{}^\mid\Box^\mid[\Box\Box
\;\;+\;\;\Box]{}^\mid\Box^\mid[\Box\\
={}^\mid\Box^\mid\Box\Box\;\;+\;\;\Box^\mid\Box^\mid\Box\;\;+\;\;
\Box\Box^\mid\Box^\mid\\
-\;\;\Box]({}^\mid\Box^\mid\Box\;\;+\;\;\Box^\mid\Box^\mid)
-({}^\mid\Box^\mid\Box\;\;+\;\;\Box^\mid\Box^\mid)[\Box\\
-\Box\Box]{}^\mid\Box^\mid\;\;-\;\;{}^\mid\Box^\mid[\Box\Box\\
+\Box]^\mid\Box^\mid[\Box.
\eeas
Here we used $\Box^\mid{}^\mid\Box=0$, as these expressions correspond to
graphs where the exterior momentum does not flow at all through the
graph.
For a massless graph in DR we then have a vanishing tadpole graph,
cf.~Eq.(\ref{dr3}).
\subsection{Links and ladders}
Now let us derive the same result from knot theory.
Let us start again with the two-loop example. According to our rules
in the previous section we obtain:  
{\Large 
\beas
\epsfysize=4mm \epsfbox{knot29.ps}\rightarrow
{}_\Omega\cil\crs\cir_\Delta+{}_\Delta\cil\crs\cir_\Omega\\
=-<\cic>\cic_\Omega-\cic_\Omega<\cic>\\
+\cic\;\cil\wrt\cir_\Omega+\cil\wrt\cir_\Omega \;\cic
\eeas}
This time the correspondence between a Feynman graph and a link diagram 
involves a sum over possibilities opened up by
the overlapping topology.
We have not changed the rules concerning overcrossings, but now
one link will correspond to a vertex correction $\Delta$, the
other link to a self-energy $\Omega$.
The first question we have to answer is: which link serves as the
skeleton and which one as the subdivergence? In case of the 
nested or disjoint subdivergences
considered so far we could always decide this in a unique manner. By
the very definition of an overlapping divergence both links above can play 
the role of the skeleton or subdivergence. 


Accordingly, as indicated above, let us sum over both possibilities.
Applying  now the formalism developed in the previous section, 
taking into account that
the skeleton is necessarily a two point function, the subdivergence 
necessarily a vertex correction, we obtain
\bea
Z_2^{(2)} & = & -<\Delta>\Omega+\Delta\;{}_1\!\Omega
-<\Delta>\Omega+\Delta\;{}_1\!\Omega\nonumber\\
 & = & 2(-<\Delta>\Omega+\Delta\;{}_1\!\Omega),
\label{e63}\eea
which is the desired answer.
Extending this result we investigate the three-loop case:  
{\Large
\beas
{}_\Omega\cil\crs\crs\cir+\cil\crs_{\makebox[0cm]{${}_\Omega$}}\crs\cir
+\cil\crs\crs\cir_\Omega \;\rightarrow\\
-<\cic>[{}_\Omega\cil\crs\cir+2\cil\crs\cir_\Omega]\\
+\cic\;[{}_\Omega\cil\wrt\crs\cir+2\cil\wrt\crs\cir_\Omega\;]\\
\rightarrow \ldots
\eeas}
{\small A three-loop example for the overlapping 
ladder topology.\\[5mm]}  
We indicated the skeleton loop by an $\Omega$.
Explicitly in our one-loop functions we have
\bea
Z_2^{(3)} & = &
2[-<\Delta>\Delta\;{}_1\!\Omega -<\Delta\;{}_1\!\Delta>\Omega\nonumber\\
 & & +
<<\Delta>\Delta>\Omega+
\Delta\;{}_1\!\Delta\;{}_2\!\Omega ]\nonumber\\
& & -2<\Delta>\Delta\;{}_1\!\Omega+<\Delta><\Delta>\Omega
+\Delta\Delta\;{}_2\!\Omega\nonumber\\
 & = & [2\Delta\;{}_1\!\Delta\;{}_2\!\Omega + \Delta\Delta
\;{}_2\!\Omega]\nonumber\\
 & & -[4<\Delta>\Delta\;{}_1\!\Omega]\nonumber\\
 & & -2[<\Delta\;{}_1\!\Delta>\Omega-<<\Delta>\Delta>\Omega]\nonumber\\
 & & +<\Delta><\Delta>\Omega.
\label{e64}\eea
Note that the case where we read the three-loop overlapping
ladder as having a one-loop correction on each side corresponds to two
disjunct subdivergences. That was the reason why we discussed disjunct 
subdivergences in some detail in the previous section.


A comparison with Eq.(\ref{fini}) shows that it is again the correct result.
Indeed, the sum over all possible assignments of the skeleton property to one ring
equals the sum over all the possibilities over which ring has to 
carry the external momentum
flow, and we see that the knot theoretic approach gives the correct answer
immediately. 
\subsection{A generalization}
\begin{figure}[ht]
\epsfysize=7cm \epsfbox{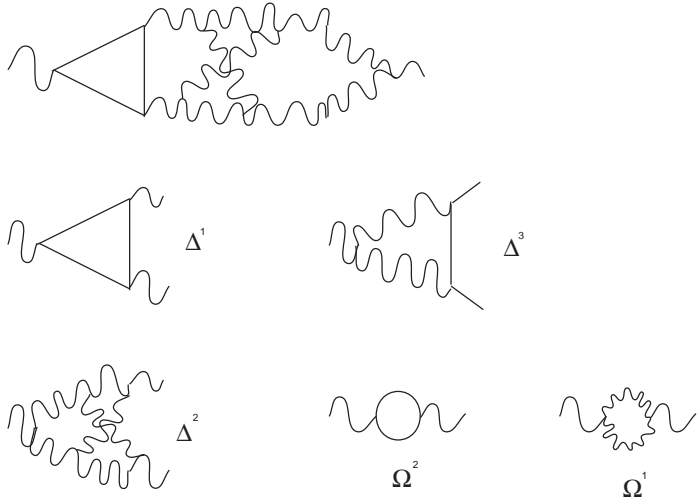}
\caption[Proper notation.]{\label{hab68}\small A more general case which demands a proper notation for 
various types
of diagrams and topologies.}\end{figure}

In analogy to the previous sections we generalize this result to the case of 
various
topologies. To calculate cases like
Fig.(\ref{hab68})  
one has to introduce a proper notation for the various
one-loop functions.
We find as the result in an obvious
generalization of the results of Eq.(\ref{fini}):
\bea
Z_2(\mbox{Fig.(\ref{hab68})}) & = & [\Delta^1\;{}_1\!\Delta^2\;{}_3\!\Omega^1 
+\Delta^2\;{}_2\!\Delta^3\;{}_3\!\Omega^2\nonumber\\
& & + \Delta^1\Delta^2
\;{}_3\!\Omega^1]\nonumber\\
& & -[<\Delta^1>\Delta^2\;{}_2\!\Omega^1
-<\Delta^2>\Delta^3\;{}_1\!\Omega^2\nonumber\\
& & 
-<\Delta^1>\Delta^2\;{}_2\!\Omega^1
-<\Delta^2>\Delta^1\;{}_1\!\Omega^1
]\nonumber\\
& & -[<\Delta^1\;{}_1\!\Delta^2>\Omega^1-<<\Delta^1>\Delta^2>\Omega^1\nonumber\\
& & 
  +<\Delta^2\;{}_2\!\Delta^3>\Omega^2-<<\Delta^2>\Delta^3>\Omega^2
  ]\nonumber\\
& & +<\Delta^1><\Delta^2>\Omega^1.
\label{e65}\eea
These are the 12 expected terms ($2^{n-1}=4$, from the skein
relation, times 3 possibilities to assign the skeleton property).

\begin{figure}[ht] 
\epsfysize=5cm \epsfbox{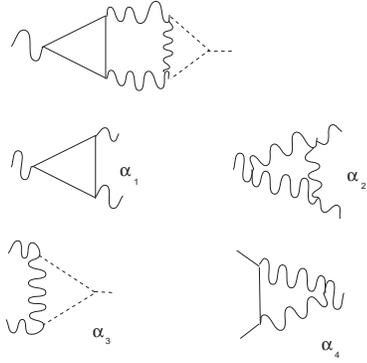}
\caption[Unsymmetric cases.]{\label{hab69}\small From the left and from the right we get different 
$\alpha$'s.}\end{figure}

In case that the incoming particle and the
outgoing particle are different this would blow up the possibilities,
e.g. Fig.(\ref{hab69}).  
This completes our treatment of overlapping topologies. 
We succeeded in representing a sufficient large class of Feynman graphs
in terms of concatenated one-loop functions. In general,
we arrived at the result that the algebraic
structure underlying overlapping divergences is similar
to the previous cases, as long as we sum over all
possibilities how to identify the skeleton graph. 
Referring to the research problem quoted in
section three, we learn that once we have solved how to
concatenate higher terms in the skeleton expansion,
this knowledge would serve for nested as well as overlapping divergences.

\clearpage
\section{No Knots, No Transcendentals}

What has been achieved so far? In previous sections we expressed familiar
results of renormalization
theory in terms of a relative simple one-loop algebra.
Occasionally we compared this with the algebraic structure dictated
by link diagrams, associated to the Feynman graphs in some
prescribed manner. In fact, the results of the previous sections give a
complete account on how
to treat general ladder-type topologies. 
But at some stage we want to assign knots to  Feynman
graphs.
More specifically, we want to establish that the transcendentals
obtained by calculating the UV-divergences of a topological 
sufficiently complicated
Feynman graph characterize the knot which one obtains when
considering the link diagram, associated to the graph
according to the rules of section three.
The following observation which is the subject of this section justifies
the identification of knots with Feynman diagrams: 
topologically simple rainbow or ladder graphs are free
of transcendentals. This sets the stage for identifying knots with
Feynman graphs. We will report on the successes of such an identification
in the next section, where we will investigate the
triangle {\em Knot Theory  - Field Theory  - Number Theory}
more closely.


We now have to investigate
all Feynman graphs which correspond to simple topologies.
We define a topology as simple when it corresponds to a 
braid group expression
of at most second degree in all braid generators. These braids,
after applying the skein relation $n$-times, will generate an unknot
with writhe number $n$, amongst even simpler (splitted) terms. 
This property of
being free of knots should reflect itself in the $Z$-factors as we want to 
have a chance to identify $Z$-factors with knots.
We claim that there is a one-to-one correspondence between the transcendental
numbers in the divergent part of a diagram and the knot in its link
diagram, so we better prove that such transcendentals do not
appear in simple topologies.


To give the reader an idea what sort of cancellations are necessary we give
in the appendix a result for a seven loop ladder topology.
The graph itself has transcendentals in abundance in its 
divergent part, but when we add its counterterm expressions they all
disappear.


We proceed in the following way. We first investigate how this
cancellation of transcendentals appears for some special choice
(ladder topology in a massless Yukawa theory, see appendix for
Feynman rules)
of the $\Delta$ function. In fact, we show that for an arbitrary
loop order $n$ the highest possible transcendental $\zeta(n-1)$
will not arise. The proof for the other transcendentals is similar. 


We indicate how the proof has to be generalized for the arbitrary 
tensor case.
Then we discuss this cancellation of transcendentals from a broader 
viewpoint
and give a general inductive proof.
Recently, Bob Delbourgo et.al.\cite{bobnew} have given an independent
proof for the statement. They transferred the Schwinger Dyson Equation
for bare ladder diagrams to a differential equation, which they
solved for general  $D$, where $D$ refers to the generalized
dimension of dimensional regularization. Upon expanding
their non-perturbative solution in the coupling constant one
confirms the results of this section.


The result here is the following statement: Any Feynman
graph in a renormalizable theory which corresponds to a simple topology
in its link diagram  gives rise to only rational divergences
when calculated together with its counterterm graphs. 
Here rational divergences means
that only rational numbers appear as coefficients of the proper 
Laurent part in
the DR expansion parameter $\epsilon$.
\subsection{A combinatorical proof}
Let us start considering the following function, which exhibits all typical
properties of the observed cancellations:
\beas
{}_j\!\Delta: &   & 
\int d^Dk \frac{(k^2)^{-\varepsilon j}}{k^2 (k+q)^2}
=:(q^2)^{-\varepsilon(j+1)}{}_j\!\Delta,\\
\Rightarrow {}_j\!\Delta & = & 
\frac{\Gamma(1+(j+1)\varepsilon)\Gamma(1-\varepsilon)
\Gamma(1-(j+1)\varepsilon)}{(j+1)\varepsilon(1-(j+2)\varepsilon)
\Gamma(1+j\varepsilon)\Gamma(1-(j+2)\varepsilon)}.
\eeas
Define
\beas
P_n & := & \prod_{i=0}^{n-1} {}_i\!\Delta,\\
\Rightarrow P_n & = & \frac{(\Gamma(1-\varepsilon))^{n+1}
\Gamma(1+n\varepsilon)}{n!\varepsilon^n(1-2\varepsilon)\ldots
(1-(n+1)\varepsilon)\Gamma(1-(n+1)\varepsilon)}.
\eeas
Now use
\bea
\Gamma(1-z)=\exp(\gamma z)\exp(\sum_{j=2}^\infty \frac{\zeta(j)}{j}
z^j).\label{exp}
\eea
It follows
\beas
P_n & = & \frac{1}{n!\epsilon^n (1-2\epsilon)\ldots
(1-(n+1)\epsilon)}\exp(-n\gamma\epsilon)\times\\
 & & \exp(\sum_{j=2}^\infty
\frac{\zeta(j)}{j} \epsilon^j[n+1+(-n)^j-(n+1)^j]).
\eeas
We conclude immediately that $\zeta(2)$ can not appear in a $Z$-factor
contribution as its coefficient is $(n+1+n^2-(n+1)^2)=-n$, so it
can be absorbed in a redefined coupling constant in the same way as $\gamma$:
\beas
g\mu^{-\epsilon}\rightarrow
g\tilde{\mu}^{-\epsilon},\;\tilde{\mu}=\mu\exp(
(\gamma+\epsilon\zeta(2)/2)/2).
\eeas
As any contribution at the $n$-loop level factors
$\alpha^n$, where $\alpha:=g^2/4\pi$, we see that this modification
deletes any $\zeta(2)$ dependence,
an argument which is familiar from $\overline{MS}$ schemes.\footnote{
We remind the reader that the $\mu$-dependence comes in
as we demand a dimensionless coupling constant in dimensional
regularization. In $\overline{MS}$ schemes one uses
this $\mu$-dependence to absorb Eulers constant $\gamma$.
The demand that physical quantities are independent of $\mu$
results in the so-called renormalization group equations.}


More subtle is the cancellation of higher transcendental $\zeta$'s.
As an example let us consider $\zeta(n-1)$, the highest possible transcendental
appearing in a $n$-loop calculation.
Only the highest pole $\frac{1}{\epsilon^n}$ can generate it, so we have to
consider
\beas
\frac{1}{n!}\frac{1}{\epsilon^n}\frac{\zeta(n-1)}{n-1}\epsilon^{n-1}
[n+1+(-n)^{n-1}-(n+1)^{n-1}],
\eeas
which is the contribution of $P_n$ to $\zeta(n-1)$.


Also the counterterm expressions have to provide their highest pole.
Let us study a 4-loop example:
\bea
Z_1^{(4)} & = & <P_4-<P_1>P_3-<P_2-<P_1>P_1>P_2\nonumber\\
  & & -<P_3-<P_1>P_2-<P_2-<P_1>P_1>P_1>P_1>.\label{ex4}
\eea
It follows for the coefficient of  $\zeta(3)$:
\beas
S_4(3) & := & \frac{1}{\epsilon}\frac{1}{4!0!}\frac{1}{3}[5+(-4)^3-5^3]\\
 & & -\frac{1}{\epsilon}\frac{1}{3!1!}\frac{1}{3}[4+(-3)^3-4^3]\\
 & & -\frac{1}{\epsilon}\frac{1}{2!2!}\frac{1}{3}[3+(-2)^3-3^3]\\
 & & -\frac{1}{\epsilon}\frac{1}{1!1!2!}\frac{1}{3}[3+(-2)^3-3^3]\\
 & & -\frac{1}{\epsilon}\frac{1}{3!1!}\frac{1}{3}[2+(-1)^3-2^3],\\
 & & \ldots
\eeas
where $\ldots$  refers to the last three terms in Eq.(\ref{ex4}) which add to 
zero.
We have  for the coefficient   of $\zeta(3)$
\beas  
S_4(3)=\frac{1}{\epsilon}\sum_{i=1}^4
(-1)^i \frac{1}{i!(4-i)!}\frac{1}{3}[i+1+(-i)^3-(i+1)^3]=0.
\eeas
It is easy to see that  in general the coefficient of
$\zeta(n-1)$ is given by
\beas
S_n(n-1):=\frac{1}{\epsilon}\frac{1}{n-1}\sum_{i=1}^n
(-1)^i \frac{1}{i!(n-i)!}[i+1+(-i)^{n-1}-(i+1)^{n-1}].
\eeas
So we have to show $S_n(n-1)=0$.
Let us use
\beas
(1-a)^n=\sum_{i=0}^{n}\frac{(-1)^i a^i n!}{i! (n-i)!}
= \sum_{i=0}^\infty \frac{(-1)^i a^i \Gamma(n+1)}{\Gamma(i+1)\Gamma(n-
i+1)},
\eeas
from which we conclude
\beas
\delta_{n,0} = \sum_{i=0}^\infty \frac{(-1)^i}{\Gamma(i+1)\Gamma(n-
i+1)}.
\eeas


Consider 
\beas
T_n(r):= \sum_{i=0}^{\infty}\frac{(-1)^i i^r}{\Gamma(i+1)\Gamma(n-i+1)},\;\;
 r>0.
\eeas
which appears in $S_n(n-1)$.
It is possible to express $i^r$ as a linear combination of terms
\beas
i(i-1)\ldots(i-r+1)+i(i-1)\ldots(i-r+2)+\ldots+i,
\eeas
so that
\beas
 i^r & = & \sum_{j=1}^{r} c_{rj} \frac{\Gamma(i+1)}{\Gamma(i-j+1)},
\eeas
and we obtain
\beas
T_n(r) & = & \sum_{i,j=1}^r (-1)^i \frac{c_{rj}}{
 \Gamma(i-j+1)\Gamma(n-i)}\\
 & = & \sum_{i,j}^r (-1)^{i+j} \frac{c_{rj}}{
 \Gamma(i+1)\Gamma(n-i-j)}\\
 & = & \sum_{j=1}^r (-1)^j c_{rj} \delta_{n-j,0}\\
 & = & (-1)^n c_{rn}.
\eeas
We have
\beas
c_{rn}=0,\; r<n,\\
c_{rn}=1,\;r=n,
\eeas and finally
\beas
T_n(r)=0\;\;\mbox{for $0<r<n$},\\
T_n(n)=(-1)^n.
\eeas
By a similar argument we can show that 
\beas
U_n(r) & := & \sum_{i=0}^{\infty}\frac{(-1)^i (i+1)^r}{i! (n-i)!}
=0 \;\;\mbox{for $0<r<n$},\\
U_n(n) & = & (-1)^n.
\eeas
Combining everything we find
\beas
S_n(r)=0,\; r\leq n,
\eeas
which includes the desired result for $S_n(n-1)$.
This proof was first obtained in collaboration with 
Bob Delbourgo.


Some comments might be appropriate. We investigated the above case for
$D=4-2\epsilon$ dimensions, where the function ${}_j\!\Delta$
corresponds for example to a Green function in 
 Yukawa theory. Nevertheless it is the generic example
for all possible cases in a renormalizable theory. In the case of
tensor integrals or dimensions other than four our function
would only be modified by a polynomial of $\epsilon$
which multiplies it. This leaves the basic structure of the sums above
unaffected and the reasoning remains unchanged.
We could now proceed to show the absence of the other
transcendentals $S_n(r)$ along similar lines.
In fact, by inspection one sees that the cases for $ \zeta(n-i),
i>1$ follow the same pattern. 


Rather we prefer to give a more general argument establishing
rational contributions for $Z$ factors from simple topologies for
all renormalizable theories.
\subsection{A global argument}
We argue for the case of a vertex function, the propagator follows
similarly.


Let a generalized one-loop function ${}_j\Delta$ be given.
Assume we can render it finite with the help of another
function ${}_j\bar{\Delta}$
\begin{equation}
<{}_j{\Delta}-{}_j\bar{\Delta}>=0.
\end{equation}
According to the Schwinger Dyson equation, we would get the
$Z$ factor as
\begin{equation}
Z=<G_R \circ {}_j \! \Delta>=<G_R \circ {}_j\!\bar{\Delta}>,\label{conc}
\end{equation}
where we solved the Schwinger Dyson equation Eq.(\ref{sd}) for $Z$ 
and discarded the finite
term containing the full renormalized vertex. It is crucial
that in the Schwinger Dyson equation only finite renormalized
Green functions appear. 


The following examples explain the meaning of our concatenation
$\circ$ in Eq.(\ref{conc}):
\begin{eqnarray*}
\mbox{2 loops:} & & G_R^{(1)}={}_0\!\Delta-<{}_0\!\Delta>\;\; \Rightarrow
Z^{(2)}=<{}_0\!\Delta {}_1\!\Delta -<{}_0\!\Delta>{}_0\!\Delta>,\\
\mbox{3 loops:} & & G_R^{(2)}= 
{}_0\!\Delta {}_1\!\Delta
-<{}_0\!\Delta> {}_0\!\Delta -<{}_0\!\Delta {}_1\!\Delta
-<{}_0\!\Delta> {}_0\!\Delta>\;\;\Rightarrow\\
 & & Z^{(3)}=
<{}_0\!\Delta {}_1\!\Delta {}_2\!\Delta
-<{}_0\!\Delta> {}_0\!\Delta {}_1\!\Delta -<{}_0\!\Delta {}_1\!\Delta
-<{}_0\!\Delta> {}_0\!\Delta>{}_0\!\Delta>.
\end{eqnarray*}



Now assume that ${}_j\bar{\Delta}$ enjoys the following
property:
\begin{equation}
<(\prod_{i=0}^{j-1}{}_i\!\Delta){}_j\bar{\Delta}> \in {\bf Q}.
\label{rat}
\end{equation}
It removes all transcendentals from lower level products
of $\Delta$ functions.


Then we would easily achieve our aim by induction: 
$Z^{(n+1)}$ is rational if and only if $Z^{(n)}$ is rational.
By the very definition of $Z^{(n+1)}$ it is given as a
sum of products of lower level $Z$ factors with
concatenated $\Delta$-functions, where we replace the last one
by the barred form, as in the above  examples. 
Eq.(\ref{conc}) would give the $Z$ factor immediately
as a product of rational functions, if Eq.(\ref{rat}) holds.


So the rationality of the $Z$ factor is proven if we can
establish the existence of a family of functions
${}_j\bar{\Delta}$ such that
\begin{itemize}
\item
${}_j\Delta - {}_j\bar{\Delta}$ is finite by power-counting
\item
Eq.(\ref{rat}) holds.
\end{itemize}
The structure of the generalized one-loop functions allows us to write
\begin{equation}
{}_m\!\Delta=:f_m(D-4)\;\exp{g_m(D-4)}
\Rightarrow \prod_{i=0}^{j-1}{}_i\!\Delta=:p_j(D-4)\;
\exp{[q_j(D-4)]},
\end{equation}
which defines $f_m,g_m,p_j,q_j$. Here $f_j, p_j$ are rational functions
of $(D-4)$.
We set
\begin{equation}
\bar{{}_j\!\Delta}:=f_j(D-4)\exp{(-q_{j-1}(D-4))}.
\end{equation}
We immediately conclude
that Eq.(\ref{rat})
holds:
\begin{equation}
<(\prod_{i=0}^{j-1}{}_i\!\Delta){}_j\bar{\Delta}> 
=p_j(D-4)f_j(D-4).
\end{equation}
The demand that $\bar{{}_j\!\Delta}$ removes the
overall divergence by powercounting 
is easy to achieve. We give an explicit construction in 
appendix E. We only have to guarantee that the function
has the correct asymptotic behaviour.
It follows that 
\begin{equation}
<G_R\circ({}_j\!\Delta -\bar{{}_j\!\Delta})>=0,
\end{equation}
which completes our proof. Again we explored the fact that 
a function which is overall finite as
$\sim ({}_j\!\Delta -\bar{{}_j\!\Delta})$
is finite when all its subdivergences are subtracted. 


For the propagator we can argue in a similar way. 
We only have to take care of the modification in the 
Schwinger Dyson equations for this case.


We further observe that the difference 
\begin{equation}
\left[\prod_{i=0}^{j-1}{}_i\!\Delta\right]{}_j\bar{\Delta}-\prod_{i=0}^{
j}{}_i\!\Delta\;,
\end{equation}
is of order
$\epsilon$. So we can allow for one further
dressing to be present at some other line or vertex.
The difference would be finite for this case.
This is still sufficient to replace the ultimate
${}_j\!\Delta$ by ${}_j\bar{\Delta}$, and thus the proof still
goes through.
On the other hand any further dressing spoils our argument. 
This is precisely the behaviour observed in
previous sections. Ladder and rainbow topologies
considered there were free of transcendentals.
We can compare it also with the
results in \cite{newplb}, which were already briefly mentioned
in section five. There, the examples in Table(\ref{t1})
confirm these findings. 
We expect transcendentals of the type $\zeta(n), n\geq 3$
to appear as soon as we have three or more dressings. This 
is indeed the case, and their
appearance matches all ways we can assign link diagrams
to the Feynman graph. But the real confirmation
of a connection between knots and numbers
via Feynman graphs will be given in the next section, where
we report on recent results exploring this idea by calculating
diagrams mostly free of subdivergences, but of complicated topology.


In this section we have learned that the the very special sort of
topological simple Feynman diagrams will only generate rational contributions
to $MS$ $Z$-factors.

\clearpage
\section{Knots and Transcendentals}
\newcommand{\fk}[4]{\st{#1}\so{#2}\st{#3}\so{#4}}
\newcommand{\sk}[6]{\so{#1}\st{#2}\so{#3}\se{#4}\st{#5}\se{#6}}
\newcommand{\se}[1]{\sigma_3^{#1}}
\newcommand{\st}[1]{\sigma_2^{#1}}
\newcommand{\so}[1]{\sigma_1^{#1}}
\newcommand{\tm}{\cdot}


In this section we describe results which were obtained recently
\cite{plb,pisa,bdk,bgk,db,newplb}.\footnote{Note added in
proof: While this paper was written, further results
along the same lines were obtained in
\cite{BK15,newplb,4TR}.} In these six publications
a fascinating connection
between field theory, number theory, and knot theory emerges.
The starting point of this connection is the result 
of the previous section connecting topologically simple
graphs with the absence of knots in their link diagram,
and with the corresponding absence of transcendentals in their 
counterterms.


So the previous section suggests that the transcendental coefficients
of the divergences are related to the topology of the 
diagram. This relation should be via knot theory. 
As we will see, field theory initiated the invention of a 
knot-to-number dictionary, which in turn spurred new findings in
number theory and opened a new route for calculations in
field theory.



In the following, we will first give elementary examples how
to obtain $(2,q)$ torus knots in topologically non-simple diagrams, and 
compare with the transcendentals $\zeta(q)$ in their counterterms.


Then we will comment in detail on all the six publications
mentioned above.
Each of them gives new insights and support to the connection
between renormalization, knot theory and number theory:
\begin{itemize}
\item[(\cite{plb})]
Here it was shown that at the six loop level for the first time
a $(3,4)$ torus knot was obtained. This matches the appearance
of the transcendental $\zeta(3,5)$ in the counterterm.
\item[(\cite{pisa})]
The identification of knots with numbers is used to calculate the scheme
independent part of the $\beta$-function in $\phi^4$
theory to the seven loop level.
Knots obtained from the diagrams match and thus predict the transcendentals
apparent in the counterterms.
\item[(\cite{bdk})]
The existence of a skein relation together with the Ward identities
explains the cancellations of transcendentals long known
for the quenched $\beta$-function in QED.
\item[(\cite{bgk})]
Bubble insertions in a basic two-loop topology deliver 
Euler double sums. The restricted class of knots 
obtained from these diagrams governs the number of independent
double sums, and informs number theory.
\item[(\cite{db})]
All this led Broadhurst to conjecture a formula for the number
of independent Euler sums, and thus solving a major number-theoretic
problem. The findings are supported by  striking numerical evidence, 
and agree with
field/knot-theoretic expectations whenever a comparison is possible.
\item[(\cite{newplb})]
Generalized one-loop functions behave under renormalization
in accordance with the demands of link theory. This confirms
the identification of knots and numbers in \cite{bgk}
and predicts the entanglement for more general cases.
\end{itemize}


\subsection{The $(2,q)$ torus knots and $\zeta(q)$}
The most prominent transcendentals stem
from the Riemann $\zeta$ function evaluated at
odd integer argument.  They arise
from the expansion of the $\Gamma$ function near unit argument
(cf.~Eq.(\ref{exp})) \cite{smirnov}.

\begin{figure}[ht]
\epsfysize=6cm \epsfbox{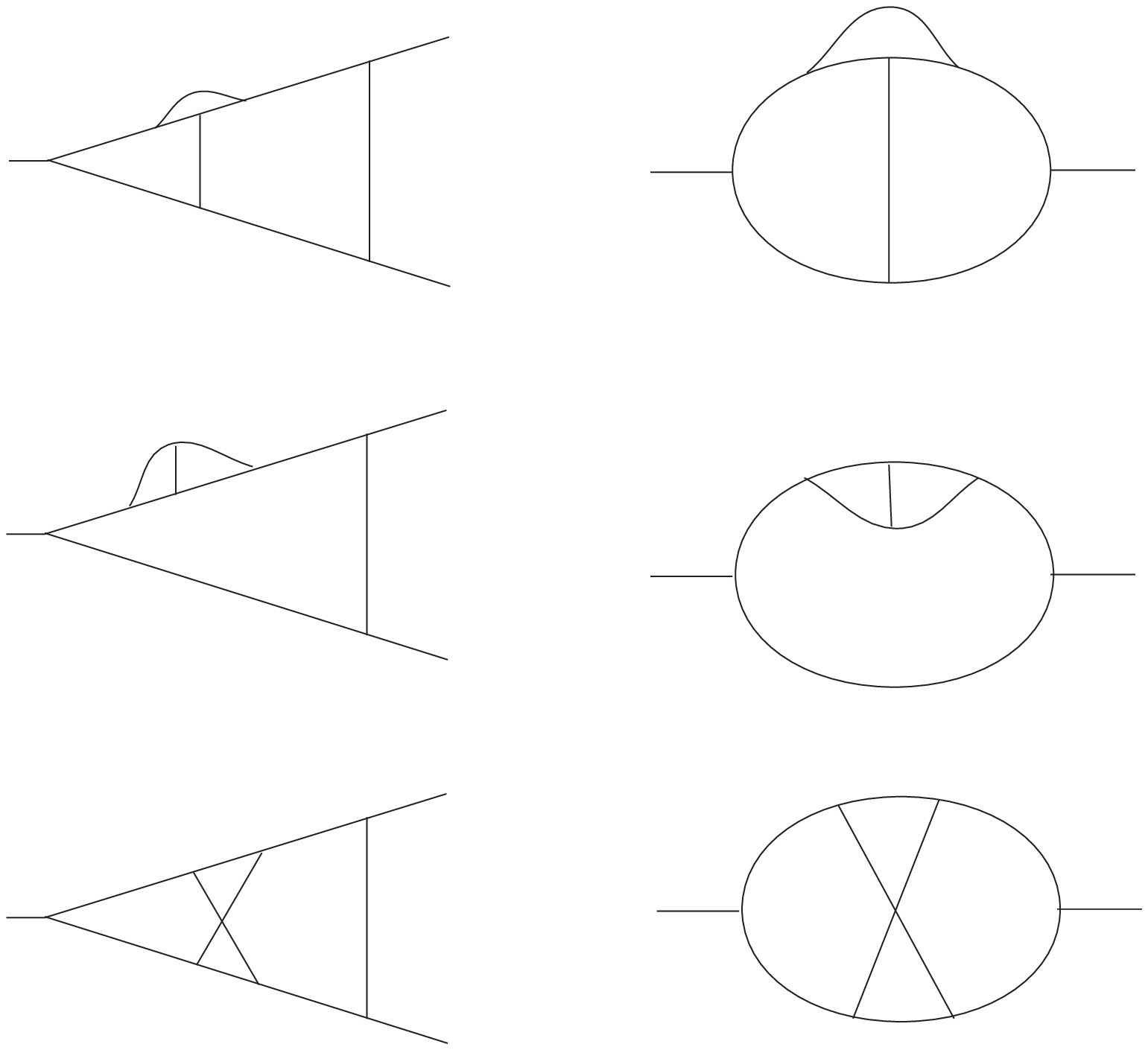}
\epsfbox{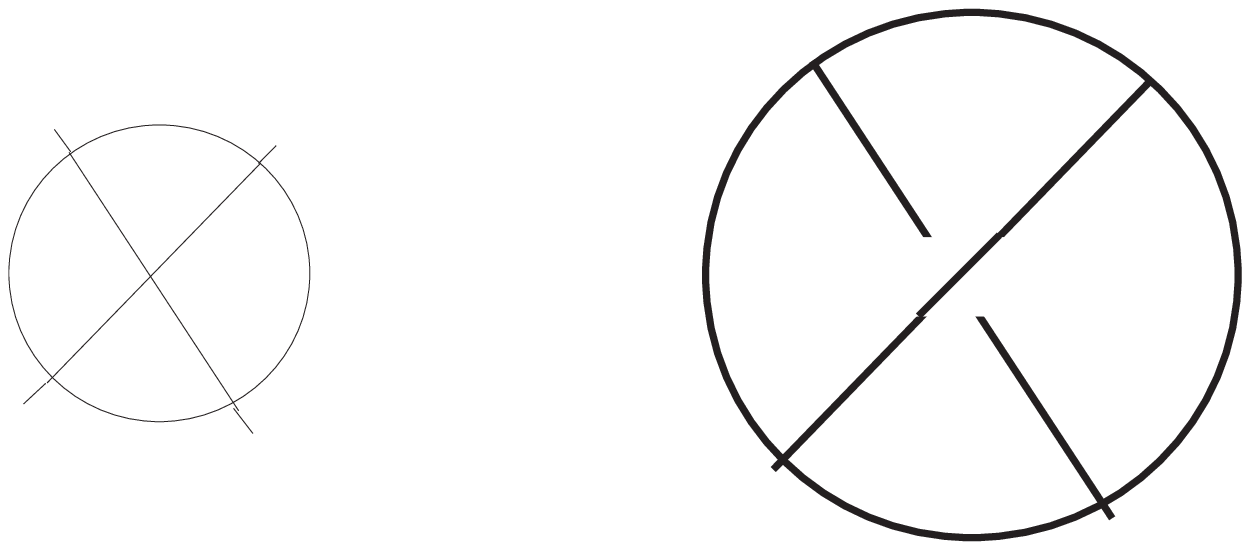}
\caption[$\zeta(3)$ in three-loop graphs.]{\label{hab81}\small 
These three-loop graphs involve $\zeta(3)$ in their
divergent part. They are all derived from the general
topology on the rhs, by coupling two, three or four external particles
at different places. The last graph on the lhs below is
a $\phi^4$ graph and free of subdivergences.}\end{figure}
\begin{figure}[ht]
\epsfysize=3cm \epsfbox{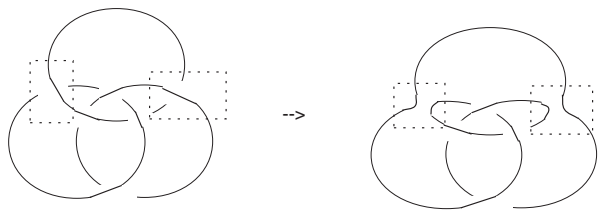}
\caption[Link diagrams for $\zeta(3)$.]{\label{hab82}\small 
The corresponding link diagrams. 
The dashed rectangles indicates
where the $Y$ part of the skein relation has been applied twice.}
\end{figure}

The first non-trivial candidate would be $\zeta(3)$.
Let us consider some three-loop graphs
as shown in Fig.(\ref{hab81}).
It is well known that all the above graphs, even after adding their
counterterms, give us a non-vanishing coefficient for $\zeta(3)$. 
For example the overlapping massless two-point
function, a prominent example for generating
a transcendental series in $\varepsilon$
\cite{exp}, appears as a subgraph in all of the above graphs.
Its finite value is $6\zeta(3)$ (for a scalar theory, with the coefficient
6 changing to other rational values for other theories).
The counterterm subtracts only the divergent part of
the subgraph, and thus there remains
the term proportional to $\zeta(3)$ which is then
multiplied by the divergence from the
final loop integration.
The coefficient of $\zeta(3)$ naturally depends on the actual theory
under consideration;
we restrict ourselves to scalar theories at the moment for convenience.
In appendix D the reader will find an explicit calculation of the
counterterm for the $\phi^4$ diagram in Fig.(\ref{hab81}),
deriving the result $6\zeta(3)/\epsilon$ for it.
In all cases the corresponding link diagrams have a topology 
which for the first time
involves a mutual entanglement of three links, see Fig.(\ref{hab82}).
We see that the trefoil knot appears when we apply the $Y$ part of
the skein relation twice.
This topology gives us the braid expression
\bea
\sigma_1\sigma_2\sigma_1^2\sigma_2\sigma_1.
\label{e68}\eea
For the first time we encounter a more complicated word in braid group
generators. 
Applying a skein relation, it splits into various terms.
Terms which separate into disjoint components should
correspond to counterterm graphs, according to our
previous experience. If there are no subdivergences, all
the subgraphs are finite, and thus nullified by the projection
onto their divergence. To see the pure knot content,
consider for example the $\phi^4$ graph in Fig.(\ref{hab81}). 
It is free of subdivergences. So the only remaining term
in the skein relation comes from considering the application of the $Y$
part. Our aim is to apply it twice in a way which delivers a
one-component knot. No matter how we do this, the
result is always of the form $\sigma_1\sigma_2\sigma_1^2$.
This is the trefoil knot, with an extra curl in it.

\begin{figure}[ht]
\epsfysize=2cm \epsfbox{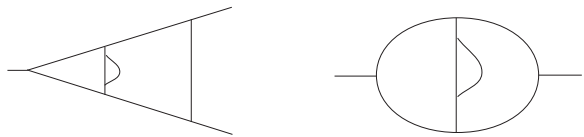}
\caption[No $\zeta(3)$ in these examples.]{\label{hab83}\small 
Two examples which do not produce $\zeta(3)$,
nor the trefoil knot.}\end{figure}

Let us comment on the following two graphs of Fig.(\ref{hab83}).

\begin{figure}[ht]
\epsfysize=4cm \epsfbox{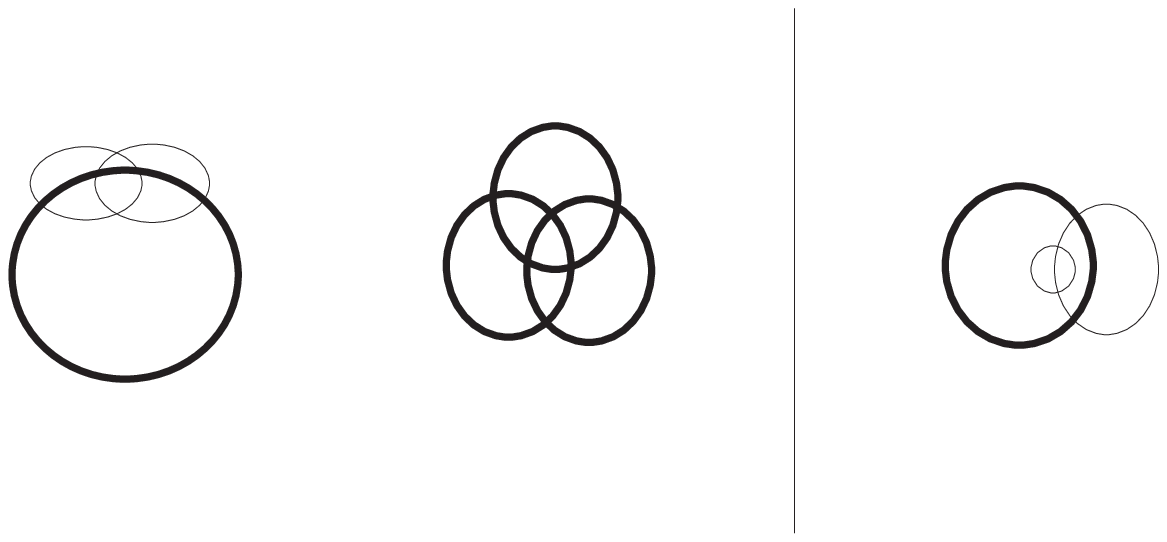}
\caption[How to get the trefoil.]{\label{hab82a}\small 
We get the trefoil via these link diagrams. The second one on the lhs
is the one for the $\phi^4$ diagram (the whole diagram
is a skeleton graph), while all the other ones in
Fig.(\ref{hab81}) provide a one-loop skeleton (thick line), dressed
with two-loop subdivergences. On the rhs
we give the link diagram for the graphs of Fig.(\ref{hab83}).}
\end{figure}

In Fig.(\ref{hab82a}) we compare the generation of the link diagram
for them with the generation of link diagrams
for the examples in Fig.(\ref{hab81}).


Note that in Fig.(\ref{hab81}) the example
for $\phi^4$ theory is unique in the respect that it is a
pure skeleton graph itself, without any subdivergences.
Summarizing,
{\em we are now prepared to assign $\zeta(3)$ to the trefoil topology}.



\begin{figure}[ht]
\epsfysize=1in \epsfbox{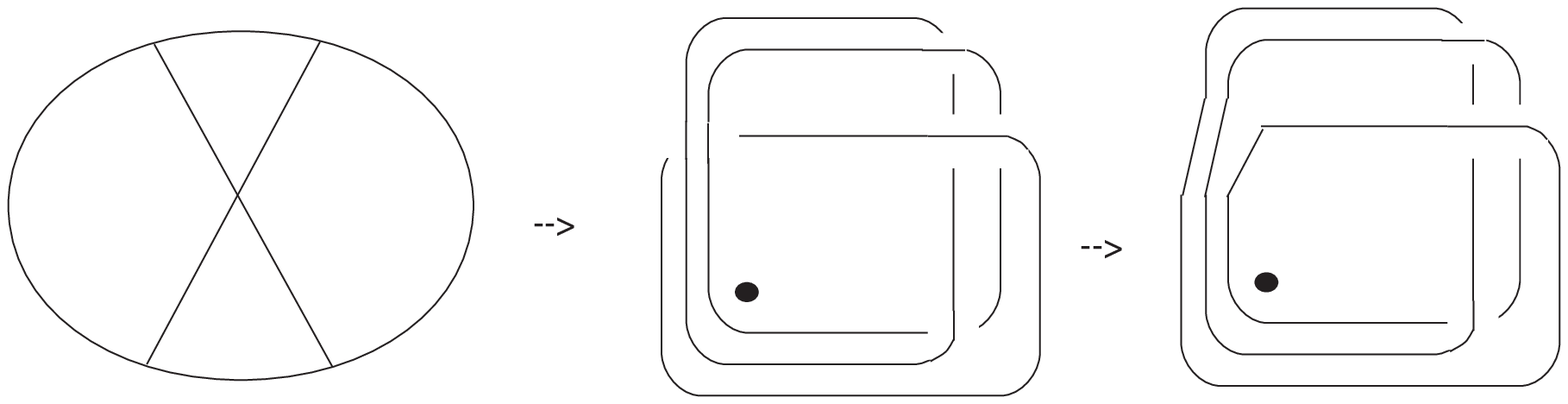}
\epsfysize=2cm \epsfbox{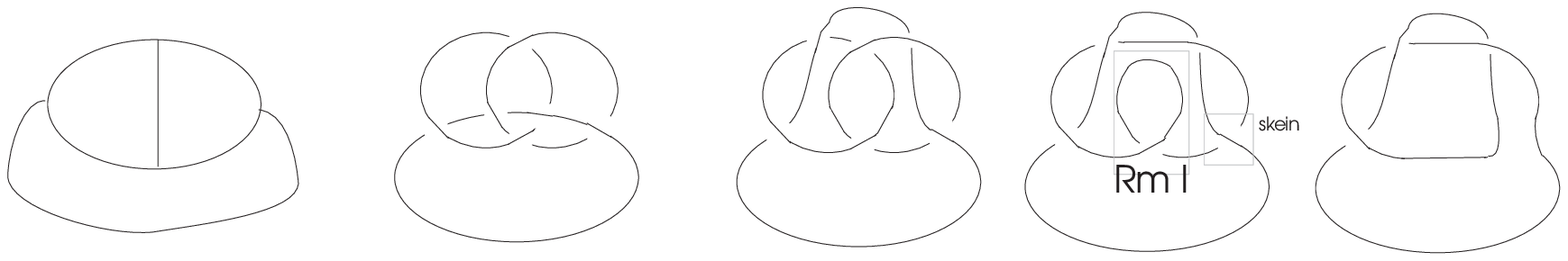}
\caption[A shortcut for the trefoil.]{\label{hab85}\small  
An economic way from a Feynman diagram to a link
diagram to a knot. We omit to indicate exterior momenta. So the Feynman
graph on the left is generic to all the graphs of
Fig.(\ref{hab81}). Note that the topological non-simple structure of the graph,
the crossing of its internal propagators,
directly reflects itself in the braid group expression
Eq.(\ref{e68}), by having the generator
$\sigma_1$ in fourth power. We give two equivalent link diagrams,
resulting from different momentum flows.}\end{figure}

\begin{figure}[ht]
\epsfxsize=5in \epsfbox{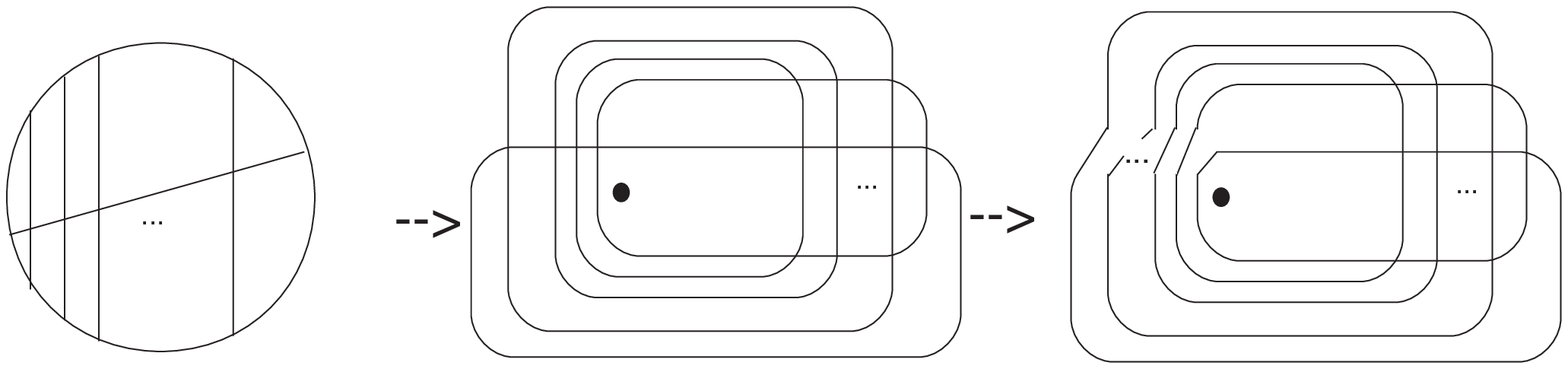}
\epsfbox{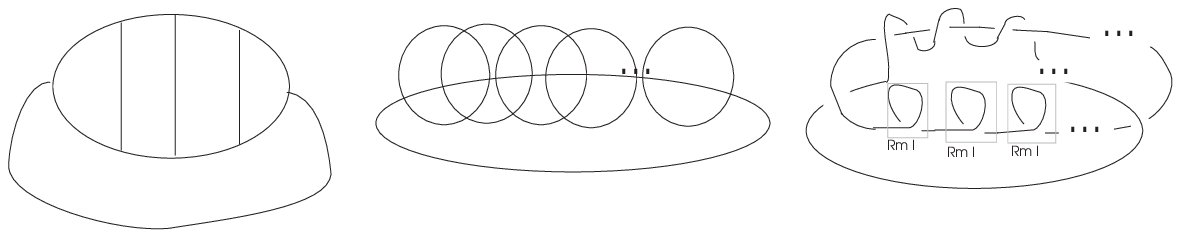}
\caption[The $n$-loop shortcut.]{\label{hab86}\small 
The same for the $n$-loop case. The topology of the
Feynman diagram is a slashed ladder. External particles may couple
at arbitrary places. We observe
that  $n-1$ loops
realize the ladder topology, and one loop goes through the middle of
this $n-1$ cable. We omitted over/undercrossings
in the above picture, as they are clear from the orientation of the
loop momenta. We give two equivalent link diagrams.
They produce the same $(2,2n-3)$ torus knot at the end.}\end{figure}

We note that we can assign loop momenta always in a way that they
all encircle a given point $\bullet$ inside the diagram
counterclockwise. This point 
corresponds to an axis in a closed braid diagram where all strands
are oriented to encircle it in this manner. 
Having defined such an axis in the
Feynman graph we replace the momentum flow by strands according to our rules
of section three
and read off the braid group expression from the link diagram.
Now, as we have more complicated topologies, over/undercrossings in our
link
diagram must either correspond to a vertex in the Feynman diagram or to
a crossing of propagators in the Feynman diagram, cf.~Fig.(\ref{hab81}).
For the trefoil example we find Fig.(\ref{hab85}).
The trefoil knot is the $(2,3)$ torus knot \cite{knots}.
We can readily generalize this example to the $n$-loop case, demonstrated
in Fig.(\ref{hab86}).


It is a well known fact that the Feynman graphs with a
slashed ladder topology are proportional
to $\zeta(2n-3)$ at the $n$ loop level \cite{ussu}. The corresponding
link diagram generates the knot $\sigma^{2n-3}$, which is the
$(2,2n-3)$ torus knot. To see this we read off from Fig.(\ref{hab86})
the braid group expression:
\beas
\sigma_{n-1}\ldots\sigma_1\sigma_2\ldots\sigma_{n-1}\sigma_1
\ldots\sigma_{n-2}=
\sigma_{n-2}^{2n-3},
\eeas
after applying Markov- and Reidemeister-moves.
For example, chosing $n=4$, we calculate
\beas
\sigma_3\sigma_2\sigma_1\sigma_2\sigma_3\sigma_1\sigma_2=
\sigma_2\sigma_1\sigma_2\sigma_1\sigma_3\sigma_2\sigma_3=\\ 
\sigma_2\sigma_1\sigma_2\sigma_1\sigma_2\sigma_3\sigma_2= 
\sigma_2\sigma_2\sigma_1\sigma_2\sigma_1\sigma_2\sigma_3=\\ 
\sigma_2\sigma_2\sigma_1\sigma_2\sigma_1\sigma_2=
\sigma_2\sigma_2\sigma_2\sigma_1\sigma_2\sigma_2=\\
\sigma_2^5\sigma_1=\sigma_2^5,
\eeas
and one easily proves the result for arbitrary $n$ by utilizing
$\sigma_i\sigma_{i+1}\sigma_{i}=\sigma_{i+1}\sigma_i\sigma_{i+1}$. 
For the second representation given in Fig.(\ref{hab86})
the appearance of the $(2,2n-3)$ torus knot is verified immediately.

\begin{figure}[ht]\epsfysize=8cm \epsfbox{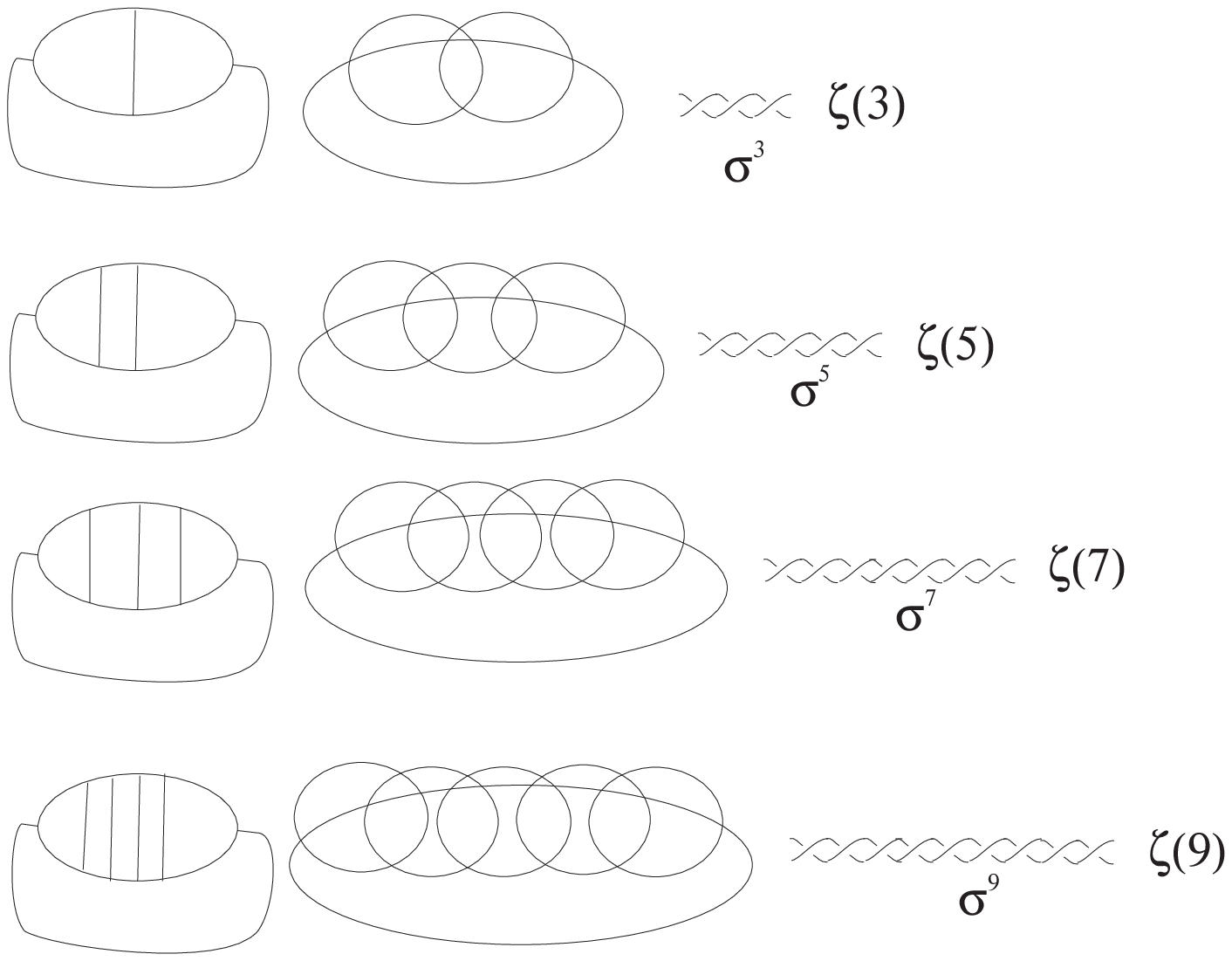}
\caption[The correspondence between $(2,q)$ torus knots
and $\zeta(q)=\sum_{n=1}^\infty\frac{1}{n^q}\sim \sigma^q$.]{\label{hab87}
\small 
A beautiful coincidence between $(2,q)$ torus knots and transcendental
counterterms $\zeta(q)=\sum_{n=1}^\infty\frac{1}{n^q}\sim \sigma^q$. 
We note that the powers of braid generators match
the powers in the infinite sum representation of the
transcendental $\zeta(2n-3)$.}
\end{figure}

So we have the beautiful correspondence {\em $\zeta(2n-3)$ in the
Feynman graphs $\leftrightarrow$ $(2,2n-3)$ torus knot in the link
diagram.}
Fig.(\ref{hab87})
lists some Feynman diagrams, link diagrams, braid words for the knots
in them, and numbers obtained from explicit calculation.
We stress that we have transcendentals which are single sums,
and braid words which rely on a single generator.
The power of the generator matches the exponent in the sums.
This is the first class of Feynman diagrams whose transcendentals
are described by a simple knot class. The counterterms of Feynman diagrams
which have the topology  of a crossed ladder -Fig.(\ref{hab86})- deliver 
countertems which have the form $\sim \zeta(2n-1)/\epsilon$,
and deliver the $(2,q)$ torus knots from the link diagrams.
The coefficients of the transcendental depend on the particles
realizing the topology. This opens an interesting problem in
representation theory: how do the rational coefficients
of $\zeta(2n-3)$ change with the particle content of the diagram?
This is under investigation at the moment.
\clearpage
\subsection{The $(3,4)$ torus knot and the first Euler double sum}
We have classified ladder (no knots) and crossed-ladder
topologies ($(2,q)$ torus knots). In the next step
we generalize this to more complicated topologies. We will see that at the
six-loop level for the first time a new knot appears,
and that at this level for the first time a new transcendental 
is obtained. It has the form of an Euler double sum,
and thus is a generalized $\zeta$-function of two integer
arguments. We identify the
new knot with the new transcendental. From now on we 
denote such transcendentals as knot-numbers, 
with $\zeta(q)$ being the knot-numbers of the $(2,q)$ torus knots.


An important observation is that all the subdivergence-free 
diagrams considered here and in the following allow
for a Hamiltonian circuit representation.
This means that we can find a closed non-singular curve (no
self-intersections)
in the diagram which traverses all internal vertices.
In general, not every three-valent graph
allows for a Hamiltonian circuit, but the failure appears at loop
orders and topologies which are not relevant for our purposes
\cite{tutte}.
In fact, this statement is only true for graphs without subdivergences. 
For graphs with subdivergences
it is very easy to construct examples which do not allow a 
Hamiltonian circuit.

\begin{figure}[ht]\epsfxsize=5in \epsfbox{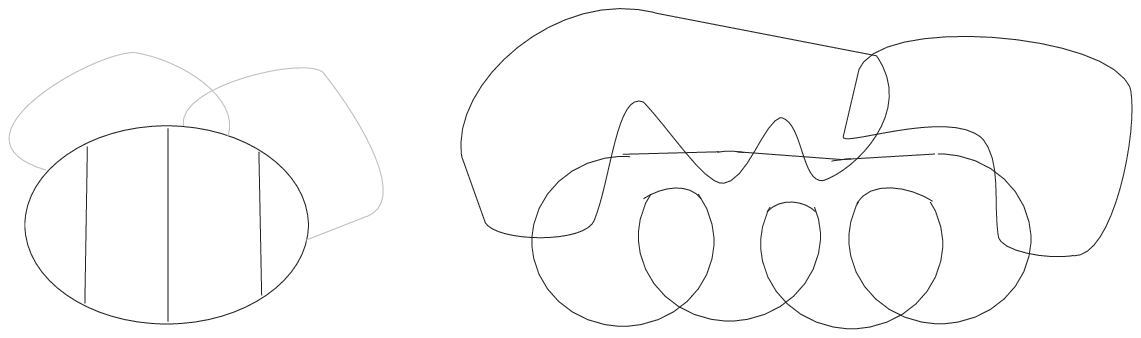}
\caption[How to obtain link diagrams.]{\label{hab89}\small 
A systematic method to obtain link diagrams.
This method proved successful in the investigation of
subdivergence-free diagrams. For this example, we would get three $RMI$
moves on the rhs for free.} 
\end{figure}

Having identified such a Hamiltonian circuit it was used in
\cite{plb} to simplify the search for corresponding 
link diagrams and knots:
assume
we draw all propagators in the interior of the Hamiltonian circle.
Now let us remove as many propagators as necessary to make the diagram
planar (it reduces to a ladder topology).
Using our standard momentum routing for this reduced diagram, let us
begin to attach the non-planar propagators again, this time always using
Reidemeister III moves to avoid crossings with the inner 
propagators of the
Hamiltonian circuit. These inner propagator constitute a ladder topology,
which is then disturbed by any further propagator which makes it
into a non-trivial topology. Assume we have $r$ such inner propagators.
Skeining the $r$ components of the ladder, it is clear that
we get $r-1$ Reidemeister I moves for free, removing $r-1$ crossings
in the link diagram. There remains a link diagram with $n-r$ components,
which still has to be skeined. In  Fig.(\ref{hab89}) we demonstrate
the method, which was also used in the second link
diagram in Fig.(\ref{hab86}).


We conclude that as far as the knot content of the link diagram is
concerned
the planar rungs are irrelevant.
By chosing the orientations of the
components in an appropriate manner we will only obtain
positive knots with this method. Topologically
non-trivial propagators are drawn in the exterior of the Hamiltonian
circuit. They are connected to vertices in the most economic way,
avoiding as many crossings as possible. 


But in multi-loop calculations one finds occasionally
new transcendentals, independent from $\zeta(i)$.
A first and prominent example is the transcendental found
by David Broadhurst in a six-loop calculation at transcendentality
level 8 \cite{david}. 
According to our experience with $\zeta$-transcendentals,
we would expect a knot with 8 crossings to appear in the
corresponding Feynman graph, which is confirmed in  
Figs.(\ref{hab811},\ref{hab812}).
The corresponding $\phi^4$ graph was mapped to a graph with three-point
couplings by using a Lagrange multiplier field. From the
three ways how to map a four-point coupling to a three-point coupling
only those possibilities were considered which did not change
the topology of the graph. 
\begin{figure}[ht]\epsfysize=3cm \epsfbox{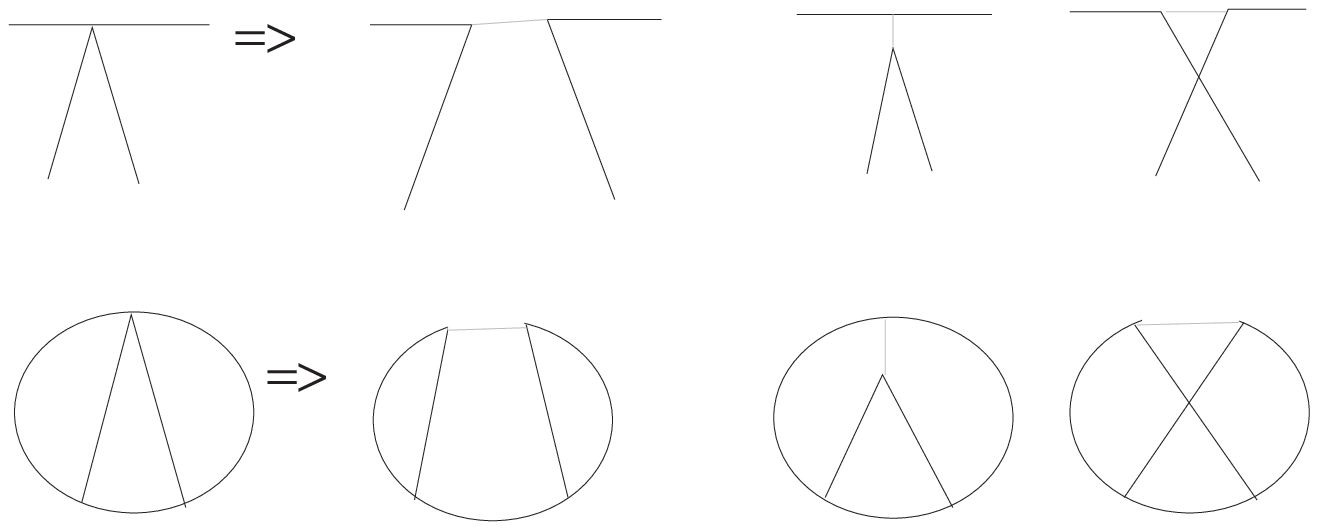}
\caption[From 4pt to 3pt coupling.]{\label{hab810}\small
There are three ways to map a four point coupling to
two three-point couplings. Some change the topology
of the corresponding link diagram. The results exclude the
choices which generate extra crossings. In the example,
only the first one gives the correct answer, 
the other two generate extra crossings in the link diagram.}
\end{figure}
The first graph which gave a non-$\zeta$
transcendental is the six-loop graph in Fig.(\ref{hab811}).
\begin{figure}[ht]
\epsfysize=4cm \epsfbox{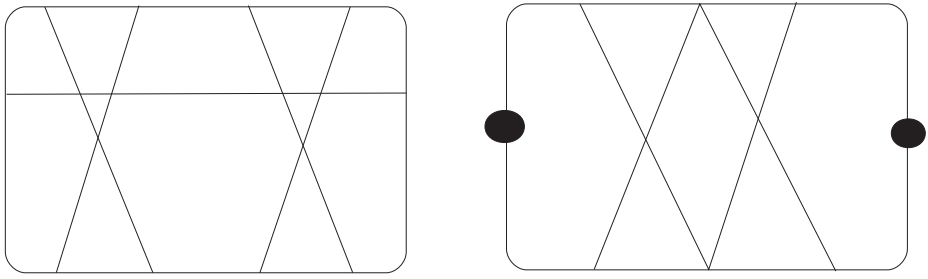}
\caption[A six-loop graph.]{\label{hab811}\small We like to investigate this
six-loop Feynman graph. We also give a $\phi^4$ graph which
is topologically equivalent, and which was investigated by 
Broadhurst. 
The two dots in this
graph have to be identified. It can be obtained
from the graph on the lhs by shrinking three propagators.}
\end{figure}
Let us map this graph to a link diagram as in Fig.(\ref{hab812}) to
identify the $(3,4)$ torus knot.
\begin{figure}[ht]
\epsfysize=8cm \epsfbox{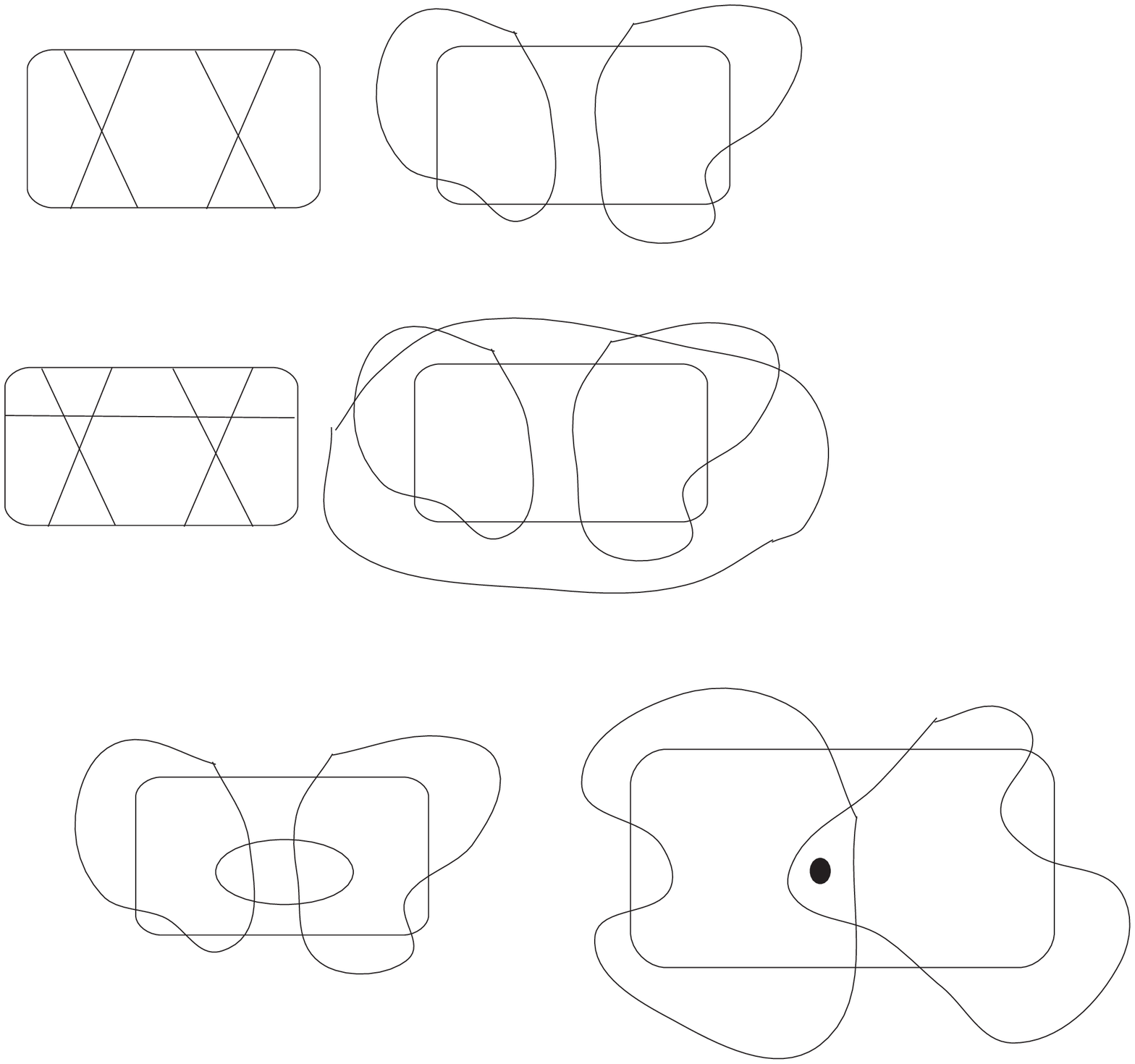}
\caption[The $(3,4)$ torus knot appears.]{\label{hab812}\small 
The generation of the $(3,4)$ torus knot.
In the first line we have removed one propagator to generate the 
$\zeta(3)\zeta(3)$ factor knot. Then we attach the last 
propagator in the most economic way, giving us the link
diagram
on the bottom rhs. We used Reidemeister II and III
moves to get from the second to the third line.
We end up with the
braid word $\sigma_1^4\sigma_2\sigma_1^4\sigma_2$.
(All components encircle the dot in the middle
counterclockwise, so that we
can read off the braid word.)
After skeining the two kidneys we find a knot.
It can be identified as
the $8_{19}$ knot in the standard tables \cite{knots}, which is
the $(3,4)$ torus knot.}\end{figure}
The identification of the $(3,4)$ torus knot is achieved
by reading off the braid word
\beas
\sigma_1\sigma_2^3\sigma_1\sigma_2^3=
(\sigma_1\sigma_2)^4
\eeas in Fig.(\ref{hab812}).
In general, the $(p,q)$ torus knot has braid word \cite{jones}
\beas
(\sigma_1\ldots\sigma_{p-1})^q.
\eeas
We expect only
positive braid words to appear. Positive braid words
have only positive powers of braid generators.
The positive braid words up to nine crossings
are the $(2,q), q\in \{3,5,7,9\}$ torus knots plus the
$(3,4)$ torus knot. Crossing
number nine is the transcendentality level nine
which is exhausted by graphs up to six loops \cite{david}.
An investigation of the results in \cite{david}
confirmed this pattern. All six-loop graphs considered
deliver knots which faithfully describe the transcendentals
in their counterterms.


Further, as the $(3,4)$ torus knot is the only
non-$\zeta$-ish transcendental at level 8, we conclude that
our knot theoretic approach predicts a relation between
the value of the transcendental $M$ which Broadhurst reported
in \cite{david} and the level 8 transcendental $U_{6,2}$ which was
found in the expansion of the master function \cite{exp}.
Such a relation was meanwhile established and is given in \cite{pisa}.

Meanwhile, number theorists proved that the new
knot number is an independent transcendental indeed \cite{db}. 
Again there is a striking match between
the structure of the possible braid word 
and the structure of the transcendental.
The transcendental comes as a double sum, where the two exponents
match the powers of the two braid generators:
\begin{equation}
\mbox{braid word}\;\; \sigma_1 \sigma_2^3 \sigma_1 \sigma_2^3
\Leftrightarrow \sum_{i_1 <i_2}^\infty \frac{1}{i_1^2 i_2^6}.
\end{equation}
Other representations of this knot have other 8-crossing braid words,
but then other level 8 double sums are equivalent to the one given.
This is so as one can prove that there is only one independent
level 8 double sum.
\clearpage
\subsection{$\phi^4$ theory: more knots and numbers}
This section reports on a major work, which initiated
most of the subsequent activities in the field.
To test the ideas developed so far a calculation of the $\beta$-function in
$\phi^4$ theory was undertaken. For us it is 
sufficient to mention that the $\beta$-function
is a function which is determined by the
divergences of the Feynman graphs. In fact, in the
MS renormalization scheme, it solely relies on the
divergences of the graph. It determines
the dependence of the coupling constant on the mass scale.\footnote{In
DR regularized quantities depend on a scale parameter $\mu$ of mass
dimension, see the remarks in the appendix.} 
Being calculable in terms of divergences, the $\beta$-function
is a quantity
well-suited to test our ideas.


Having satisfyingly matched all six-loop graphs with
the correct knot-numbers, the process was inverted to obtain
a method to calculate diagrams analytically, when they were only
available numerically. 
By mapping Feynman diagrams to knots, and by matching knots
to knot-numbers, we obtain a conjecture for the transcendentality
content of a graph.


Let us say we suspect that a graph $G$ should have the
form $G=a_1 \zeta_1 +a_2 \zeta_2$, for knot-numbers $\zeta_1,\zeta_2$ and
rational numbers $a_1,a_2$. The transcendentals
$\zeta_1,\zeta_2$ are delivered by the
link diagram associated with the graph, and $a_1,a_2$ are yet to be
determined. Then, by having a numerical value for the graph $G_{num}$,
we can try to determine the rational numbers $a_1,a_2$ by solving
$G_{num}=a_1 \zeta_1+a_2 \zeta_2$, such that $a_1,a_2$ are simple
rational numbers, and then check if the solutions remain unchanged
if we increase the accuracy of $G_{num}$. Numerical as well as 
number theoretic methods give a precise answer and error estimate
for
the validity of this approach.


Exploring this method, the authors of \cite{pisa}
calculated all the primitive divergences contributing to the 7--loop
$\beta$\/--function of $\phi^4$ theory, i.e.\ all 59 diagrams
that are free of subdivergences and hence give scheme--independent
contributions. 
The method delivered analytic results
for 56 out of 59 diagrams. 
Three diagrams, associated
with the knots $10_{124}$, $10_{139}$, and $10_{152}$, were only
obtained numerically to 10 significant figures. The accuracy
was not sufficient to apply the method in a unique manner.\footnote{
For the 56 graphs for which the method was successful numerical accuracy was
such that the chance of an error is $<10^{-30}$.}
Only one hyperbolic knot with 11 crossings was encountered and the
transcendental number associated with it was found. It is
the only positive knot with 11 crossings other than the torus knot 
$(2,11)$.


In addition, the series of `zig--zag' counterterms, 
$\{6\zeta(3)$, $20\zeta(5)$, 
$\frac{441}{8}\zeta(7)$, $168\zeta(9),\,\ldots\}$, 
was obtained.
This series is the equivalent of the crossed ladder
topology in the $\phi^4$ theory
and thus should produce $(2,2n-3)$ torus knots, 
corresponding to $\zeta(2n-3)$.
And indeed, it was found that  
the $n$\/--loop zig--zag term is $4C_{n-1}
\sum_{p>0}\frac{(-1)^{p n - n}}{p^{2n-3}}$,
where $C_n=\frac{1}{n+1}{2n \choose n}$
are the Catalan numbers, familiar in knot theory.


At seven loops, transcendentals of degree $\leq 11$
can occur, and thus we expect knots with up to eleven
crossings. 
With the exceptions of the hyperbolic knots
$10_{139}$, $10_{152}$, and $11_{353}$, all
positive knots with up to 11 crossings
are torus knots.\footnote{In \cite{pisa} these three hyperbolic knots
were erroneously denoted as satellite knots.} 
Apart from the $(2,q)\;q\in \{3,5,7,9,11\}$ torus knots
there are two more torus knots encountered at up to 7 loops:
$8_{19}=(4,3)$ and $10_{124}=(5,3)$. 
We know already that $8_{19}$
corresponds to a double sum whose structure matches
the structure of the braid word, and for $(3,5)$ a similar result
was obtained.


Also it was found that 
the triple--sum $N_{3,5,3}\equiv\sum_{l>m>n>0}l^{-3}m^{-5}n^{-3}$
is associated with the knot $11_{353}$. Again the corresponding knot
had a braid word which matches: $11_{353}\sim \sigma_1 \sigma_2^2 \sigma_1^2
\sigma_3^2\sigma_2^3\sigma_3$. This  knot has a braid word in three 
different generators, with the total powers of the generators matching
the exponents in the sum.

\begin{figure}[ht]\epsfysize=5cm\epsfbox{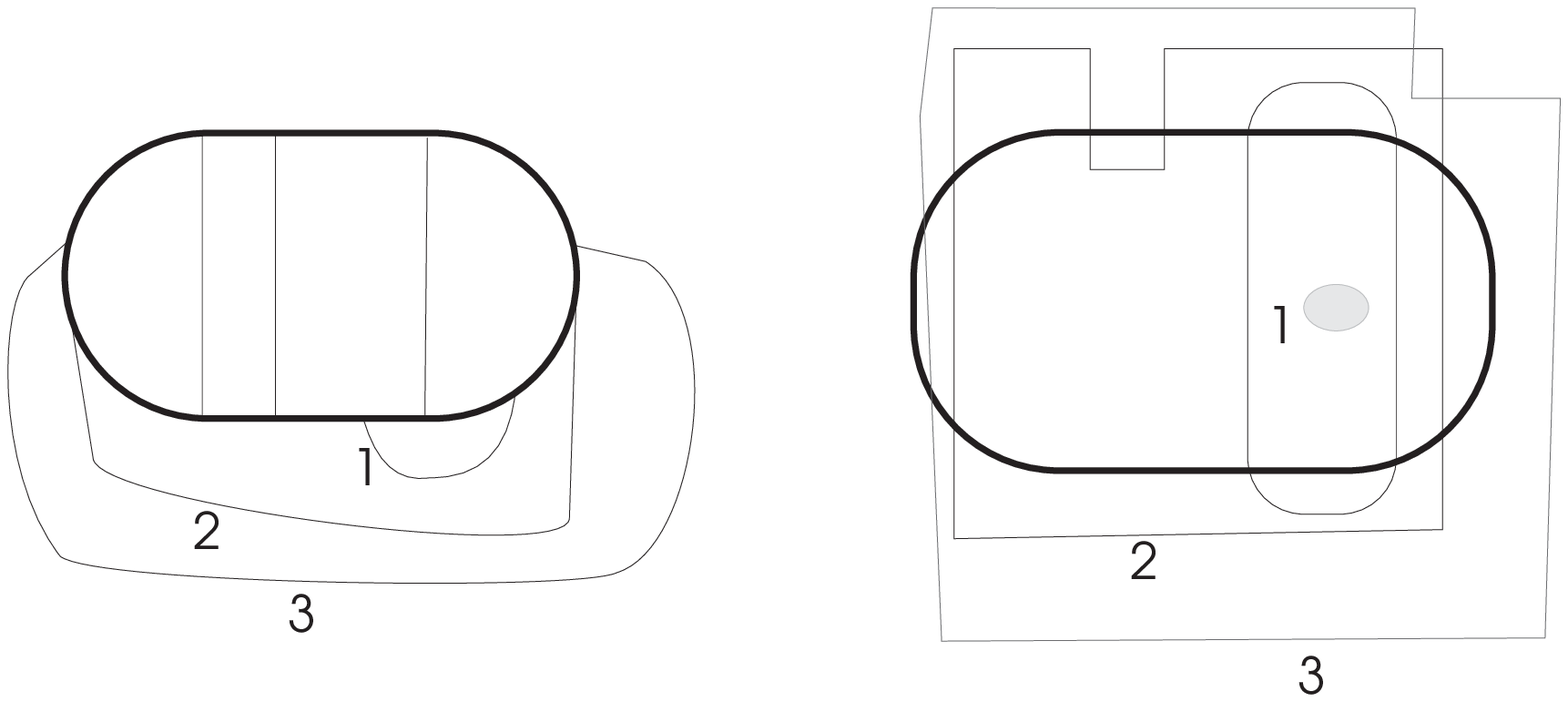}
\caption[$11_{353}$.]{\label{hab813}\small 
The identification of the knot $11_{353}$.
There are also $\zeta(3)^2\zeta(5)\sim (2,3)^2(2,5)$ factor knots
as well as the $\zeta(11)=(2,11)$ torus knot in this diagram.
We only give one possible link diagram which delivers
$11_{353}$. The analytic result for this diagram 
contains $N_{3,5,3}$ as well as
$\zeta(11)$ and $\zeta(3)^2\zeta(5)$. It can be
found in \cite{pisa}, where this diagram belongs to class $C7$.}
\end{figure}
\begin{figure}[ht]\epsfysize=8cm\epsfbox{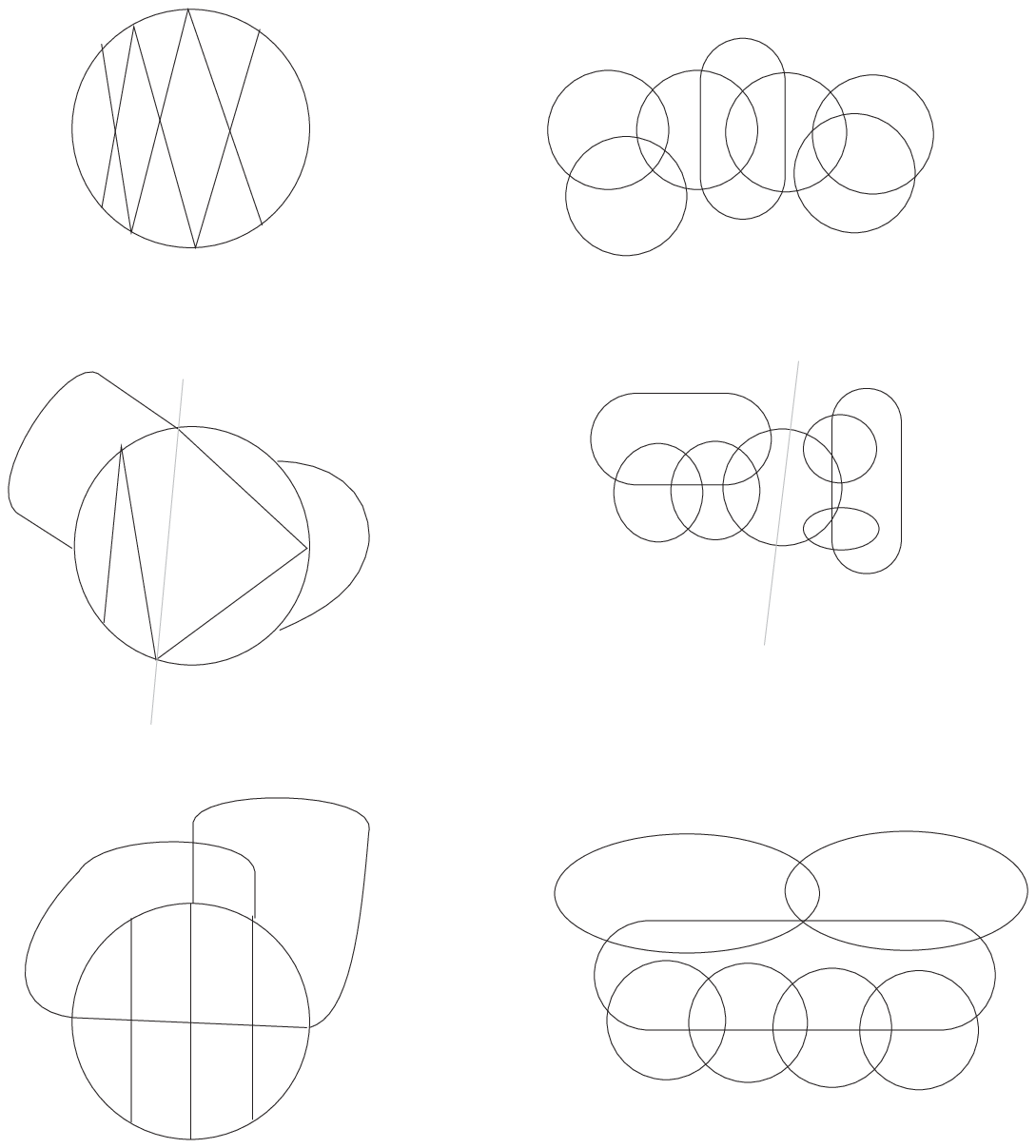}
\caption[Factor knots.]{\label{hab814}\small
These are easy factor knots. From top to bottom,
they deliver $(2,3)^3\sim \zeta(3)^3$,
$ (2,5)^2\sim \zeta(5)^2$, $(2,3)\times (2,7)\sim \zeta(3)\zeta(7)$.}
\end{figure}

Fig.(\ref{hab813}) gives an example for the appearance of this knot.
We give a Feynman graph which was obtained from a $\phi^4$
graph by replacing four-point couplings
by three-point couplings, as in Fig.(\ref{hab810}).
Fig.(\ref{hab814}) gives further examples for  simple factor knots
at the seven loop level.


The knot $10_{124}=(5,3)$ 
was also obtained in \cite{pisa} and identified with the transcendental
double sum $U_{7,3}\equiv N_{7,3}$ \cite{db}.
The reader will find a list of all results in \cite{pisa},
and further details and generalizations will figure in \cite{book}.
\clearpage
\subsection{Rationality and the $\beta$-function of quenched $QED$}
In this paper, again a $\beta$-function was considered.
This time it is the $\beta$-function of QED. In fact, the restricted
class of diagrams to be considered were 
diagrams containing only a single fermion 
loop, which constitute the so-called quenched $\beta$-function. 
The aim was
to give a knot-theoretic explanation for the rationality of the quenched
QED $\beta$ function. The following paragraphs are more or less
a direct quote from \cite{bdk}, where the reader can find details
of the calculation, as well as a confirmation of the mechanism proposed here
for usual QED 
in scalar QED, which is the theory of a photon coupling to scalar particles.


In QED one has a gauge symmetry. On the quantized level
this is a powerful tool which establishes relations between
various Feynman diagrams. Here we only need to know
that the divergences of the fermion propagator
match the divergences of the vertex-function,
which is expressed as $Z_1=Z_2$, where $Z_1$ and
$Z_2$ are the $Z$-factors of the QED
Lagrangian, cf.~Eq.(\ref{aqed}). We refer to the identity $Z_1=Z_2$ as the 
Ward identity.\footnote{The original Ward identity relates the fermion 
propagator and the vertex function as $S^{-1}(p_1)-S^{-1}(p_2)=
(p_1^\mu-p_2^\mu)\Gamma_\mu(p_1,p_2)$, which implies $Z_1=Z_2$.}
At the link level, the Ward identity entails
cancellation of subdivergences generated by one term of the skein relation,
which in turn implies cancellation of knots generated by the other term. 

\begin{figure}[ht]\epsfysize=12cm \epsfbox{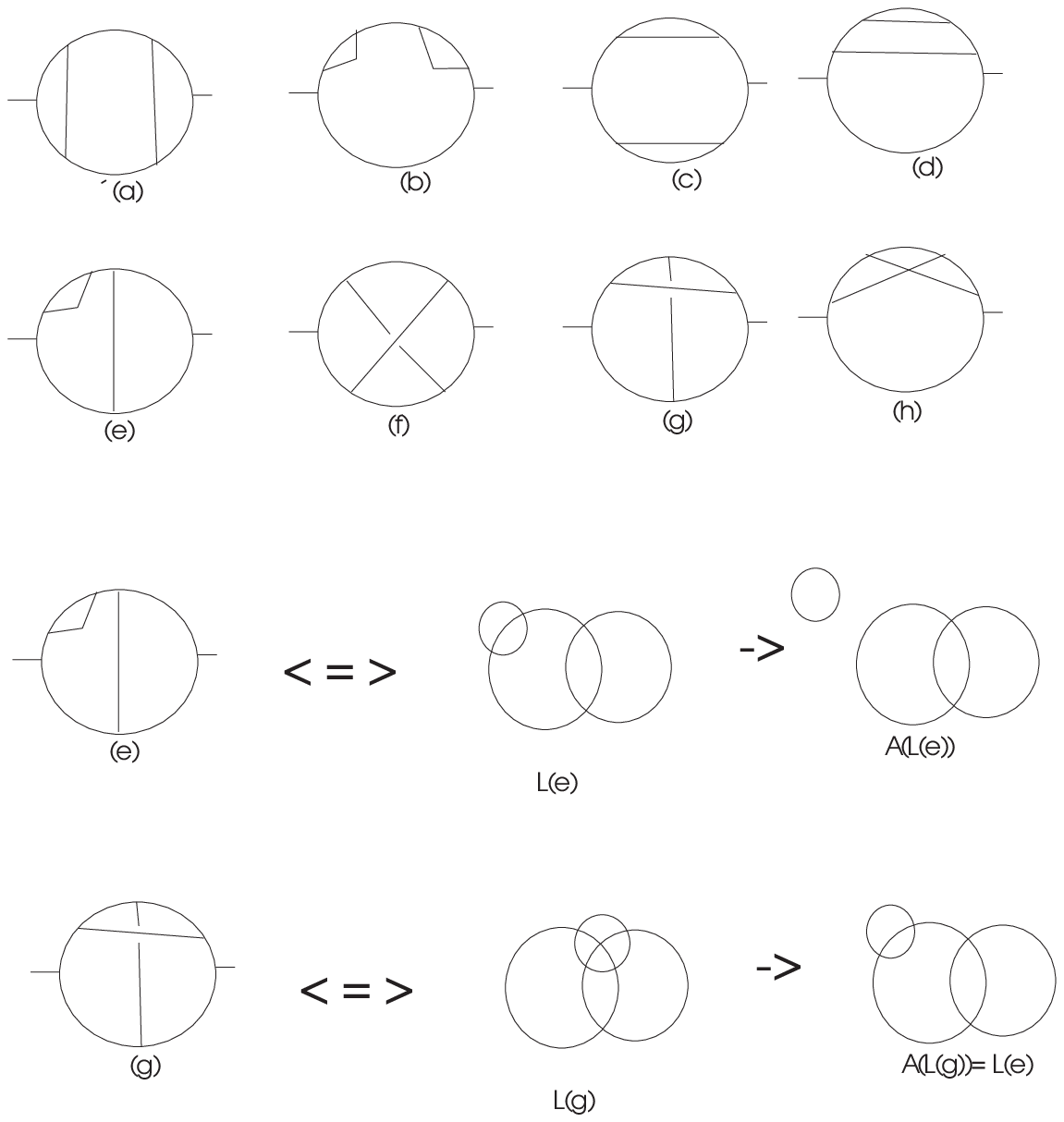}
\caption[The vanishing of transcendentals for the QED 
$\beta$-functions.]{\label{hab819}\small
This figure summarizes the argument in
\cite{bdk} which explains the absence of transcendentals in the
quenched QED $\beta$-function.}\end{figure}

Rationality of the three-~\cite{JLR} and
four-loop~\cite{4LQ} quenched (i.e.\ single-electron-loop) terms in the QED
beta function comes as a surprise. Feynman diagrams with non-trivial knots
contribute to it, so that the vanishing of transcendentals demands
an explanation.
For example the crossed-photon graphs of 
Figs.(\ref{hab819}f,g,h) all realize the link diagram
whose skeining contains the trefoil knot.
With the experience of the previous subsections in mind we expect
to find $\zeta(3)/\epsilon$ in their divergent parts. 
To explain the cancellation
of transcendentals, we must study the interplay between knot-theoretic
arguments and the gauge structure of QED.


We propose to associate the cancellation of transcendentals with the
cancellation of subdivergences in the quenched beta function of QED, which
is an immediate consequence of the Ward identity, $Z_1=Z_2$.


To see the key role of the Ward identity, consider Fig.(\ref{hab819}g). 
There is an
internal vertex correction, which is rendered local by adding the
appropriate counterterm graph. Due to the Ward identity, this counterterm
graph is the same as that which compensates for the self-energy correction
in Fig.(\ref{hab819}e). We know that the latter counterterm
could be interpreted as the $A$ part of the skein operation on
the link diagram of Fig.(\ref{hab819}e), associated with the Feynman
graph. We assume that this is a generic feature of the relationship
between skeining and renormalization and associate the corresponding
counterterm for Feynman graph Fig.(\ref{hab819}g) 
with the term obtained from applying
$A$ twice to its link diagram, which
requires two skeinings to generate the same
counterterm, along with the trefoil knot from the $B$ term.
The Ward identity thus becomes a relation between
crossed and uncrossed diagrams, after skeining:
\begin{equation}
A(A(L({\rm g})))=A(L({\rm e}))\quad\Rightarrow\quad
A(L({\rm g}))=L({\rm e})\,,
\label{Ward}
\end{equation}
where $L(e),L(g)$ refer to the appropriate link diagrams in
Fig.(\ref{hab819}).


The trefoil knot results then from
$B(B(L({\rm g})))$, which generates, in general, a $\zeta(3)/\epsilon$ term
from Fig.(\ref{hab819}g), even after the subtraction of subdivergences.
We now use the Ward identity Eq.(\ref{Ward}) at the link
level to obtain
\begin{equation}
B(B(L({\rm g}))) = B(B(A^{-1}(L({\rm e}))))\,,
\label{arg}
\end{equation}
from which we see that it relates the transcendental
counterterm from a torus-knot topology, in Fig.(\ref{hab819}g), to a
knot-free ladder topology, in Fig.(\ref{hab819}e).
From section seven we know
that ladder topologies are free of
transcendentals when the appropriate counterterms are added: 
after subtraction of subdivergences, ladder graphs, such
as in Figs.(\ref{hab819}a,e), give rational terms in the Laurent expansion in 
powers of
$1/\epsilon$. 


We conclude that $\zeta(3)/\epsilon$ 
should be absent from the bare diagram of
Fig.(\ref{hab819}g) when it is calculated in the gauge where the bare
diagram of Fig.(\ref{hab819}e) is free of $\zeta(3)/\epsilon$, 
i.e.\ in the Landau
gauge, where the latter is free of subdivergences.\footnote{In gauge
theories, the local gauge symmetry allows such welcome accidents:
that certain graphs may be free of divergences in a special gauge.} 
In Fig.(\ref{hab819}) we summarize
this argument. 
By a similar argument, we also conclude that the other graphs
with the trefoil topology, namely Figs.(\ref{hab819}f,h), should be free of
$\zeta(3)/\epsilon$ in the Landau gauge. Which was indeed
confirmed by explicit calculation in \cite{bdk}.
\clearpage
\subsection{Euler double sums and the expansion in the
two-loop master function} 
\begin{figure}[ht]\epsfysize=6cm \epsfbox{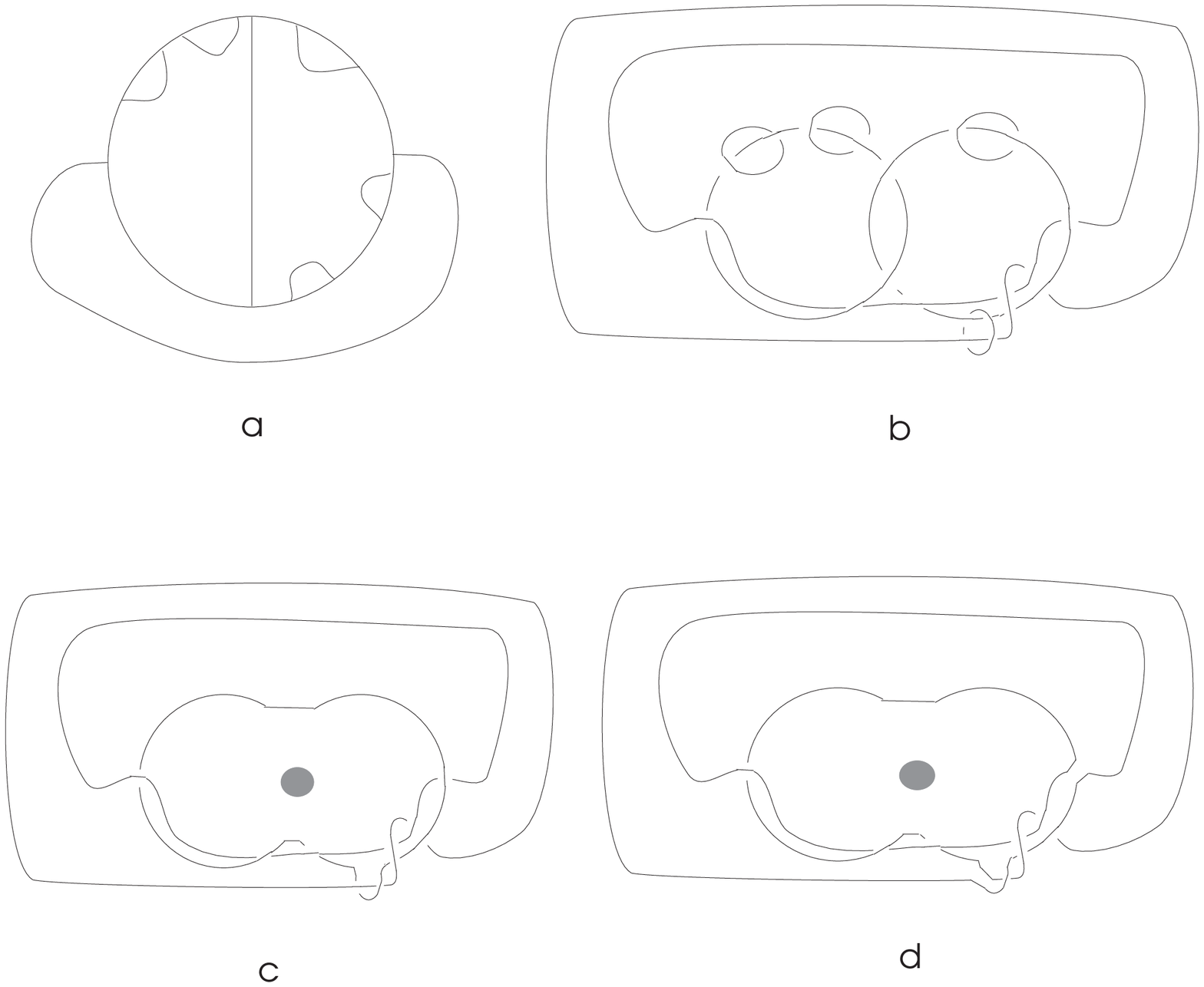}
\caption[Dressed trefoils graphs.]{\label{hab820}\small 
These graphs which come
from dressing the basic trefoil topology generate double sums
and correspond to three-braids on the other hand. The treatment of the bubbles coincide
with the experience of section five and the results in
\cite{newplb}.}\end{figure}

Here the investigation of knots and numbers was pushed further to level 12.
The diagrams considered are the dressed master function,
which incorporates one-loop bubbles in the basic trefoil
topology of Fig.(\ref{hab820}). The results are in agreement
with the results of section five, and with the upcoming
results in \cite{newplb}.

\begin{figure}[ht]\epsfxsize=5in\epsfbox{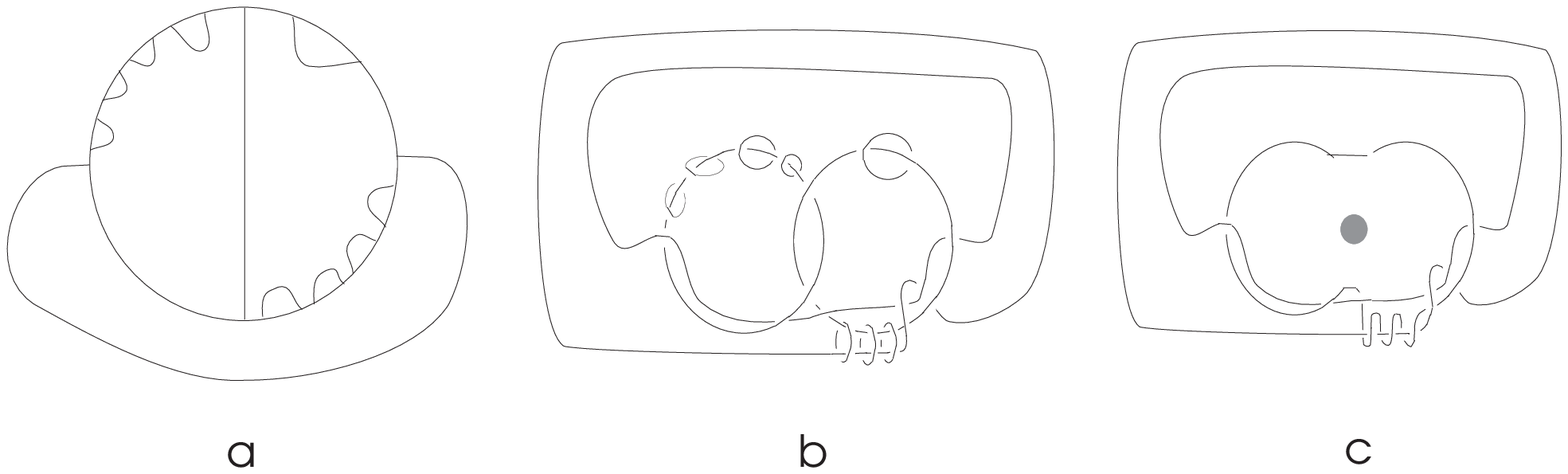}
\caption[A 12 crossing example.]{\label{hab821}An example
for the characteristic knot class for double sums.}\end{figure}

At first sight the analysis in \cite{bgk} seemed to be in conflict
with
the number theory in~\cite{BBG}.
In Fig.(\ref{hab821}) we demonstrate the appearance of twelve
crossing knots in the diagrams considered in \cite{bgk}.


The
`rule of 3' discovered in~\cite{BBG}, for
non-alternating sums $\zeta(a,b)=\sum_{n>m}1/n^a m^b$
of level $l=a+b$, demands,
at even level $l=2p+2$, as much as $\lfloor p/3\rfloor$ 
irreducibles.\footnote{$\lfloor x\rfloor$ is the nearest integer
$\leq x$.}
On the other hand the number of positive three-braid knots
generated from topologies as in Fig.(\ref {hab820})
grows by a rule of two,
giving $\lfloor p/2\rfloor$ knots at $l=2p+4$ crossings.
So there are two 12-crossing knots, while \cite{BBG}
only have one level 12 irreducible. 


This problem initiated a major finding in number theory
by David Broadhurst \cite{db}, which we will comment
in the next section. For now, we directly quote from there:
"Faced with a 12-crossing knot in search of a number,
we saw two ways to turn: to study 4-fold non-alternating
sums, or 2-fold sums with alternating signs.
The first route is numerically intensive: it soon emerges
that well over 100 significant figures are needed to find
integer relations between 4-fold sums at level 12.
The second route
is analytically challenging; it soon emerges that
at all even levels $l\geq6$ there are relations between alternating double
sums that cannot be derived from any of the identities
given in~\cite{BBG}.


Remarkably, these two routes lead, eventually,
to the {\em same\/} answer.
The extra 12-crossing knot is indeed associated with
the existence of a 4-fold
non-alternating sum, $\zeta(4,4,2,2)$ $=$ $\sum_{n>m>p>q}1/n^4m^4p^2q^2$,
which {\em cannot\/} be reduced to non-alternating sums
of lower levels. It is, equivalently, associated
with the existence of an irreducible
{\em alternating\/} double sum,
$U_{9,3}=\sum_{n>m}$$[(-1)^n/n^9]$$[ (-1)^m/m^3]$.
The equivalence stems from the unsuspected circumstance that the combination
$\zeta(4,4,2,2)-(8/3)^3 U_{9,3}$, and only this combination,
is reducible to non-alternating double sums.
Moreover, $l=12$ is the lowest level at which the reduction
of non-alternating 4-fold sums necessarily entails an alternating
double sum. The `problem pair' of knots are a problem no more.
Their entries in the knot-to-number dictionary~\cite{book}
record that they led to a new discovery in number theory:
{\em the reduction of non-alternating sums necessarily entails
alternating sums}."



The family of positive knots~\cite{bgk}
that gave rise to this investigation has
braid
words
\begin{equation}
\{\sigma_1\sigma_2^{2a+1}\sigma_1\sigma_2^{2b+1}\mid a\geq b\geq1\}\,,
\label{braid}
\end{equation}
whose enumeration matches that of the irreducible
$\{U_{2a+3,2b+1}\mid a\geq b \geq 0\}$.\footnote{Conventions
as in \cite{db}. These irreducibles furnish all the transcendentals
associated with the family of three-braids in Eq.(\ref{braid})
which was  confirmed up to level $a+b\leq 20$ in \cite{db}.}
Fig.(\ref{hab821}) gives an example for the level 12 members of this family.
The observation that more bubbles will only generate the
knot family given in Eq.(\ref{braid}), together with the reduction of the
quantity $I(\mu)$ to a double sum in \cite{bgk}, initiated
the findings in \cite{db}.
\clearpage
\subsection{General Euler sums}
We mentioned already Broadhurst´s investigation of the
number of irreducible Euler sums \cite{db}.
Guided by knot theory, he finds \cite{db}
"a generating function for the number, $E(l,k)$, of irreducible
$k$-fold Euler sums, with all possible alternations of sign, and exponents
summing to $l$. Its form is remarkably simple:
$\sum_n E(k+2n,k)\,x^n=\sum_{d|k}\mu(d)\,(1-x^d)^{-k/d}/k$,
where $\mu$ is the M\"obius function. Equivalently, the size
of the search space in which $k$-fold Euler sums of level $l$
are reducible to rational linear combinations of irreducible basis terms
is $S(l,k)=\sum_{n<k}{\lfloor(l+n-1)/2\rfloor\choose n}$."


Irreducible Euler sums are found to occur in explicit analytical
results, for counterterms with up to 13 loops, yielding transcendental
knot-numbers, up to 23 crossings.
In \cite{db} we can find a dramatic account
how these results were obtained, supported by and support
the connection between knot theory and field theory as 
proposed here.
\subsection{Identities demanded by knot theory}
In this work we consider results from one and two-loop bubble
insertions in various Feynman diagrams which were iterated from
the first term in the skeleton expansion. We obtain
identities reflecting the demands of link theory and show
that the diagrams fulfil these identities. Some basic results
are already given in section five, and are extended in
\cite{newplb} to higher loop orders. The findings confirm the
analysis of \cite{bgk}.
\subsection{Field Theory, Knot Theory, Number Theory}
From the results of \cite{plb,pisa,bgk,db} it
is clear that for the case of field-theory counterterms up to $L$ loops,
the irreducible Euler sums
\begin{eqnarray}
&&\{\zeta{(2a+1)}\mid L-2\geq a>0\}\,,\label{2b}\\
&&\{U_{2a+1,2b+1}\mid a>b>0,\ L-3\geq a+b\}\,,\label{3b}\\
&&\{\zeta(2a+1,2b+1,2c+1)\mid a\geq b\geq c>0,\ a>c,\ L-3\geq a+b+c\}
\label{4b}\,,
\end{eqnarray}
appear \cite{db}.\footnote{$\zeta(2a+1,2b+1,2c+1):=\sum_{n>m>r}\frac{1}{n^{2a+1}
m^{2b+1}r^{2c+1}}\equiv N_{2a+1,2b+1,2c+1}$.}


It is not yet known whether Euler sums exhaust the
transcendentals in counterterms at $L\geq7$ loops. In addition to
$10_{124}
\simeq U_{7,3}$, there are two further
positive knots with 10 crossings, namely~\cite{jones}
$10_{139}$ and
$10_{152}$, which
7-loop analysis~\cite{pisa} suggests are {\em not}
associated with Euler sums of depth $k<4$. 
All those, and {\em only\/} those,
subdivergence-free 7-loop $\phi^4$ diagrams whose link diagrams skein
to $10_{139}$ and $10_{152}$ appear, on the basis of numerical evidence,
to give counterterms that cannot be reduced to Euler sums
with depth $k\leq3$. Hence it is an open (and fascinating) question
whether they belong to the class of Euler sums at all or provide different
transcendentals.
Table(\ref{t3}) gives a summary of the achievements so far.
It summarizes the work of Broadhurst and Kreimer
and is taken from \cite{db}.
\begin{table}
\caption[Positive prime knots and corresponding numbers.]{\label{t3}\small
Positive prime knots related to Euler sums,
via field-theory counterterms.}
\[\begin{array}{|l|l|l|}\hline
\mbox{crossings}&\mbox{knots}&\mbox{numbers}\\[3pt]\hline
2a+1
&(2,2a+1)=\sigma_1^{2a+1}&\zeta(2a+1)\\\hline
8
&(3,4)=8_{19}=\sigma_2\sigma_1^3\sigma_2\sigma_1^3 & N_{5,3}\\\hline
9
&\mbox{none}&\mbox{none}\\\hline
10
&10_{124}=(3,5)=\sigma_2\sigma_1^3\sigma_2\sigma_1^5&N_{7,3}\\
&10_{139}=\sigma_2\sigma_1^3\sigma_2^3\sigma_1^3&?\\
&10_{152}=\sigma_2^2\sigma_1^2\sigma_2^3\sigma_1^3&?\\\hline
11
&\sigma_1^2\sigma_2^2\sigma_1\sigma_3\sigma_2^3
\sigma_3^2&N_{3,5,3}\\\hline
12
&\fk{}7{}3&N_{9,3}\\
&\fk{}5{}5&N_{7,5}-\frac{\pi^{12}}{2^5\tm10!}\\
&\fk{}353&?\\
&\fk{}335&?\\
&\fk2235&?\\
&\fk2334&?\\
&\fk3333&?\\\hline
13
&\sk{}32{}42&N_{3,7,3}\\                                    
&\sk2{}3223&N_{5,3,5}\\                                     
&\sk{}3{}{}52&?\\                                         
&\sk2{}33{}3&?\\                                          
&\sk{}3{}233&?\\                                          
&\so2\st2\so{}(\se{}\st3)^2&?\\                           
&(\st{}\so{}\se{}\st{})^3\so{}&?\\                        
&(\st{}\so{}\se{}\st{})^3\st{}&?\\                        
\hline\end{array}\]
\end{table}


Let us collect the successes of the empirical association of
knots and numbers, initiated by the calculation of counterterms in
renormalizable field theories.
The Euler sum $\zeta(2L-3)$, corresponding to the 2-braid
torus knot $(2L-3,2)$, first appears in anomalous dimensions or
$\beta$-functions at $L$
loops \cite{pisa}. 
Only further symmetries like gauge symmetries might guarantee its
absence \cite{bgk}. 
This is the only transcendental which  
results from subdivergence-free diagrams
with less than 6 loops, because no other positive knot has
less than 8 crossings.
On the other hand at 6 loops 3-braids start to appear \cite{plb}. 
The first of these is
$8_{19}\simeq N_{5,3}$. The irreducibility of its knot-number
was confirmed in~\cite{BBG}.
Further at 7 loops 4-braids start to appear with the knot
$11_{353}$. The irreducibility of its knot-number
was confirmed in~\cite{BG},
following communication of its appearance in field theory~\cite{pisa}.
From there on, the findings in \cite{bgk} and \cite{db}, supported
by \cite{newplb}, inform number theory.
The reader will find further comments and results in the quoted papers
as well as in \cite{book}.
\clearpage

\clearpage
\section{Conclusions}
In this paper we established a useful algebra to facilitate 
calculations for topological simple Feynman graphs. Remarkably, this
algebraic
structure could be established whenever the topology was of ladder type,
so that we could use the same one-loop algebra for nested as well
as overlapping divergences. Starting from this algebra,
we proved in section seven a striking connection between
topological simplicity and a number theoretic property.
In section eight we reported
results which strongly suggest a yet to be explored
connection between field theory, knot theory and number theory.
When musing about the results achieved a few obvios things come to mind.
The Euler sums which so pertinently appear in the counterterms of 
renormalizable field theories have relations to quantum groups and
recent results in topology. For example in \cite{Kassel} 
we can find them
in relation to the Drinfeld associator and Kontsevich´s work.
This opens fascinating questions to be explored in the future.
Also, our use of the Hamiltonian circuit brings to mind
chord diagrams. Any three-valent graph which allows for a
Hamiltonian circuit can be directly related to chord 
diagrams.\footnote{In \cite{Jarvis} the authors tried to explore this
idea to map Feynman diagrams to knots using singular
knot projections. Unfortunately, they fail. Their "connectivity rules"
where strictly designed to get the correct results up to the three loop 
level, where they
reproduce $\zeta(3)\sim$ trefoil. But their rules give wrong
results even for the four loop example. Contrary to the hope the authors express 
in their paper, the incorporation of the resolution of singular
points in terms of over- and undercrossings {\em a la
Vassiliev} \cite{vass} cannot cure their results,
which one sadly but easily confirms by explicit calculation.}
This indicates that one should explore four-term relations for their
usefulness in classifying and simplifying the divergences of 
a pQFT. This is under investigation at the moment.
We hope that in the future it might be possible not only to identify the
transcendentals
in the $Z$-factor by considering the associated link diagram, but also
hope to be able to predict the actual coefficients of these
transcendentals,
when the theory under investigation is specified. This would enable us 
to calculate  $Z$-factors to arbitrary loop order just from 
a knowledge of the
topology
of the graph and knowledge of its particle content. 
\section*{Acknowledgements}
It is a pleasure to thank Bob Delbourgo,
Peter Jarvis and Ioannis Tsohantjis. The material presented here in
sections 3-7 was obtained during a stay at the University of Tasmania.
Without their patient interest and tenacious questions, this paper would probably
have never been written. With equal pleasure I thank the whole theory
group for a wonderful time down there.

I can not imagine how the ideas presented here would have turned 
into the rocksolid results of \cite{plb,pisa,bdk,bgk} without
the possibility to collaborate with David Broadhurst.  It was and
is a most stimulating collaboration.

It is a pleasure to thank Lou Kauffman for his interest and hospitality during a 
stay at the University of Illinois at Chicago, and for good advice in
preparation of this manuscript.

I thank Roger Picken for organizing "Knots and Physics V" (Porto 1996),
which was a most interesting and stimulating conference.

Thanks are due to Kostja Chetyrkin,
Lance Dixon, Hugh Osborne, Jon Rosner, Bas Tausk and
George Thompson for discussions and interest.

I thank Karl Schilcher and the theory group in Mainz for a supporting
and stimulating atmosphere.

This work was supported by DFG.

\begin{appendix}
\section{Conventions}
In explicit results given in this paper we usually omit
trivial factors of $2^D\pi^{D/2}$, typically
generated by Eqs.(\ref{dr1},\ref{dr2}), with the only exception
being the results for $\phi^3_6$ given in section two, where we 
obey the conventions of \cite{Collins}.

We also omit to give explicitly the dependence
on the renormalization scale $\mu$. Such a scale results from
the demand to have a coupling constant which has no dimension
in DR, so that one redefines the coupling 
$g^2\rightarrow g^2(\mu^2)^{\epsilon}$, where  
$\mu$ has the dimension of a mass. Renormalized quantities
are independent of $\mu$.
\section{Feynman Rules}
We gave the Feynman rules for $\phi^4$ in four dimensions
in section two.
Here we give the Feynman rules for some other theories.
\subsection{$\phi^3_6$}
There is a  modification in the symmetry factor. The Lagrangian
comes now with a coupling $g/3!$. The vertex $-ig$ connects to three
propagators. The loop integration is with respect to a measure
\begin{equation}
\int \frac{d^6k}{(2\pi)^6}
\end{equation}
which becomes in DR $\int \frac{d^Dk}{(2\pi)^D}$, with $D=6-2\epsilon$.
One has three $Z$-factors as in $\phi^4$, one for each monomial
in the Lagrangian. 
This theory is renormalizable in six dimensions. All other theories 
considered  in this paper are renormalizable in four dimensions.
\subsection{QED}
In QED one finds the following set of rules:
\begin{itemize}
\item For each internal electron line, associate a propagator given by:\\
$iS_F(p)=\frac{i}{\pslash-m+i\eta}=\frac{i(\pslash+m)}{p^2-m^2+i\eta}$
\item For each internal photon line (in the Feynman gauge), associate a propagator:\\
$iD_F(p)_{\mu\nu}=-\frac{ig_{\mu\nu}}{p^2+i\eta}$
\item At each vertex place a factor of:\\
$-ie\gamma_\mu$
\item For each internal loop, integrate over:\\
$\int \frac{d^4q}{(2\pi)^4}$
\item 
Incoming electrons have momenta entering the diagram, while incoming
positrons have momenta leaving the diagram. Internal fermion lines appear with
arrows in both clockwise and counterclockwise direction. However, topologically
equivalent diagrams are counted only once. 
\end{itemize}
We see that one can almost read off these Feynman rules from the Lagrangian of
$QED$:
\begin{equation}
{\cal L}=\bar{\psi}(i\gamma^\mu(\partial_\mu+ieA_\mu)-m)\psi-\frac{1}{4}F_{\mu\nu}F^{\mu\nu}
-\frac{1}{2}(\partial_\mu A^\mu)^2\label{QED}
\end{equation}
and the propagators above appear as the inverse of the equation of motion
for the noninteractiong electron and photon fields.
The $Z$-factors are as follows.
\begin{equation}
{\cal L}=\bar{\psi}(i\gamma^\mu(Z_2\partial_\mu+iZ_1 eA_\mu)-Z_mm)\psi-
\frac{1}{4}Z_3 F_{\mu\nu}F^{\mu\nu}
-\frac{1}{2}Z_a(\partial_\mu A^\mu)^2\label{aqed}
\end{equation}
\subsection{Scalar QED}
This is the theory of an abelian gauge theory coupling to charged
scalar particles.
We only give the modifications compared to QED.
\begin{itemize}
\item At each scalar-scalar-vector vertex place a factor of:\\
$-ie (p+p^\prime)_\mu$, where $p$ and $p^\prime$ are the
momenta for the scalar line.
\item Insert a factor of $2ie^2g_{\mu\nu}$ for each four-point
vector-vector-scalar-scalar vertex.\\
\end{itemize}
The Lagrangian is
\begin{equation}
{\cal L}=(\partial_\mu-ieA_\mu)\phi^\dagger (\partial_\mu-ieA_\mu)\phi
-m^2\phi\phi^\dagger. 
\end{equation}
\subsection{Yukawa theory}
This is the coupling of fermions with scalar particles.
We have an interaction term $g\phi \bar{\psi}\psi$ and a $Z$-factor
$Z_g$ for it with vertex $-ig$,
as well as the usual free Lagrangians of scalar and fermionic particles.
Yukawa theory provides convenient examples. For example the three-point one-loop
vertex function becomes (zero-mass limit, zero momentum transfer)
\begin{equation}
\int d^Dk \frac{1}{\kslash\kslash(q+k)^2}=\int d^Dk\frac{1}{k^2(k+q)^2},
\end{equation}
and thus provides a very simple example for a one-loop function
$\Delta$.
\section{Form Factors}
\subsection{A $QED$ example}
The situation is most easily explained by studying a example.
Consider the case of a one-loop vertex correction in massless QED at zero
momentum transfer:
\beas
\int d^Dl \;\frac{1}{
(l+k)^2}\gamma_\rho\frac{1}{\lslash}\gamma_\mu\frac{1}{\lslash}
\gamma_\rho=\\
(\Delta_{11}\gamma_\mu+\Delta_{21}\kslash k_\mu
)(k^2)^{-\varepsilon},
\eeas
which allows for two form factors 
$\Delta_{11}$ and $\Delta_{21}$. Only $\Delta_{11}$ is UV divergent.

Iterating this in nested topologies both form factors
may replace the tree level vertex $\gamma_\mu$. So we also have to
consider
functions where the vertex structure is $\kslash k_\mu$ for some
momentum $k$, and this gives us two more $\Delta$ functions:
\beas
\int d^Dl \;\frac{(l^2)^{-\varepsilon i}}{
(l+k)^2}\gamma_\rho\frac{1}{\lslash}\gamma_\mu\frac{1}{\lslash}
\gamma_\rho=\\
(k^2)^{-\varepsilon (i+1)}[{}_i\!\Delta_{11}\gamma_\mu
+{}_i\!\Delta_{21}\frac{\kslash k_\mu}{k^2}],\\
\int d^Dl \;\frac{(l^2)^{-\varepsilon i}}{
(l+k)^2}\gamma_\rho\frac{1}{\lslash}\frac{\lslash l_\mu}{l^2}
\frac{1}{\lslash}=\\
(k^2)^{-\varepsilon (i+1)}[{}_i\!\Delta_{12}\gamma_\mu
+{}_i\!\Delta_{22}\frac{\kslash k_\mu}{k^2}].
\eeas
At the one-loop level we have
\beas
Z_1^{(1)}=<\Delta_{11}\gamma_\mu+\Delta_{21}\frac{\kslash
k_\mu}{k^2}>=<\Delta_{11}>\gamma_\mu,
\eeas
by definition.
At the next order we find 
\beas
Z_1^{(2)}=<(\Delta_{11}\;\; {}_1\!\Delta_{11}+\Delta_{21}
\;\; {}_1\!\Delta_{12})\gamma_\mu\\
+
(\Delta_{11} \;\;{}_1\!\Delta_{21}+\Delta_{21}
\;\; {}_1\!\Delta_{22})\frac{\pslash p_\mu}{p^2}>\\
-<<\Delta_{11}>(\Delta_{11}\gamma_\mu+\Delta_{21}\frac{\pslash
p_\mu}{p^2})>\\
=<(\Delta_{11}\;\; {}_1\!\Delta_{11}+\Delta_{21}
\;\; {}_1\!\Delta_{12})\gamma_\mu>-<<\Delta_{11}>\Delta_{11}\gamma_\mu>
\eeas
as $<\Delta_{21}\Delta_{22}>=0$.
We used
\beas
<\Delta_{11}\;\; {}_1\!\Delta_{21}>-<<\Delta_{11}>\Delta_{21}>=0,
\eeas
by Weinberg's theorem (we subtracted out the subdivergence from an
overall convergent form factor 
$\sim \pslash p_\mu)$.
UV-divergences only appear in the form factor $\sim \gamma_\mu$, as it
should be. But we see that, at higher order, contributions of the
other form factor mix into the result. Nevertheless, the concatenation
works when we take this mixing iteratively into account.
The above considerations are easily generalized to more complicated
cases
involving other form factors, $r$ say. We would obtain
an $r \times r$ matrix of $\Delta$ functions
$\Delta_{ij}$, which concatenates in the way outlined above.

The same considerations apply to the two-point case, and the results
of the example in the text concerning the two-loop fermion self-energy
were  obtained along these lines.
\subsection{Quadratic divergences}
Occasionally one has to confront quadratic divergences.
The gauge boson propagator in a gauge theory provides
such an example.
We describe how to handle these quadratic divergences.
It is sufficient to consider the case of a massless theory.
As an example we take Fig.(\ref{hab63}) for this case,
and we assume that this topology is realized as a correction
to the photon propagator in $QED$.
The quadratic divergent term is of the form 
\begin{equation}
\frac{Tr(\lslash\gamma_\mu\lslash\gamma_\rho\kslash\gamma_\nu\kslash
\gamma^\rho)}{l^2(l+q)^2(l-k)^2k^2(k-q)^2},
\end{equation} 
Using Eq.(\ref{dr3}) we subtract zero
\begin{eqnarray*}
& & \frac{Tr(\lslash\gamma_\mu\lslash\gamma_\rho\kslash\gamma_\nu\kslash
\gamma^\rho)}{l^2(l+q)^2(l-k)^2k^2(k-q)^2}\\
& & -\frac{Tr(\lslash\gamma_\mu\lslash\gamma_\rho\kslash\gamma_\nu\kslash
\gamma^\rho)}{l^2(l)^2(l-k)^2k^2(k)^2}\\
& & =
\frac{Tr(\lslash\gamma_\mu\lslash\gamma_\rho\kslash\gamma_\nu\kslash
\gamma^\rho)(l^2k^2-(l+q)^2(k-q)^2)}{l^4(l+q)^2(l-k)^2k^4(k-q)^2}\\
& & =\frac{Tr(\lslash\gamma_\mu\lslash\gamma_\rho\kslash\gamma_\nu\kslash
\gamma^\rho)(l^2(q^2-2k.q)+2l.q(k^2-2k.q+q^2)
+q^2(k^2-2k.q+q^2))}{l^4(l+q)^2(l-k)^2k^4(k-q)^2},
\end{eqnarray*} 
and no term is worse than linear divergent by powercounting,
and thus amenable to the usual treatment.
\section{Generation of $\zeta(3)$}
We want to calculate the divergent part of the 
three-loop $\phi^4$ graph in Fig.(\ref{hab81}).
We use a variant of the Gegenbauer polynomial x-space technique 
\cite{smirnov,chet}.
We use Gegenbauer polynomials in four dimensions, Chebyshev
polynomials that is.

The analytic expression for the graph is
\begin{equation}
I:=\int d^Dk d^Dl d^Dr \frac{1}{r^2l^2k^2}\frac{1}{(l-k)^2(r-l)^2(k-r)^2}.
\label{i1}
\end{equation}
We use an expansion of a propagator in Gegenbauer polynomials
$C^1_n(x)$.
\begin{equation}
\frac{1}{(p-q)^2}= \sum_{n=0}^{\infty}\frac{1}{\mbox{max}(p,q)^2}
\mbox{min}(\frac{p}{q},\frac{q}{p})^n C^1_n(\hat{p}.\hat{q}).\label{gexpa}
\end{equation}
The $C^1_n(x)$ fulfil
\begin{eqnarray}
& & C_n^1(1)=n+1,\;\int d\hat{q} C_{n_1}^1(\hat{x}.\hat{q})C_{n_2}^1(\hat{y}.
\hat{q})=\frac{1}{n_1+1}\delta_{n_1,n_2}C_{n_1}^1(\hat{x}.\hat{y}),
\nonumber\\
& & C_0^1(x)=1,\; \int d\hat{q}=1,\label{i2}
\end{eqnarray}
where $p:=\sqrt{p^2}$, $\hat{p_\mu}:=p_\mu/p$, $\hat{p_\mu}\hat{p^\mu}=1$.

The integral is totally symmetric in $l,r,k$.
Therefore we consider only one sector, $l<r<k$, say. We have six such
sectors, and the total result is just six times the result for one
such sector.
We are only interested in the divergent part of $I$,
which appears when all three loop momenta become large.
So in our chosen sector we can take as integration measure
\begin{eqnarray*}
& & I=C\int_1^\infty dk^2 \;dl^2\;dr^2\;
(l^2)^{-\epsilon}(r^2)^{-\epsilon}(k^2)^{-\epsilon}
\;\int d\hat{k}d\hat{l}d\hat{r}
\sum_{n_1,n_2,n_3}\frac{1}{k^2l^2r^2}\\
& & \frac{1}{\mbox{max}(l,r)^2}
\mbox{min}(\frac{l}{r},\frac{r}{l})^{n_1} C^1_{n_1}(\hat{l}.\hat{r})\\
& & \frac{1}{\mbox{max}(l,k)^2}
\mbox{min}(\frac{l}{k},\frac{k}{l})^{n_2} C^1_{n_2}(\hat{l}.\hat{k})\\
& & \frac{1}{\mbox{max}(k,r)^2}
\mbox{min}(\frac{k}{r},\frac{r}{k})^{n_3} C^1_{n_3}(\hat{k}.\hat{r}).
\end{eqnarray*}
$C$ incorporates trivial factors as the volume of the
four-dimensional unit-sphere.
Using Eqs.(\ref{i2})
we can do the angular integrations easily.
Then we evaluate the radial integrations in our specified  sector
\begin{eqnarray}
& & I=C\int_1^\infty dk^2 \;dl^2\;dr^2\;
(l^2)^{-\epsilon}(r^2)^{-\epsilon}(k^2)^{-\epsilon}\;
\frac{1}{r^2}\frac{1}{k^4}\sum_{n=0}^\infty\frac{1}{n+1}
\left[\frac{l^2}{k^2}\right]^{n}\\
& & \sim\;\sum_{n=0}^\infty \frac{1}{(n+1)^3}\sim \zeta(3),
\end{eqnarray}
in DR. Here we used that we can do the angular integrations in four dimensions.
One easily verifies that the result in $D$ dimensions only modifies
the finite part. We also used that in DR all radial integrations
vanish at infinity. Inserting all factors one confirms the well-known
result $6\zeta(3)$.
\section{A Regulator Function}
We want to define a function ${}_j\!\bar{\Delta}$ which fulfils
Eq.(\ref{rat}) and renders ${}_j\!\Delta$ finite by powercounting.
This is easily achieved by constructing a function
which has the same behaviour for large loop momentum.
Let $F(k)$ be the integrand of the generalized one-loop
function ${}_j\!\Delta$. Then $\Theta(k^2-1)F(k)$
provides an integrand which has the same asymptotic behaviour for large $k$.
We then can always define a function $G(k)$ so that 
\begin{equation}
{}_j\!\bar{\Delta}:=\int d^Dk [\Theta(k^2-1)F(k)+\Theta(1-k^2)G(k)],
\end{equation}
has the desired properties.
\section{Cancellations of Transcendentals}
We give an example for the concatenated product
$B^6(\Delta)$, which corresponds to a seven loop calculation,
compared with $(-A+B)^6(\Delta)$, for the example
\beas
\Delta=\int d^D k \; \frac{1}{k^2(k+q)^2}\mid_{q^2=1}.
\eeas
We also give the lower order terms, as it is quite interesting to see
the conspiracy of rationals appearing. So in the following
$G(r)$ is a $r+1$ loop Green function, and $Z(r)$ the corresponding
counterterm expression. `ge' is the Euler
$\gamma$ and `zet(i)' means $\zeta(i)$
and $x:=D-4$. 
We do not use a renormalization which would absorb the
$\gamma$ and $\zeta(2)$, and also do not use that
$\zeta(2)^3$, $\zeta(4)\zeta(2)$ and $\zeta(6)$ are all
dependent ($\sim \pi^6$), as we want to exhibit the generated
rationals in their purest form. 
{\small
\begin{verbatim}


         1   -1    1   -2
Z(1):=---*x   - ---*x
         2         2


       -1           5      1   -2
G(1):=x  *( - ge + ---) + ---*x
                    2      2


         2   -1    1   -2    1   -3
Z(2):=---*x   - ---*x   + ---*x
         3         2         6


       -1      1            3    2    9        55
G(2):=x  *( - ---*zet(2) + ---*ge  - ---*ge + ----)
               4            4         2        6


          -2      1        3      1   -3
       + x  *( - ---*ge + ---) + ---*x
                  2        2      6


         5   -1    19   -2    1   -3    1    -4
Z(3):=---*x   - ----*x   + ---*x   - ----*x
         4         24         4         24


       -1      23            1               7            4    3
G(3):=x  *( - ----*zet(3) + ---*zet(2)*ge - ---*zet(2) - ---*ge
               9             3               6            9


               14    2    125        455
            + ----*ge  - -----*ge + -----)
               3           6         12


          -2      1             1    2    7        125
       + x  *( - ----*zet(2) + ---*ge  - ---*ge + -----)
                  12            3         3        24


          -3      1        7       1    -4
       + x  *( - ---*ge + ----) + ----*x
                  6        12      24


         14   -1    19   -2    11   -3    1    -4     1    -5
Z(4):=----*x   - ----*x   + ----*x   - ----*x   + -----*x
         5          12         24         12         120


       -1      133            335               335
G(4):=x  *( - -----*zet(4) + -----*zet(3)*ge - -----*zet(3)
               96             72                18


                5         2    25           2    25
            + -----*zet(2)  - ----*zet(2)*ge  + ----*zet(2)*ge
               192             96                12


               245            125    4    125    3    1225    2
            - -----*zet(2) + -----*ge  - -----*ge  + ------*ge
               48             576         36           48


               595        6727      -2      67
            - -----*ge + ------) + x  *( - ----*zet(3)
                6          40               72


             5                5             25     3    25    2
          + ----*zet(2)*ge - ----*zet(2) - -----*ge  + ----*ge
             48               12            144         12


             245        119
          - -----*ge + -----)
             24          6


          -3      1             5     2    5        49
       + x  *( - ----*zet(2) + ----*ge  - ---*ge + ----)
                  48            48         6        24


          -4      1         1       1    -5
       + x  *( - ----*ge + ---) + -----*x
                  24        6      120


           -1    1313   -2    47   -3    25    -4    1    -5
Z(5):=7*x   - ------*x   + ----*x   - -----*x   + ----*x
                 360          48         144         48


             1    -6
         - -----*x
            720


       -1      512            183               1647
G(5):=x  *( - -----*zet(5) + -----*zet(4)*ge - ------*zet(4)
               75             80                160


               23                   23           2    207
            + ----*zet(3)*zet(2) - ----*zet(3)*ge  + -----*zet(3)*ge
               30                   5                  5


               4991            3         2       27         2
            - ------*zet(3) - ----*zet(2) *ge + -----*zet(2)
                45             80                160


               3            3    81           2    217
            + ----*zet(2)*ge  - ----*zet(2)*ge  + -----*zet(2)*ge
               20                40                20


               903             9     5    81    4    217    3
            - -----*zet(2) - -----*ge  + ----*ge  - -----*ge
               40             100         40         10


               2709    2    19369        62601      -2
            + ------*ge  - -------*ge + -------) + x  *(
                20           40           80


             61             23               69
          - -----*zet(4) + ----*zet(3)*ge - ----*zet(3)
             160            15               10


              1         2    3            2    27
          + -----*zet(2)  - ----*zet(2)*ge  + ----*zet(2)*ge
             160             40                40


             217            3     4    27    3    217    2    903
          - -----*zet(2) + ----*ge  - ----*ge  + -----*ge  - -----*ge
             120            40         20         20          20


             19369      -3      23            1
          + -------) + x  *( - ----*zet(3) + ----*zet(2)*ge
              240               90            40


             9             1     3    27    2    217        301
          - ----*zet(2) - ----*ge  + ----*ge  - -----*ge + -----)
             80            20         40         60         40


          -4       1             1     2    9         217
       + x  *( - -----*zet(2) + ----*ge  - ----*ge + -----)
                  240            40         40        360


          -5       1         3        1    -6
       + x  *( - -----*ge + ----) + -----*x
                  120        80      720


         132   -1    277   -2    839   -3    19   -4     7    -5
Z(6):=-----*x   - -----*x   + -----*x   - ----*x   + -----*x
          7          30          360         48         144


             1    -6     1     -7
         - -----*x   + ------*x
            240         5040


       -1      20641            49567               49567
G(6):=x  *( - -------*zet(6) + -------*zet(5)*ge - -------*zet(5)
               4320             3600                 720


               1687                   11809           2
            + ------*zet(4)*zet(2) - -------*zet(4)*ge
               5760                   5760


               11809               28679            102487        2
            + -------*zet(4)*ge - -------*zet(4) + --------*zet(3)
                576                 480             12960


               5929                      5929
            - ------*zet(3)*zet(2)*ge + ------*zet(3)*zet(2)
               4320                      864


               41503           3    41503           2
            + -------*zet(3)*ge  - -------*zet(3)*ge
               12960                 864


               100793               44891             49          3
            + --------*zet(3)*ge - -------*zet(3) - -------*zet(2)
                360                  72              34560


                343         2   2    343         2       833        2
            + -------*zet(2) *ge  - ------*zet(2) *ge + -----*zet(2)
               11520                 1152                960


               2401            4    2401           3
            - -------*zet(2)*ge  + ------*zet(2)*ge
               34560                1728


               5831           2    2597               49049
            - ------*zet(2)*ge  + ------*zet(2)*ge - -------*zet(2)
               480                  48                 480


               16807     6    16807    5    40817    4    18179    3
            + --------*ge  - -------*ge  + -------*ge  - -------*ge
               518400         17280         2880           144


               343343    2    116039        4753177      -2
            + --------*ge  - --------*ge + ---------) + x  *(
                480             48           1260


             7081            1687               1687
          - ------*zet(5) + ------*zet(4)*ge - ------*zet(4)
             3600            2880               576


             847                    5929           2
          + ------*zet(3)*zet(2) - ------*zet(3)*ge
             4320                   4320


             5929               14399             49         2
          + ------*zet(3)*ge - -------*zet(3) - ------*zet(2) *ge
             432                 360             5760


              49         2    343            3    343           2
          + ------*zet(2)  + ------*zet(2)*ge  - -----*zet(2)*ge
             1152             8640                576


             833               371            2401     5    2401    4
          + -----*zet(2)*ge - -----*zet(2) - -------*ge  + ------*ge
             240               48             86400         3456


             5831    3    2597    2    49049        16577      -3
          - ------*ge  + ------*ge  - -------*ge + -------) + x  *(
             720           48           240          48


             241             847                847
          - ------*zet(4) + ------*zet(3)*ge - -----*zet(3)
             2880            2160               432


              7          2     49            2    49
          + ------*zet(2)  - ------*zet(2)*ge  + -----*zet(2)*ge
             5760             2880                288


             119             343     4    343    3    833    2
          - -----*zet(2) + -------*ge  - -----*ge  + -----*ge
             240            17280         864         240


             371        7007      -4      121
          - -----*ge + ------) + x  *( - ------*zet(3)
             24         240               2160


              7                  7              49     3    49     2
          + ------*zet(2)*ge - -----*zet(2) - ------*ge  + -----*ge
             1440               288            4320         288


             119        53
          - -----*ge + ----)
             120        24


          -5       1               7      2     7         17
       + x  *( - ------*zet(2) + ------*ge  - -----*ge + -----)
                  1440            1440         144        120


          -6       1          1        1     -7
       + x  *( - -----*ge + -----) + ------*x
                  720        144      5040


\end{verbatim}}
\end{appendix}

\listoffigures
\listoftables

\begin{thebibliography}{99}
\bibitem{schweber}
S.~Schweber, {\em QED and the men who made it: Dyson, Feynman,
Schwinger,
Tomonaga}  Princeton Univ.~Press (1994), and references there.
\bibitem{thooft}
G.'t Hooft, Nucl.Phys.B33 (1971) 173; Nucl.Phys.B35 (1971) 167.
\bibitem{Lee}
B.W.~Lee, in {\em Methods in Field Theory}, ed.~R.~Balian, 
J.~Zinn-Justin (1976) North-Holland.
\bibitem{Dy}
F.J.~Dyson, Phys.Rev.75 (1949) 486, 1736.
\bibitem{Salam}
A. Salam, Phys.Rev.82 (1951) 217, 84, 426.
\bibitem{Weinberg}
S.~Weinberg, Phys.Rev.118 (1960) 838.
\bibitem{BPHZ}
N.N.~Bogoliubov, D.V.Shirkov, {\em Introduction to the theory of
quantized fields}, Wiley (1980).
\bibitem{Zimm}
W.~Zimmermann, Comm.Math.Phys.15 (1969) 208.
\bibitem{org}
D.~Kreimer, {\em Renormalization and Knot Theory},
Univ.~of Tasmania preprint UTAS-PHYS-94-25 (1994).
This preprint contains most of the material of sections 3,4,6,7 here.
Following a suggestion of Lou Kauffman, we added section 2.
Section 5 is new.
Section 8 reviews recent results and contains a few new results as well.
\bibitem{Kaku}
M.Kaku, {\em QFT}, Oxford UP 1993.
\bibitem{Collins}
J.Collins, {\em Renormalization}, Cambridge UP 1984.
\bibitem{knots}
C.C.~Adams, {\em The Knot Book}, Freeman (1994).\\
D.~Rolfsen, {\em Knots and Links}, Publish or Perish (1976).\\
L.H.~Kauffman, {\em Knots and Physics}, World Scientific (1991).
\bibitem{plb}
D.~Kreimer, Phys.Lett.B354 (1995) 117.
\bibitem{pisa}
D.J.~Broadhurst, D.~Kreimer, Int.J.Mod.Phys.C6 (1995) 519.
\bibitem{bdk}
D.J.~Broadhurst, R.~Delbourgo, D.~Kreimer, Phys.Lett.B366 (1996) 421.
\bibitem{bgk}
D.J.~Broadhurst, J.~Gracey, D.~Kreimer, Z.Phys.{\bf C}75 (1997) 559;\\
see also D.J.~Broadhurst, A.V.~Kotikov, {\em
Compact analytical form for non-zeta terms in critical exponents
at order $1/N^3$}, hep-th/9612013.
\bibitem{db}
D.J.~Broadhurst, {\em On the enumeration of irreducible
$k$-fold Euler sums and their roles in knot theory and field theory},
to appear in J.Math.Phys., OUT-4102-62 (hep-th/9604128).
\bibitem{newplb}
D.~Kreimer, {\em On Knots in subdivergent Diagrams},
to appear in Z.Phys.{\bf C}, MZ-TH/96-31 (hep-th/9610128).
\bibitem{LSZ}
H.~Lehmann, K.~Symanzik, W.~Zimmermann,
Nuovo Cim.1, (1955) 1425.
\bibitem{Haag}
R.Haag, {\em Local Quantum Physics}, Springer 1992.
\bibitem{exp}
D.J.~Broadhurst, Z.Phys.C32 (1986) 249;\\
D.T.~Barfoot, D.J.~Broadhurst, Z.Phys.C41 (1988) 81.
\bibitem{mpla}
D.~Kreimer, Mod.Phys.Lett.A9 (1994) 1105.
\bibitem{smirnov}
V.A.~Smirnov, {\em Renormalization and Asymptotic Expansions},
Birkh\"auser (1991).
\bibitem{chet}
K.G.~Chetyrkin, A.L.~Kataev, F.V.~Tkachov, Nucl.Phys.B174 (1980)
345, see also \cite{smirnov} and references there.
\bibitem{Itz}
C.~Itzykson, J.-B.~Zuber, {\em Quantum Field Theory},
McGraw-Hill (1980).
\bibitem{bobnew}
R.~Delbourgo, A.~Kalloniatis, G.~Thompson, Phys.Rev.{\bf D54}
(1996) 5373;\\
R.~Delbourgo, D.~Elliott, D.S.~McAnally, 
Phys.Rev.{\bf D55} (1997) 5230.
\bibitem{ussu}
V.V.~Belokurov, N.I.~Ussyukina, J.Phys.A16 (1983) 2811.\\
A.I.~Davydychev, N.I.~Ussyukina, Phys.Lett.B298 (1993) 363;\hfill
Phys.Lett.B305 (1993) 136.\\
D.J.~Broadhurst, Phys.Lett.B164 (1985) 356; Phys.Lett.B307 (1993) 132.
\bibitem{david}
D.J.~Broadhurst, {\em Massless Scalar Feynman Diagrams: Five Loops
and Beyond}, (1985) Open Univ.~preprint OUT-4102-18.
\bibitem{tutte}
D.Barnette, {\em Map Coloring, Polyhedra, and the Four-Color Problem},
The Dolcani Mathematical Expositions (1983).
\bibitem{jones}
V.~Jones, Annals of Math.126 (1987) 335.
\bibitem{book}
D.~Kreimer, {\em Knots and Feynman diagrams}, Cambridge UP, in preparation.
\bibitem{JLR}
J.L.~Rosner, Phys.Rev.Lett.17, (1966) 1190; Ann.Phys.44 (1967) 11.
\bibitem{4LQ}
S.G.~Gorishny, A.L.~Kataev, S.A.~Larin, L.R.~Surguladze,
Phys.Lett.B256 (1991) 81.
\bibitem{BBG}
D.~Borwein, J.M.~Borwein, R.~Girgensohn, Proc.Edin.Math.Soc.38 (1995)
273.
\bibitem{BG}
J.M.~Borwein, R.~Girgensohn, CECM research report 95-053 (1995),
Electronic J.~Combinatorics {\bf 3} (1996) R23, with an appendix by
D.J.~Broadhurst.
\bibitem{Kassel}
C.~Kassel, {\em Quantum Groups}, Springer 1995.
\bibitem{Jarvis}
I.~Tsohantjis, A.C.~Kalloniatis, P.D.~Jarvis, {\em
Chord diagrams and BPHZ subtractions} UTAS-PHYS-96-03,
(hep-th/9604191).
\bibitem{vass}
Dror Bar-Natan, Topology 34,2 (1995) 423.
\bibitem{BK15}
D.J.~Broadhurst, D.~Kreimer, 
Phys.Lett.{\bf B393} (1997) 403;\\
D.J.~Broadhurst, {\em Conjectured Enumeration of irreducible Multiple
Zeta Values, from Knots and Feynman Diagrams}, hep-th/9612012.
\bibitem{4TR}
D.~Kreimer, {\em Weight Systems from Feynman Diagrams},
to appear in JKTR, MZ-TH/96-36, (hep-th/9612010);\\
D.J.~Broadhurst, D.~Kreimer,
{\em Feynman Diagrams as a Weight System: Four-Loop Test of a Four-Term
Relation}, MZ-TH/9612011, (hep-th/9612011), subm.~to PLB.
\end{thebibliography}
\end{document}